\documentclass[alpha-refs]{wiley-article}

\usepackage{graphicx}
\usepackage[space]{grffile}
\usepackage{latexsym}
\usepackage{textcomp}
\usepackage{longtable}
\usepackage{tabulary}
\usepackage{booktabs,array,multirow}
\usepackage{amsfonts,amsmath,amssymb, amstext, amsxtra, amscd }
\usepackage{natbib}
\usepackage{placeins}
\usepackage{url}
\usepackage{hyperref}
\hypersetup{colorlinks=false,pdfborder={0 0 0}}
\usepackage{etoolbox}
\newif\iflatexml\latexmlfalse

\AtBeginDocument{\DeclareGraphicsExtensions{.pdf,.PDF,.eps,.EPS,.png,.PNG,.tif,.TIF,.jpg,.JPG,.jpeg,.JPEG}}

\usepackage[utf8]{inputenc}
\usepackage[english]{babel}

\usepackage{siunitx}
\usepackage{epsfig}
\usepackage{comment}

\def\rate{r_{\epsilon}}
\def\windex{\upsilon}
\def\Windex{\Upsilon}
\def\trace{\mathrm{trace}}
\def\diag{\mathrm{diag}}
\def\supp{\mathrm{supp}}
\def\range{\mathrm{range}}
\def\Cov{\mathrm{Cov}}
\def\Var{\mathrm{Var}}

\newtheorem{assumption}{Assumption}

\iflatexml

\else

\paperfield{Bayesian inverse problems}
\corraddress{University of Edinburgh, Kings Buildings, Peter Guthrie Tait Rd, Edinburgh EH9 3FD, United Kingdom}
\corremail{N.Bochkina@ed.ac.uk}
\fundinginfo{Jenovah Rodrigues was funded for his PhD by the MIGSAA Doctoral Training Centre, EPSRC grant EP/L016508/01. }
\fi

\papertype{Original Article}

\title{Bayesian inverse problems with heterogeneous variance}

\author{Natalia Bochkina}
\author{Jenovah Rodrigues}

\affil{University Of Edinburgh and Maxwell Institute, UK}

\runningauthor{Bochkina and Rodrigues}

\begin{document}

\maketitle
\selectlanguage{english}

\section{Introduction}

{\label{707961}}\par\null

\subsection{Bayesian linear inverse problem}
Consider the following probability model for $Y$ which are noisy indirect observations of an unknown function $\mu$:
\begin{equation}
Y= K\mu +\epsilon  W,  \label{f2model}
\end{equation}
where $\mu \in H_0$, a separable Hilbert space and a known, injective, continuous, linear operator $K$ maps $\mu$ into another separable Hilbert space, $H_1$. Here a scalar $\epsilon$ represents the level of noise, and $W$ is a random Gaussian process in $H_1$. See \citet{knapik_bayesian_2011} for a statistical interpretation of this model. We consider the setting where $W$ is not necessarily a white noise. This occurs for differential operators \cite{agapiou_posterior_2013} and econometric problems \cite{FlorensSimoni}. Also, model (\ref{f2model}) can be viewed as a continuous experiment used to study asymptotically equivalent discrete models \cite{JohnstoneSilverman1997, JSchmidtHieber}.


We consider ill-posed problems when the solution, even of a noise free problem (with $\epsilon=0$), does not depend continuously on observations. This happens, for instance when  the eigenvalues of operator $K^T K$ decay to 0 \cite{cavalier_nonparametric_2008}.
 Typically, most methods for solving linear ill-posed inverse problems involve \emph{regularising} the solution space, by constraining the set of solutions using some a priori information such as a small norm, sparsity or smoothness, normally leading to a unique solution in a noise free case. For further details, see \citet{engl_regularization_1996}. Most regularised solutions can be interpreted as a Bayesian estimator where the regularisation is reflected as the prior information. For a more detailed discussion of the correspondence between the penalised likelihood and Bayesian approaches, see \citet{bochkina_consistency_2013}.

However, the Bayesian perspective brings more than merely a
different characterisation of a familiar numerical solution.
Formulating  a statistical inverse problem as one of inference in
a Bayesian model has great appeal, notably for what this brings in
terms of coherence, the interpretability of regularisation
penalties, the integration of all uncertainties, and the
principled way in which the set-up can be elaborated to encompass
broader features of the context, such as measurement error,
indirect observation, etc. The Bayesian formulation comes close to
the way that most scientists intuitively regard the inferential
task, and in principle allows the free use of subject knowledge in
probabilistic model building \cite[etc]{kaipio_inverse_1999, rover_coherent_2007, voutilainen_model_2009, cotter_bayesian_2009, auranen_bayesian_2005}. For an interesting philosophical view on inverse problems, falsification, and the role of Bayesian
argument, see \citet{tarantola_popper_2006}. Various Bayesian methods to solve practical inverse problems have been proposed by \citet{efendiev_bayesian_2010, cotter_bayesian_2009, dashti_besov_2012}, among many others.

\subsection{Posterior consistency and contraction rate}

The solution to an inverse problem in the presence of noise is usually analysed by taking the limit of the noise $\epsilon \to 0$. In a Bayesian approach, the solution is a probability measure $\mathbb{P}(\mu \in \cdot \mid Y)$ over a set of functions $\mu$ which depends on observations $Y$, making it a random probability measure over a set of functions. \citet{ghosal_convergence_2000} proposed to study  {\it contraction rate} of the posterior distribution which is defined as the smallest $\rate= \rate(\epsilon)$ such that for every $M \rightarrow \infty$,
\begin{equation}\label{eq:defCOntractionRate}
  \mathbb{P}(\mu: \, ||\mu-\mu_0|| \geq M \rate|Y) \stackrel{\mathbb{P}_{\mu_0}}{\to} 0,  \text{~as~~} \epsilon \rightarrow 0,
\end{equation}
uniformly over the true solution $\mu_0$ in a relevant functional class (e.g. a Sobolev class). Here $\mathbb{P}_{\mu_0}$ is the true distribution of $Y$ data under model (\ref{f2model}) with $\mu=\mu_0$.

In case of the inverse problem under the white noise model, the non-adaptive rate of contraction of the posterior distribution in a sequence space with Gaussian prior was studied by \citet{knapik_bayesian_2011}. More generally, \citet{ray_bayesian_2013} studied general adaptive priors, including wavelets, with suboptimal rates.  When the covariance operator of the noise is not identity, this problem was studied by \citet{agapiou_posterior_2013} and \citet{FlorensSimoni}, with the particular types of covariance operators motivated by the respective areas of application.  Motivated by inverse problems arising in PDEs, \citet{agapiou_posterior_2013} considered a case where the covariance operators do not necessarily commute. In such a challenging setting, to obtain conditions for the contraction of the posterior distribution, the authors assumed the unknown function to be continuous, and expressed  the smoothness of the unknown function in terms of the prior covariance operator; their contraction rates were slower than the minimax optimal ones in the case of white noise \cite{knapik_bayesian_2011}, although in many cases the exponent in the rate could be arbitrarily close to the optimal exponent. \citet{FlorensSimoni} investigated contraction of the posterior distribution in a challenging case of non-trivial covariance operators motivated by inverse problems arising in econometrics; to overcome the challenges, the authors assumed the covariance operator of the noise is trace class and true functions have monotonically decreasing coefficients in some basis (resulting in a subclass of Sobolev spaces) where they showed that the posterior contracts at the minimax rate, up to a log factor.

 In the adaptive setting, where the posterior distribution adapts to the unknown smoothness of $\mu_0$, a direct problem in a sequence space with non-decreasing known variances (which is equivalent to an ill-posed inverse problem with white noise) was considered by \citet{BelitserCredible2017},  under a restricted bias condition using a sieve prior. \citet{AgapiouMathe2017} proposed to choose a data-dependent prior mean as a way of making their posterior distribution  adaptive and to contract at the minimax optimal rate; such an approach is common in optimisation but in its current form may be less intuitive to a Bayesian statistician. Other ways to achieve adaptation are to use an empirical Bayes and a full Bayesian estimator of, either a prior hyperparameter under a white noise error model with a self-similarity assumption \cite{knapik_bayes_2016},    the truncation number under the white noise model \cite{JohannesAdaptive2020}, or a plug-in estimator of the prior scale parameter for direct models under white noise and H\"older spaces  \cite{SzaboEmpiricalBayes}, that lead to the posterior contraction rate that achieves the minimax rate, up to a constant.  \citet{FlorensSimoni} proposed an estimator of the prior scale parameter in inverse problems with heterogeneous noise, which resulted in a suboptimal rate of contraction. In a different setting, \citet{KnapikContraction2018} showed that adaptive hierarchical mixture models lead to optimal posterior contraction rates. We emphasize that for inverse problems, only the approach of \citet{JohannesAdaptive2020} leads to the posterior distribution contracting at the optimal rate (under white noise) without a restriction to a subset of functions.

\subsection{Our contribution}

In this paper, we focus on  linear inverse problems with Gaussian noise where it is possible to reduce model (\ref{f2model}) to the sequence space using a Riesz basis of $H_0$ so that its image under the forward operator is a basis of $H_1$, for which there exists a biorthogonal basis. This generalises a standard assumption of reduction to a sequence space where the covariance operator $V$ of the Gaussian noise $W$ and operator $K K^T$ are simultaneously diagonalisable. We consider covariance operators that do not have to belong to the trace class (e.g. white noise or the generalised derivative of fractional Brownian motion, called fractional noise), nor do we constrain our unknown function of interest to be continuous or to have monotonically decreasing coefficients in some basis. Motivation for this approach comes from wavelet-based approaches to inverse problems \cite{AbramovichSilverman}, that apply to homogeneous operators $K$, and representation of fractional noise in terms of biorthogonal wavelet bases \cite{MeyerSellanTaqqu}, leading to a novel approach we refer to as vaguelette-vaguelette  that we then generalise to non-wavelet bases. Our first contribution is the derivation of the posterior contraction rate for mildly ill-posed problems under fractional noise using the vaguelette-vaguelette approach, and identifying nonadaptive priors leading to optimal posterior contraction rate, in the minimax sense, for the true function belonging to a Sobolev space. The results are derived in a more general setting.

Our second contribution is to investigate the case where the variances in the sequence space may be unknown,  and how using their plug-in estimator affects the posterior contraction rate. We illustrate this on an example with repeated observations. We also study the problem of error in operator and identify conditions when the rate of posterior contraction of $\mu$ is not affected.

Our third contribution is the study of the empirical Bayes posterior with a plugged-in maximum marginal likelihood estimator of the prior scale under an appropriate Gaussian prior, in the sequence space.  We show that  the corresponding empirical Bayes posterior distribution contracts at the optimal rate (in the minimax sense), adaptively over a set of Sobolev classes, under some conditions on prior smoothness. Also, we show that under an additional self-similarity condition, the  maximum marginal likelihood estimator of the prior scale concentrates around its oracle value.

The paper is organised as follows. In Section~\ref{sec:BayesianModel} we introduce the Bayesian formulation of an inverse problem, fractional Brownian motion and its generalised derivative, mildly ill-posed inverse problems and problems with error in operator. In  Section~\ref{sec:LikInvProbfBMWavelets} we review a wavelet-based formulation of inverse problems. In  Section~\ref{sec:RatesWaveletFBM} we consider mildly ill-posed inverse problems under fractional noise, formulate the vaguelette-vaguelette approach and study the posterior contractions rates for the problems discussed in Section~\ref{sec:BayesianModel}. We also study properties of an empirical Bayes posterior distribution for this problem that adapts to unknown smoothness of $\mu_0$.  In Section~\ref{sec:GeneralRates} we formulate this problem for more general bases that lead to a sequence space representation of (\ref{f2model}) with arbitrary sequences of coefficients of ill-posedness, prior and noise variances, and study the posterior contractions rates for the problems discussed in Section~\ref{sec:BayesianModel}. In Section~\ref{sec:PolyEV} these general results are applied to the case where all sequences (of coefficients of ill-posedness, prior and noise variances) decay polynomially, which includes mildly ill-posed problems with fractional noise. We also study properties of an adaptive empirical Bayes posterior distribution in this setting. The effects of the choice of the prior parameters and sample size on the posterior contraction rate is illustrated on simulated (synthetic) data, with both fixed and estimated prior scale, for mildly ill-posed inverse problems (Section~\ref{sec:SimData}). We conclude with a discussion. The proofs are given in the appendix. 

Notation:  $a_n \asymp b_n$ denotes the existence of $c, C$: \, $0 < c \leq a_n/b_n \leq C < \infty$ for all $n \geq 1$. Here $||\cdot||$ denotes the Euclidean vector norm if the argument is a vector, and the $L^2$ norm if the argument is a function. For an operator $K: \, L^2 \to L^2$, where $L^2$ is viewed as a Hilbert space, define operator $K^T$ as follows: for any $f, g \in L^2$, $\langle K^T f, g \rangle= \langle  f, K g \rangle$.





\section{Bayesian inverse problem with heterogeneous Gaussian noise}\label{sec:BayesianModel}

\subsection{Bayesian inverse problem}

Note that model (\ref{f2model}) in Hilbert space $L^2[0,1]$ can be viewed as the following model
 \begin{equation}
d \tilde Y_t= (K\mu)(t) dt +\epsilon  d \tilde W_t, \,\, t\in [0,1],  \label{f2modelDiff}
\end{equation}
with $Y = (d \tilde Y_t/dt, \, t \in [0,1])$ and $W = (d \tilde W_t/dt, \, t \in [0,1])$ where the derivatives are generalised, i.e. in both cases the models can only be used to compute scalar products with appropriate test functions \cite{knapik_bayesian_2011}. It has been shown that model (\ref{f2modelDiff}) is  asymptotically equivalent, in Le Cam distance, to a finite dimensional model with growing sample size $n$; if the errors  are independent then such a finite dimensional model is equivalent to (\ref{f2modelDiff}) with $\tilde W_t$ being Brownian motion and $\epsilon = n^{-1/2}$ (the white noise model), while if the errors are correlated it is equivalent to (\ref{f2modelDiff}) with $\tilde W_t$ being fractional Brownian motion with parameter $H$  and $\epsilon = n^{-(1-H)}$ for $H>1/4$ under an assumption on correlation errors, for details on the direct problem with $K=I$ see \citet{JSchmidtHieber}.

We assume that $W_t$ is a Gaussian process with covariance operator $V: \, H_1 \to H_1$ such that for any bounded function $f, g \in L^2[0,1]$,  $ \langle V f, g\rangle = E \left[\int f(s) g(t)  W_s   W_t ds dt\right]$. If $W_t$ is a derivative of a Brownian motion, i.e. (\ref{f2model}) is a white noise model, then $V = I$. Note that $V$, just like the identity element of $H_1$, does not have to be of trace class.

We consider the prior distribution of $\mu$ to be a centred Gaussian process 
\begin{equation}\label{eq:priorOp}
\mu \sim N(0,\Lambda),
\end{equation}
 with covariance operator $\Lambda: \, H_0 \to H_0$ belonging to a trace class ($\trace(\Lambda) <\infty$). Here we also consider $H_0 = L^2[0,1]$.

%
Then, the posterior distribution of $\mu$ is given by
$$
\mu \mid Y \sim N\left(\left( K^T  V^{-1} K +\epsilon^{2}\Lambda^{-1}\right)^{-1} K^T V^{-1} Y,   \left(\epsilon^{-2} K^T  V^{-1} K +\Lambda^{-1}\right)^{-1}\right).
$$
 Similar results were presented in \citet{agapiou_posterior_2013} and  \citet{FlorensSimoni}.
   The posterior distribution is proper if $\trace\left( \left(\epsilon^{-2} K^T  V^{-1} K +\Lambda^{-1}\right)^{-1}\right) <\infty$.

We assume that the true distribution of $Y_t$, $t\in [0,1]$ is given by model (\ref{f2model}) with $\mu=\mu_0$ that we denote by $P_{\mu_0}$; sometimes we write $P_{\mu_0,V}$ if we need to be specific about the covariance operator of the noise. The expectation with respect to this distribution is denoted by $\mathbb{E}_{\mu_0}$.

Our aim is to determine the contraction rate of the posterior distribution of $\mu$ given $Y$ around the true value $\mu_0$ as the error $\epsilon \to 0$.  Next we introduce our motivating example, fractional Brownian motion noise, and mildly ill-posed inverse problems.

\subsection{Fractional Brownian motion and fractional noise}\label{sec:DefBrownianMotion}

In this section we focus on a particular case where $W_t$ in (\ref{f2model}) is a generalised derivative of fractional Brownian motion (fBM), which we call fractional noise, defined below. 

\begin{definition}
Let $B = (B_t)_{t\in[0,1]}$ be a centered Gaussian process. $B$ is called a fractional Brownian motion (fBm) with Hurst exponent $H\in (0,1)$ if it has the following covariance structure:
$$
 \mathbb{E} (B_s B_t) = \frac 1 2 \left( |t|^{2H}  + |s|^{2H} - |t- s|^{2H} \right).
$$
\end{definition}
fBm is a self-similar process with stationary increments. When $H=1/2$, it coincides with the standard Brownian motion. Sample paths of fBM are H\"older continuous of any order strictly less than $H$. See  \citet{picard} for further details.


Series representations of fBM in various orthonormal bases and spanning frames of $L^2[0,1]$ are available. These bases can be normalised harmonics $\sin(x_{k,H} t)$ and $1-\cos(y_{k,H} t)$ with specific frequencies $y_{k,H}$ and $x_{k,H}$  \citet{Dzhaparidze} also extended it to the case of multidimensional support, and more involved bases in \citet{JSchmidtHieber} and \citet{KuhnLinde2002}. A general series representation is given in \citet{GILSINGSOTTINEN}, as well as other bases such as shifted Legendre polynomials and modified hypergeometric functions. A good introduction to fractional Brownian motion is given by \citet{picard}.
Efficient simulation of fractional Brownian motion is discussed in \citet{Ndaoud2016}, who gives a series representation of a fBM with Hurst exponent $H$ in terms of a trigonometric series, which implies that its generalised derivative (fractional noise) can be decomposed in a series with functions $(1, \cos(\pi k  t), \sin(\pi k  t), \, k \in \mathbb{N})$.

\citet{MeyerSellanTaqqu} showed that one-dimensional fractional Brownian motion with Hurst exponent $H \geq 1/4$ can be written as a wavelet series which we discuss in detail in Section~\ref{sec:LikInvProbfBMWavelets}, alongside a  wavelet-based representation of inverse problem (\ref{f2model}).

\subsection{Mildly ill-posed inverse problems}\label{sec:MildlyIP}

In this paper we focus on mildly ill-posed inverse problems. An inverse problem (\ref{f2model}) is called mildly ill-posed if the eigenvalues $\{k_i^2\}$ of operator $K^T K$, arranged in the decreasing order, decrease polynomially.
\begin{assumption}\label{AssumeEVIP}
There exist $p \geq 0, C_1 \geq 1$ such that the eigenvalues $(k_i^2)$  satisfy $C_1^{-1}i^{-p} \leq k_i \leq C_1 i^{-p}$ for all $i=1,2,\ldots$.
\end{assumption}
Examples of such inverse problems are, for instance, Volterra operator (Section~\ref{sec:SimData}) and Radon operator used in tomography. They also arise in econometrics \cite{FlorensSimoni}.

\subsection{Error in operator}\label{sec:ErrorOperatorIntro}

In this paper we also study the problem when operator $K$ is observed with error, in particular how this affects the posterior contraction rate of $\mu$. \citet{HoffmannReiss} studied the inverse problem (\ref{f2model}) with white noise (i.e. $V=I$) where the operator is observed with errors, namely
 \begin{eqnarray}\label{eq:ErrorOperator}
\hat K = K + \delta Z,
\end{eqnarray}
where $Z$ is the operator white noise, independent of observational noise $W$ and $\delta$ is the observation error. We will extend this approach to a model with dependent errors (Section~\ref{sec:ErrorOperatorGeneral}), and in particular to the case of fractional noise (Section~\ref{sec:ErrorInOperatorFBM}).



\section{Inverse problems and fractional noise via wavelets}\label{sec:LikInvProbfBMWavelets}


\citet{Donoho1995} proposed two wavelet-based approaches to linear inverse problems with operator $K$ based on an orthonormal wavelet basis $(\psi_{jk})$, called vaguelette-wavelet and wavelet-vaguelette which we introduce below.

\subsection{Wavelets and Riesz basis}

A one-dimensional, orthonormal, compactly-supported wavelet basis is given by $\{\phi,\psi_{j,k}, j=0,1,\ldots, \, k=0,1,\ldots, 2^j-1\}$  with $\psi_{j,k}(x) = 2^{j/2} \psi(2^j x-k)$, where $\phi$ is the scaling function (father wavelet) and $\psi$ is the wavelet function (mother wavelet), both having a compact support. To simplify the notation, it is common to denote $\psi_{-1,0}= \phi$. Wavelets are said to be of regularity $r$ if they have $r$ continuous derivatives and vanishing moments up to order $r$, i.e. $\int x^k \psi(x) dx =0$ for all $k=1,2,\ldots, r$. We denote an orthonormal wavelet basis on $[0,1]^d$ for a fixed $d$ by $(\psi_{\windex}, \windex \in \Windex)$, for a countable $\Windex$. In the one-dimensional setting, with the wavelet and the scaling functions supported on $[0,1]$, $\windex=(j,k)$ and $\Windex = \{(j,k): \, j=0,1,2,\ldots, \, k=0,1,\ldots, 2^j-1\}$.  See  \citet{Vidakovic1999} for the compactly supported wavelets and their statistical applications, including higher dimensional wavelets.

A wavelet approach to inverse problems involves generalisation of an orthonormal basis to a Riesz basis that we define below.
\begin{definition}\label{def:RieszBasis}
Sequence of functions $(v_{\windex})$ forms a Riesz basis in $L^2[0,1]$ if $(v_{\windex})$ spans  $L^2[0,1]$ and
  there exist $0 < A \leq B <\infty$ such that for any square summable sequence $(c_{\windex})$,
\begin{equation}\label{eq:RieszBasis}
A \sum_{\windex \in \Windex} c_{\windex}^2 \leq || \sum_{\windex \in \Windex} c_{\windex} v_{\windex}||^2 \leq B \sum_{\windex \in \Windex} c_{\windex}^2.
\end{equation}
\end{definition}
Typically a Riesz basis is an image of an orthonormal basis under some bounded invertible operator. \citet{Donoho1995}  showed that this also applies to what he called weakly invertible operators for a given orthonormal basis. An orthonormal basis is trivially a Riesz basis with $A=B=1$ due to Parseval's identity.

We will also rely on the concept of biorthogonal bases which provides an infinite dimensional generalisation of the singular value decomposition of matrices.
\begin{definition}\label{def:BiorthogonalBasis}
Given a set of functions $(v_{\windex})_{\windex \in \Windex}$ spanning $L^2[0,1]$,  a biorthogonal basis $(w_{\windex})_{\windex \in \Windex}$ to $(v_{\windex})_{\windex \in \Windex}$ satisfies
\begin{equation}\label{eq:DefBiorthogonal}
\langle v_{\windex}, w_{\windex}\rangle =1 \text{ for all } \windex \in \Windex  \text{ and } \langle v_{\windex_1}, w_{\windex_2}\rangle =0 \text{ for all } \windex_1 \neq \windex_2, \,\,\,  \windex_1, \windex_2 \in \Windex.
\end{equation}
\end{definition}
For more details on biorthogonal bases, including construction, see  \citet{Vidakovic1999}.

\subsection{Sobolev classes in terms of wavelets}\label{sec:SobolevWavelets}

Note that the Sobolev norm can be written equivalently in terms of wavelet coefficients  for wavelets with regularity $r > \beta$,
\begin{equation}\label{eq:DefSobolevWavelet}
 S^{\beta}(A)=\left\{f(x) = \sum_k a_k \phi(x-k)+ \sum_{j=0}^{\infty} \sum_{k} b_{jk} \psi_{jk}: \, || a ||+  \left[\sum_{j=0}^{\infty}  2^{2 \beta j} ||b_{j, \cdot}||^2 \right]^{1/2}\leq A \right\},
\end{equation}
\cite{JohnstoneSilverman1997}. 
\citet{CohenDeVoreBiorthogonalBesov} show the equivalence of the Sobolev norm and the above sequence norm for biorthogonal compactly supported  wavelet bases with regularity $r > \beta$. 


\subsection{Wavelet representation of fractional noise}\label{sec:WaveletFBM}


\citet{MeyerSellanTaqqu} proved that a fractional derivative, or antiderivative $Z_{d+1/2}(t) = D^{-d} Z(t)$, of white noise  $Z(t)$ has the same distribution as the generalised derivative of the fBM with the Hurst exponent $H=d+1/2$ for $d<1/2$, and the same distribution as the following wavelet random series:
\begin{equation}\label{eq:WaveletFracNoise}
W_t =\sum_{k \in D} \tilde z_k \phi_{\Delta}^{(d)}(t-k) + \sum_{\windex \in \Windex} \sigma_{\windex} z_{\windex} \psi_{\windex}^{(d)}(t),
\end{equation}
with  iid $z_{\windex} \sim N(0,1)$, $ \sigma_{\windex}=C 2^{-j(H-1/2)}$, $\windex = (j,k)$.  Here $\tilde z_k = \sum_{\ell=0}^{\infty} \gamma_\ell^{(-d)} \zeta_{k-\ell}$ is a fractional ARIMA(0, d, 0) process bounded with high probability, where the coefficients $\gamma_{k}^{(-d)}$ are defined by
\begin{equation}\label{eq:DefGamma}
(1- x)^{-d} = \sum_{k=0}^{\infty} \gamma_k^{(-d)} x^k,\quad  \text{ for } |x|<1,
\end{equation}
with property $|\gamma_k^{(-d)}|=O(k^{-1-d})$ for large $k$.
 Here the scaling function $\phi_{\Delta}^{(d)}$ and the wavelet function $ \psi^{(d)}$ are constructed from wavelet and scaling functions $\psi$ and $\phi$
   in terms of their Fourier transform: $\widehat{\phi_{\Delta}^{(d)}}(\xi) = \left(\frac{e^{i  \xi}}{1-e^{i \xi}} \right)^d \widehat{\phi}(\xi)$ and $\widehat{\psi^{(d)}}(\xi) =  e^{i  \xi d}  \widehat{\phi}(\xi)$ where $\widehat{f}(\xi) = \int e^{-i t \xi} f(t) dt$, so that all of their integer moments vanish (see \citet{MeyerSellanTaqqu} eq (3.1) for the full set of conditions and Proposition 3 for the statement).
   This gives a Riesz wavelet basis $(\phi^{(d)}_{\Delta}(\cdot -m), \psi^{(d)}_{\windex})$, with $\psi^{(d)}_{\windex}$ being dilations and translations of a single function $\psi^{(d)}$, which is biorthogonal to the wavelet Riesz basis $(\phi^{(-d)}_{\Delta}(\cdot -m), \psi^{(-d)}_{\windex})$.  \citet{AyacheTaqqu2003} showed this  under weaker conditions on wavelets $\psi$, including compactly supported Daubechies wavelets of regularity $r \geq 6$, as well as alternative wavelet representations of fractional noise (precise conditions on wavelets $\psi$ are given by the authors  on p. 456). 
 \citet{AbrySellan} constructed the corresponding filters that allow fast practical implementation of wavelet transform and its reverse using the cascade algorithm. To simplify the notation, we denote $\Windex = \Windex_{\phi} \cup  \Windex_{\psi}$ with 
  $\windex  \in \Windex_{\phi} = \{(-1,k), \, k\in D\}$, and
 $\psi_{\windex}^{(d)} = \phi_{\Delta}^{(d)}(\cdot-k)$ for $\windex =  (-1,k) \in \Windex_{\phi}$. For compactly supported wavelets, $D$ and $\Windex_{\phi}$ are finite, and $\Windex_{\psi} = \{(j,k): \, j=0,1,\ldots, k \in K_j, \, |K_j|\leq C_{\psi} 2^j \}$.










\subsection{A wavelet-vaguelette approach}\label{sec:wavevague}

Wavelets transformed by some fixed operator $A$  are called vaguelettes. In this section we describe the wavelet-vaguelette approach, and in the next we discuss the vaguelette-wavelet approach.

A wavelet-vaguelette approach applies to operators $K$ such that wavelets $\psi_{\windex}$ are in their image, and involves the numerical evaluation of $K^{-1} \psi_{\windex}$ \cite{Donoho1995}. Then, we can write $(K \mu_0)(t) = \sum_{\windex \in \Windex} c_{\windex} \psi_{\windex}(t)$ and hence
$$
\mu_0(t) = \sum_{\windex \in \Windex} c_{\windex} \tilde\beta_{\windex}^{-1} \left(K^{-1}\psi_{\windex}\right)(t)\tilde\beta_{\windex}= \sum_{\windex \in \Windex} \mu_{0;\windex}  \tilde v_{\windex}(t),
$$
with $\mu_{0; \windex}=c_{\windex} /\tilde\beta_{\windex}$ and normalised vaguelettes $\tilde v_{\windex}(t)=\left(K^{-1}\psi_{\windex}\right)(t)\tilde\beta_{\windex}$. Here  $\tilde\beta_{\windex}$ are  normalisation factors so that $||\tilde v_{\windex}||=1$.  This approach is applicable if the normalised vaguelettes $\tilde v_{\windex}(t)$ span $L^2[0,1]$ and form a Riesz basis.

If the noise $W$ can be written as
\begin{equation}\label{eq:GaussPWaveletSeries}
W_t = \sum_{\windex \in \Windex} \sigma_{\windex} z_{\windex} \psi_{\windex}(t),
\end{equation}
for some coefficients $\sigma_{\windex}$ and iid $z_{\windex} \sim N(0,1)$, then the model for the wavelet transform of the observations $y_{\windex} = \langle \psi_{\windex},  Y \rangle$, under (\ref{f2model}), is given by
 $y_{\windex} \mid \mu_{\windex} \sim N(\mu_{\windex}\tilde\beta_{\windex}, \epsilon^2 \sigma_{\windex}^2)$, independently for $\windex \in \Windex$.

Therefore, under the above assumptions, using the wavelet-vaguelette approach reduces problem (\ref{f2model}) to a sequence space model.  For fractional noise, as far as we are aware, representation (\ref{eq:GaussPWaveletSeries}) does not hold.


\subsection{A vaguelette-wavelet  approach}\label{sec:vaguewave}

A vaguelette-wavelet approach works if $\psi$ is in the image of $K^T$, and the normalised image of these functions under operator $K$ forms an orthonormal Riesz basis in $L^2[0,1]$ \cite{Donoho1995}. Formally these assumptions are formulated as follows.
\begin{definition}\label{def:AssumptionVaguelette}
Given a linear operator $K$ and a wavelet basis $(\psi_{\windex})_{\windex \in \Windex}$ assume that there exist scalars $\beta_{\windex}$, $\windex \in \Windex$ such that the functions 
$v_{\windex}= K\psi_{\windex}/\beta_{\windex}$
   form a Riesz basis, 
 and  there exists a biorthogonal basis $(w_{\windex})_{\windex \in \Windex}$ for $(v_{\windex})_{\windex \in \Windex}$,  defined by (\ref{eq:RieszBasis}) and (\ref{eq:DefBiorthogonal}) respectively.
  Such functions $(v_{\windex})$ are called vaguelettes.
\end{definition}
 \citet{AbramovichSilverman} give examples of linear operators that satisfy this condition, such as homogeneous operators $K$, where for all $t_0$ and $s>0$, $K[ f(s(t-t_0)) ] = s^{-p} (K f)(s(t-t_0)) ]$ for some constant $p$, called an index of an operator. Examples of homogeneous operators include integration, differentiation (including fractional integration and differentiation) and other operators such as Radon transform in an appropriate coordinate system \cite{Kolaczyk1996}. For such operators, the normalising constants are  $\beta_{jk}= C_{\beta} 2^{-p j}$, and the corresponding vaguelettes are dilations and translations of a single function but they are not necessarily mutually orthogonal. This property also holds for some convolution operators \cite{Donoho1995}.

In this approach, $\mu =  \sum_{\windex \in \Windex} \mu_{\windex} \psi_{\windex}$, and hence $K\mu =  \sum_{\windex \in \Windex} \mu_{\windex} \beta_{\windex} v_{\windex}$. If $(w_{\windex})$ is a biorthogonal basis for $(v_{\windex})$, i.e. $\int v_{\windex} w_{\windex} =0$ for all $\windex \in \Windex$, then the sequence space model can be written as
$$
y_{\windex} \sim N(\mu_{\windex}\beta_{\windex}, \epsilon^2 \sigma_{\windex}^2), \quad \windex \in \Windex,
$$
  where  $y_{\windex}  = \int Y w_{\windex}$, and $\sigma_{\windex}^2$ is such that $ \epsilon^2 \sigma_{\windex}^2 = \Var(y_{\windex}) = \int V(s,t)  w_{\windex}(t)  w_{\windex}(s) ds dt $.  $y_{\windex}$'s  are independent if $\Cov(y_{\windex}, y_{\windex'}) = \int V(s,t)  w_{\windex}(t)  w_{\windex'}(s) ds dt =0$ for all $\windex \neq \windex'$. Note that for fractional noise this approach applies with fractional wavelets which are not necessarily an image of an orthogonal wavelet basis under $K$ (Section~\ref{sec:WaveletFBM}).


  This also implies that  the normalising coefficients $\beta_{\windex}$ represent the level of ill-posedness of operator $K$. Therefore, this approach also reduces model (\ref{f2model}) to a sequence space model. However, none of these approaches apply to a fractional noise.

In the next section we extend these approaches to use representation of fractional noise in terms of vaguelettes $\psi^{(d)}$ (Section~\ref{sec:WaveletFBM}) and give results on the posterior contraction rates for mildly ill-posed inverse problems, under fractional noise, using a wavelet approach with $\mu_0$  in a Sobolev class, which in turn will rely on general results stated in Section~\ref{sec:GeneralRates}.






\section{Posterior contraction rate under fractional noise}\label{sec:RatesWaveletFBM}


In this section we study the posterior contraction rate of mildly ill-posed inverse problems introduced in Section~\ref{sec:MildlyIP} with fractional noise decomposed in a wavelet series (Section~\ref{sec:WaveletFBM}), including the cases of error in operator, unknown Hurst exponent $H$ and adaptation to  unknown smoothness of $\mu_0$. First, we introduce a new approach, called vaguelette-vaguelette, to account for the vaguelette type of wavelet used in the decomposition of a fractional noise.

\subsection{Vaguelette-vaguelette approach}\label{sec:vaguevague}

Following the series decomposition of fractional noise in terms of biorthogonal wavelets $\left(\psi^{(d)}_{\windex}\right)$
and   $\left(\psi^{(-d)}_{\windex}\right)$ with $d=H-1/2$ (Section~\ref{sec:WaveletFBM}), we need to adapt the wavelet-vaguelette approach to the case of wavelet $\psi$ being replaced by $\psi^{(d)}$, which is a particular kind of vaguelette determined by the covariance operator of noise. We call this a vaguelette-vaguelette approach.
 Similarly to the wavelet-vaguelette approach, the vaguelette-vaguelette approach applies to operators $K$ such that wavelets $\psi^{(d)}_{\windex}$ are in their image, which holds for example if   $\range(K) \subseteq \range(D^{-d})$ where $D^{-d}$ is a differentiation/integration operator (Section~\ref{sec:WaveletFBM}).

Then, we can write $(K \mu)(t) = \sum_{\windex \in \Windex} c_{\windex} \psi^{(d)}_{\windex}(t)$ with $c_{\windex} = \langle \psi^{(-d)}_{\windex}, K\mu\rangle$ and hence
$$
\mu(t) = \sum_{\windex \in \Windex} c_{\windex} k_{\windex}^{-1} \left(K^{-1}\psi^{(d)}_{\windex}\right)(t)k_{\windex}= \sum_{\windex \in \Windex} \mu_{\windex}  \tilde v_{\windex}^{(d)}(t),
$$
with $\mu_{ \windex}=c_{\windex} /k_{\windex}$ and normalised vaguelettes
$$
\tilde v^{(d)}_{\windex}/k_{\windex} =K^{-1}\psi^{(d)}_{\windex}  = K^{-1} D^{-d} \psi_{\windex},
$$
with normalisation factors $k_{\windex}$ so that $\left\lVert\tilde v^{(d)}_{\windex}\right\rVert=1$.


This approach is applicable if the normalised vaguelettes $\left( \tilde v^{(d)}_{\windex}\right)$ span $L^2[0,1]$ and form a Riesz basis. Under the former assumption, any $\mu(t) \in L^2[0,1]$ can be decomposed into this vaguelette basis, and the latter assumption allows us to study the posterior contraction rate of $\mu$ in terms of its coefficients $\mu_\windex$.

Following the wavelet-vaguelette and vaguelette-wavelet approaches, we call this inverse problem mildly ill-posed if for some $p>0$
\begin{equation}\label{eq:MildIPwavelets}
\kappa_{\windex} =   \left\lVert K^{-1}\phi_{\Delta, \windex}^{(d)} \right\rVert \asymp 1, \, \windex \in \Windex_{\phi}
, \text{ and }
\kappa_{\windex} = \left\lVert K^{-1}\psi^{(d)}_{\windex}\right\rVert \asymp 2^{-j p}, \, \,  \windex =(j,k)\in \Windex_{\psi}.
\end{equation}
Recall that this holds for homogeneous operators $K$ (Section~\ref{sec:vaguewave}).

Then, using the wavelet representation of a fractional noise (Section~\ref{sec:WaveletFBM}),  the sequence space model for the wavelet transform of observations $y_{\windex} = \langle \psi_{\windex},  Y \rangle$ under model (\ref{f2model}) with a fractional noise is given by
\begin{equation}\label{eq:SeqWaveletModel}
y_{\windex} \sim N\left(\mu_{\windex} \kappa_{\windex}, \epsilon^2 \sigma_{\windex}^2\right), \windex \in \Windex_{\psi} \text{ independently },
\end{equation}
with $\sigma_{\windex} \asymp 2^{-j(H-1/2)}$ for large $j$ and $\epsilon= n^{H-1}$, and independently of  $ y_{-1,k} = \langle \phi^{-d}_{\Delta,k},  Y \rangle$, $k \in D$.  Note that the coefficient $y_{-1,k}$ is nonzero only if  $(-k,-k+1) \cap \supp\left(\phi_{\Delta}^{(-d)}\right) \neq \emptyset$, which implies that $D$ is a finite set since $\supp\left(\phi_{\Delta}^{(-d)}\right)$ is finite. Hence, the scaling coefficients $y_{D}= \left(y_{-1,k}, \, k\in D \right)$ satisfy
$$
  y_{D} \sim N\left( K_{D} \mu_D, \epsilon^2 V_D\right),
$$
where  $\mu_D = \left(\mu_{-1,k}, \, k\in D \right)$, $K_{D} = \diag\left(\kappa_{-1,k}, \, k \in D \right)$ and the covariance matrix $V_D$ has variances $\sigma_{-d}^2$ on the diagonal and $V_{D; k,m} = \sigma_{-d}^2 \rho^{(-d)}_{|k-m|}$, where $\sigma_{-d}^2$ and $\rho^{(-d)}_{m}$ are defined in terms of the sequence $\gamma_\ell^{(-d)}$ (see eq \eqref{eq:DefGamma} ) by
$$
\sigma_{-d}^2  = \sum_{\ell=0}^{\infty} \left[\gamma_\ell^{(-d)}\right]^2 < \infty \quad \text{ and } \quad  \rho^{(-d)}_{m} = \sigma_{-d}^{-2} \sum_{\ell=0}^{\infty} \gamma_\ell^{(-d)} \gamma_{\ell+|m|}^{(-d)}.
$$
Note that the support of Daubechies wavelets is of length at most $2 r$ where $r$ is the regularity of the wavelets, hence, $D$ has at most $2r$ elements. Recall that $\Windex_{\phi}= \left\{(-1,k), k\in D\right\}$.

We consider a  prior distribution on $\mu$ defined in terms of random series $\mu = \sum_{\windex \in \Windex} \mu_{\windex}  \tilde v^{(d)}_{\windex}(t)$ with $\mu_{\windex} \sim N(0, \lambda_{\windex})$ independently for $\windex \in \Windex$. We take $\lambda_{\windex} = \tau 2^{-2 \alpha j}$ for the wavelet coefficients $\windex=(j,k) \in \Windex_{\psi}$, for some $\tau, \, \alpha >0$. Denote the prior variances for the scaling coefficients as a diagonal matrix $\Lambda_D = \diag(\lambda_{-1,k}, \, k\in D)$.
  In the context of wavelets, it is common to take a non-informative prior for scaling coefficients, i.e. $\lambda_{-1,k} = +\infty$, or its proper approximation with large prior variances.
  For  $\lambda_{-1,k} = +\infty$, $k \in D$, we have $\mu_{D} \mid y_{D} \sim N( K_D^{-1} y_D,  \epsilon^{2} K_D^{-1}  V_D  K_D^{-1})$, which holds approximately if $\lambda_{-1, k}$, $k \in D$ are large but finite.
  If $\lambda_{\windex} = \tau 2^{-2 \alpha j}$ for $\windex=(j,k) \in \Windex_{\psi}$ and  $\lambda_{\windex}$ is finite for $\windex \in \Windex_{\phi}$ , then a priori, with high probability, $\mu \in S^{\alpha'}$ for all $\alpha' < \alpha$.

For such a prior distribution, the corresponding  posterior distribution is
$$
\mu_{\windex} \mid y_{\windex} \sim
N\left(\frac{ y_{\windex} \kappa_{\windex}\lambda_{\windex}}{ \lambda_{\windex} \kappa_{\windex}^2 + \epsilon^2 \sigma^2_{\windex}},  \,\,\,\, \frac{\sigma^2_{\windex}\lambda_{\windex}}{\epsilon^{-2} \lambda_{\windex} \kappa_{\windex}^2 + \sigma^2_{\windex}}\right), \text{ independently for } \windex \in \Windex_{\psi},
$$
 and for the scaling coefficients,  independently of $\mu_{\windex} \mid y_{\windex}$ for $\upsilon \in \Upsilon_\psi$, we have
\begin{eqnarray*}
\mu_{D} \mid y_{D} \sim
 N\left( \left( K_D^T  V_D^{-1} K_D +\epsilon^{2}\Lambda_D^{-1}\right)^{-1} K_D^T V_D^{-1} y_D, \,\,\,\,  \left(\epsilon^{-2} K_D^T  V_D^{-1} K_D +\Lambda_D^{-1}\right)^{-1}\right).
\end{eqnarray*}

\subsection{Minimax rate of an inverse problem with fractional noise}\label{sec:MinimaxMildIPFBM}

In this section we derive the minimax rate of convergence, in $L^2[0,1]$ norm, of estimators of $\mu$ under model (\ref{f2model}) given the true model $P_{\mu_0}$,  $W_t$ is a generalised derivative of fractional Brownian motion and the inverse problem is  mildly ill-posed, i.e. it satisfies assumption \eqref{eq:MildIPwavelets}. \citet{ghosal_convergence_2000} proposed to use it as a benchmark for posterior contraction rates. 
  The minimax rate of convergence is defined in Section~\ref{sec:MinimaxGeneral}.
 We assume that the true function $\mu_0$ belongs to a generalised Sobolev class $S^\beta(A)$   defined by (\ref{eq:DefSobolevWavelet}) with wavelets $\psi_{\windex}$ replaced by vaguelettes $\tilde v_{\windex}^{(d)}$.



 \citet{JohnstoneSilverman1997} showed that the minimax rate of estimation in $L^2$ norm in the direct problem  with fractional Brownian motion over Sobolev spaces $S^\beta$ is $n^{-2\beta (1- H)/(2\beta + 2-2 H)}$, which is a particular subset of Besov spaces studied in that paper ($S^\beta=B_{2,2}^\beta$). Now we state the minimax rate for a mildly ill-posed inverse problem with fractional noise.

\begin{proposition}\label{prop:MinimaxRatefBM}
The minimax rate of convergence in the mildly ill-posed inverse problem (\ref{f2model}) under a fractional noise with Hurst exponent $H \in (1/4,1)$, under condition (\ref{eq:MildIPwavelets}), over functions from a Sobolev class $S^\beta$  defined by (\ref{eq:DefSobolevWavelet}) with vaguelettes $\tilde v_\windex^{(d)}$, is
\begin{equation}\label{eq:RateMildInverseFracBrown}
\rate^* = \epsilon^{\frac{2\beta}{2\beta + 2p+2-2H}} = n^{-\frac{2\beta(1-H)}{2\beta + 2p+2-2H}}. 
\end{equation}
\end{proposition}
For $p=0$ this rate coincides with the rate for the direct problem \cite{JohnstoneSilverman1997}, and for $H=1/2$ it coincides with the minimax rate under white noise. This result follows from Proposition~\ref{prop:optimalrateshetero} where the minimax rate is defined in a more general setting which applies here since the
 Sobolev norm (\ref{eq:DefSobolevWavelet}) can be written equivalently  as
\begin{eqnarray}\label{def:SobolevVaguelettes}
\mu_0 \in  S^{\beta}(A)=\left\{f = \sum_{\windex \in \Windex} \tilde v_\windex f_\windex: \, \sum_{\windex \in \Windex_{\phi}}  f_\windex^2 + \sum_{i=1}^{\infty}  i^{2 \beta } f_{i}^2 \leq A^2 \right\},
\end{eqnarray}
due to the set $D$ being finite, (shown by using the bijective mapping of the set $\{ (j,k), \, j=0,1,2, \ldots, k=0,1,\ldots, 2^j-1\}$ into $\mathbb{N}$: $(j,k) \to 2^j+k$ such that $2^j \leq i=2^j+k\leq 2 \cdot 2^j$).

If a posterior contraction rate $r_{\epsilon}$ in a mildly ill-posed inverse problem, over functions from a Sobolev class with a fractional noise, matches the rate stated in Proposition~\ref{prop:MinimaxRatefBM} up to a constant, i.e. there exist finite  $C>0$ and $\epsilon_0$ independent of $\epsilon$ such that $0 < C^{-1} \leq r_{\epsilon}/\rate^* \leq C <\infty$ for all $\epsilon \leq \epsilon_0$, then we will say that the posterior distribution contracts at the optimal rate in the minimax sense.

\subsection{Posterior contraction rate for fractional noise}\label{sec:PostContractionForFBM}


Now we study the posterior contraction rate in the setting of Section~\ref{sec:vaguevague} under the vaguelette-vaguelette approach for a mildly ill-posed inverse problem with fractional noise. We consider a centred Gaussian prior on $\mu_{\windex}$ with $\lambda_{\windex} = \tau 2^{- 2  \alpha j}$ for wavelet coefficients $\windex \in \Windex_{\psi}$, and finite $\lambda_{\windex} \geq 1$ for the scaling coefficients $\windex \in \Windex_{\phi}$. We derive conditions on this prior, namely on its parameters $\tau$ and $\alpha$, so that the corresponding posterior contraction rate is optimal, in the minimax sense. First we state results with known $\beta$ and operator $K$, extending it to the cases when operator $K$ is observed with error and when smoothness $\beta$ is unknown. All the stated results follow from the general results on posterior concentration rates, which apply to more general  Gaussian processes $W$ and operators $K$ (Section~\ref{sec:GeneralRates}), with particular cases considered in Section~\ref{sec:PolyEV}.


\begin{proposition}\label{cor:vaguelettes}
Consider the inverse problem (\ref{f2model}) under fractional noise with Hurst exponent $H \in (1/4,1)$, with the   vaguelette-vaguelette approach described in Section~\ref{sec:vaguevague} under condition (\ref{eq:MildIPwavelets}), for $\mu_0 \in S^{\beta}(A)$ defined by (\ref{eq:DefSobolevWavelet}). We consider  a centred Gaussian prior on $\mu_{\windex}$ with $\lambda_{\windex} = \tau 2^{- 2  \alpha j}$ for wavelet coefficients, and finite $\lambda_{\windex} \geq 1$  for the scaling coefficients.

Then, the contraction rate of the posterior distribution of $\mu$ uniformly over $\mu_0 \in S^\beta(A)$ is given by
$$
\rate \asymp  n^{ -\alpha(1-H)/(p+1-H+\alpha)} + n^{-   \beta (1-H)/(p+1-H+\alpha)} +
   \tau n^{- \alpha (1-H)/(p+1-H+\alpha)} + \tau^{-1}   n^{-\min( 2(p+1-H+\alpha),  \beta)(1-H)/(p+1-H+\alpha)}.
$$
\end{proposition}
Therefore, the vaguelette-vaguelette approach leads to the posterior contraction rate of the same order as if they were an orthonormal basis. This result follows from Theorems~\ref{th:contractiontheoremGeneral} and \ref{contractiontheorem2}.

\begin{corollary}\label{cor:frBgeom0}
In the setting of Proposition~\ref{cor:vaguelettes}, assume $\alpha \geq \beta/2 -p -1 +H$. Then, for $\tau^2   =  \left[ n^{2(1-H)}\right]^{(\alpha-\beta)/(1-H+\beta+p)}$,  the posterior distribution contracts at the minimax convergence rate uniformly over $S^\beta(A)$.

If $\alpha < \beta/2 -p -1 +H$, the contraction rate is suboptimal for any choice of $\tau_{\epsilon}$.
\end{corollary}

It is possible to extend this result to derive the posterior contraction rate for wavelet-based estimation in inverse problems over Besov spaces. However it is known that for some Besov spaces linear estimators are suboptimal, therefore the considered Bayesian  model with a Gaussian prior is not appropriate for a general estimation over Besov spaces,  (see  \citet{AgapiouBesovPrior}  for a Bayesian model with Besov prior). Therefore here we restrict ourselves to Sobolev spaces.

\subsection{Error in operator}\label{sec:ErrorInOperatorFBM}

Now we consider the case when the operator $K$ may be observed with error: $\hat K = K + \delta B$, where the operator $B$ has the standard operator normal distribution, and study its effect on the posterior contraction rate of $\mu$. This problem is discussed in more detail in Section~\ref{sec:ErrorOperatorGeneral}. 


\begin{proposition}
Consider the inverse problem (\ref{f2model}) under fractional noise with Hurst exponent $H \in (1/4,1)$, with the vaguelette-vaguelette approach described in Section~\ref{sec:vaguevague} under condition (\ref{eq:MildIPwavelets}), for $\mu_0 \in S^{\beta}(A)$ defined by (\ref{eq:DefSobolevWavelet}). We consider  a centred Gaussian prior on $\mu_{\windex}$ with $\lambda_{\windex} = \tau 2^{- 2  \alpha j}$ for wavelet coefficients, and finite $\lambda_\windex \geq 1$  for the scaling coefficients. Suppose $\kappa_{\windex}$ are unknown but we observe $\hat \kappa_{\windex} = \langle \hat K, \psi^{(-d)}_{\windex}\rangle$, $\windex \in \Windex_N$ where $\Windex_\phi\subset \Windex_N \subset \Windex$ of size $|\Windex_N| = N$ and $\hat K$ satisfies \eqref{eq:ErrorOperator}.


Then, the rate of contraction of the posterior distribution of $\mu$ given $y$ with plugged in $(\hat \kappa_{\windex}, \, \windex \in \Windex_N)$  is not affected by the error in operator  if
$$
\delta  =o\left( \left[\tau_{\epsilon} n^{2(1-H)}\right]^{-p/(2p+2-2H  +2\alpha)} [\log(\tau_{\epsilon} n^{2(1-H)})]^{-1/2}\right). 
$$
In particular, for non-adaptive models (i.e. where $\beta$ is known), the posterior contraction rate of $\mu$ is optimal, in the minimax sense,  for $\alpha \geq \beta/2 -p -1 +H$ and $\tau_{\epsilon}=[ n^{2(1-H)}]^{(\alpha-\beta)/(1-H+\beta+p)}$ if
$$
\delta =o\left(n^{-(1-H)p/(p+1-H +\beta)} [\log n]^{-1/2}\right). 
$$
\end{proposition}
The proof is a straightforward application of Lemma~\ref{lem:MildIPerrorOp}. 




\subsection{Empirical Bayes posterior distribution of $\mu$}\label{sec:EBtauFBM}

As we have seen in Section~\ref{sec:PostContractionForFBM}, for the posterior distribution to contract at the optimal rate of covergence (in the minimax sense), prior parameters $\alpha$ or $\tau$ need to be chosen appropriately (Corollary~\ref{cor:frBgeom0}), and their choice depends on smoothness $\beta$ of the true function which in practice may be unknown. To adapt to the unknown smoothness of $\mu_0$, we  estimate the scale $\tau$ of the prior covariance operator. In particular, we consider the maximum marginal likelihood estimator of $\tau$ based on the marginal distribution of $Y$ given $\tau$, and study the contraction rate of the posterior distribution of $\mu$ with this estimator as the  plugged in value of $\tau$, which is usually called an empirical Bayes posterior distribution.

 The  marginal distribution of data $y_{\windex}$ for $\windex \in \Windex_\psi$ given $\tau$ is
$$
y_{\windex} \mid \tau \sim N\left(0, \kappa_{\windex}^2 \lambda_{\windex} + \epsilon^2 \sigma^2_{\windex}\right),\,\, \text{ independently for } \windex \in \Windex_\psi,
$$
and $\hat\tau$ is the value maximising this marginal likelihood of $\tau$. See Section~\ref{sec:EBpoly} for further details.
\begin{proposition}\label{prop:vagueFBM_EB}
Consider the inverse problem (\ref{f2model}) under fractional noise with Hurst exponent $H \in (1/4,1)$, with the vaguelette-vaguelette approach described in Section~\ref{sec:vaguevague} under condition (\ref{eq:MildIPwavelets}), for $\mu_0 \in S^{\beta}(A)$ defined by (\ref{eq:DefSobolevWavelet}).
We consider  a centred Gaussian prior on $\mu_{\windex}$ with $\lambda_{\windex} = \tau 2^{- 2  \alpha j}$ for wavelet coefficients, and finite $\lambda_{\windex} \geq 1$  for the scaling coefficients.

Then, 
   the MMLE of $\tau$ satisfies, almost surely for small $\epsilon$,
\begin{equation}\label{eq:AsympTauFracBM}
\hat\tau   \leq   \begin{cases}
\epsilon^{-2(\alpha-\beta)/(1 + p -H + \beta) } (1+o_P(1)), \quad & \text{ if $ \alpha   +1/2 \geq  \beta$},\\
\epsilon^{1/(3/2+ p  -H +\alpha)} (1+o_P(1)), \quad & \text{ if  $\alpha   +1/2 <  \beta$}.
\end{cases}
\end{equation}
Therefore, the posterior distribution with the plugged in estimator $\hat\tau$ is consistent. The contraction rate is optimal in the minimax sense  uniformly over $S^\beta(A)$ for $0<\beta \leq B_0<\infty$ as long as
 \begin{equation}\label{eq:EBalphaFBM}
\alpha \geq \max( B_0 -1/2, B_0/2 -p -1 +H).
\end{equation}
\end{proposition}
Note that for $ \alpha   +1/2 \geq  \beta$, under the assumptions of Proposition~\ref{prop:vagueFBM_EB}, the MMLE almost surely coincides with the oracle  value of $\tau$ which is of the same order as the optimal value given in Corollary~\ref{cor:frBgeom0}. This implies that the corresponding posterior distribution of $\mu$ with the plugged in value of $\tau$ contracts at the optimal rate, in the minimax sense. This proposition is a straightforward  consequence of Theorem~\ref{th:EBtau1} with $\gamma=1/2-H$.


\section{Bayesian inverse problem with correlated Gaussian noise}\label{sec:GeneralRates}

Following the setup for the fractional noise and wavelet based approaches to inverse problems discussed in Section~\ref{sec:RatesWaveletFBM}, we can generalise this approach to other types of correlated noise and other bases.

\subsection{Assumptions}\label{sec:AssumptionsBasesGen}

Recall that in model (\ref{f2model}), $\mu \in H_0 = L^2[0,1]$, and a linear operator $K$ maps $\mu$ into another separable Hilbert space, $H_1= L^2[0,1]$, which are both Hilbert spaces with scalar product $\langle f,g\rangle = \int_0^1 f(x) g(x) dx$.



Now we generalise the vaguelette-vaguelette assumption to a more general set of biorthogonal bases.
\begin{assumption}\label{assumptionB1}
Assume that there exists a countable set of functions $(e_{1, \windex})_{\windex \in \Windex}$  such that
\begin{enumerate}
\item functions $(e_{1,\windex})$ is a basis of $H_1$ such that there exists a biorthogonal basis $(\tilde e_{1,\windex})$,
\item functions $(e_{0,\windex})$ form a Riesz basis of $H_0$ where $e_{0,\windex} = \kappa_{\windex} K^{-1} e_{1,\windex} $ with  $\kappa_{\windex}=1/\left\lVert K^{-1} e_{1,\windex}\right\rVert $.
\end{enumerate}
\end{assumption}
We assume that functions $(e_{1,\windex})$ and $(\tilde e_{1,\windex})$ are also normalised, i.e. $||e_{1,\windex}||=||\tilde e_{1,\windex}||=1$.
 Here the set of indices $\Windex$ is countable, i.e. it is isomorphic to the set of natural numbers $\mathbb{N}$.
 This assumption describes all three considered approaches: wavelet-vaguelette with orthonormal wavelet basis $(e_{1,\windex})$, vaguelette-wavelet with orthonormal wavelet basis $(e_{0,\windex})$, and the vaguelette-vaguelette approach introduced in Section~\ref{sec:vaguevague}.
This assumption  holds also when $(e_{1,\windex})$ form an orthonormal eigenbasis of $K^T K$ which implies that $(e_{0,\windex})$ form an orthonormal eigenbasis of $K K^T$ and $(\kappa_{\windex}^2)$ are the corresponding eigenvalues.

Under this assumption, we can write
\begin{equation}\label{eq:SignalMuDecompB1}
\mu = \sum_{\windex} \mu_{\windex} e_{0,\windex},\quad K\mu = \sum_{\windex} \mu_{\windex} \kappa_{ \windex} e_{1,\windex},
\end{equation}
where $  \mu_{\windex} = \langle K \mu, \tilde e_{1, \windex} \rangle/\kappa_{ \windex}$. Note that here $\mu_{\windex}$ may not be equal to $\langle \mu, e_{0,\windex} \rangle$.

Assume that  noise $W$ in (\ref{f2model}) is such that coefficients $\xi_{\windex} = \langle \tilde e_{1,\windex}, W \rangle$ satisfy
\begin{eqnarray}\label{eq:WnoiseDecomp}
\xi_{\windex} \sim N(0, \sigma_{\windex}^2) \text{ independently for} \windex \in \Windex_{I},\quad \xi_{D} \sim N(0, V_{D}),
\end{eqnarray}
  and $\xi_{D}$ is independent of $\xi_{\windex}$ for all $\windex \in \Windex_{I}\subseteq \Windex$, for a finite set $D =  \Windex\setminus  \Windex_{I}$, some sequence $(\sigma_{\windex})$ and an invertible covariance matrix $V_D$. Note that we do not require $\sum_{\windex\in \Windex} \sigma_{\windex}^2 $ to be finite. Normality and zero mean of $\xi_{\windex}$, $\windex \in \Windex$, follow as long as $W$ is a centred Gaussian process.

\begin{remark} Note that this approach is related to the model of \citet{agapiou_posterior_2013} where given operator $K$ and Assumption~\ref{assumptionB1}, covariance operator of the noise and the prior covariance operator  are assumed to be of the following type
\begin{eqnarray}
V(t,s) &=& \sum_{\windex \in \Windex_I} \sigma_{\windex}^2 e_{1,\windex}(t) e_{1,\windex}(s) + \sum_{\windex \in D} \sum_{\windex' \in D} V_{D,\windex, \windex'} e_{1,\windex}(t) e_{1,\windex'}(s) ,\quad  \sigma_{\windex} \in \mathbb{R}^{\Windex}, \,\, V_D \in \mathbb{R}^{D \times D}, \\
\Lambda(t,s) &=& \sum_{\windex \in \Windex} \lambda_{\windex}  \kappa_{\windex}^2 \, \left(K^{-1} e_{1,\windex}\right)(t) \, \left(K^{-1} e_{1,\windex}\right)(s), \quad  \lambda_{\windex} >0 \, \& \sqrt{\lambda_{\windex}} \in \ell^2(\Windex),
\end{eqnarray}
where $V_D$ is positive definite and $D \subset \Windex$ is finite. In the notation of \citet{agapiou_posterior_2013}, $C_0 = \Lambda$, $C_1 = V$ and $A=K$. In Section~\ref{sec:PolyEV} we compare our posterior contraction rate to the rate given in \citet{agapiou_posterior_2013} for mildly ill-posed inverse problems.
\end{remark}

\subsection{Sequence space formulation}\label{sec:SequenceSpace}

Under  Assumption~\ref{assumptionB1} and the noise in model (\ref{f2model}) satisfying \eqref{eq:WnoiseDecomp}, model (\ref{f2model})
 can be written in terms of  $y_{\windex} = \langle Y, \tilde e_{1,\windex} \rangle$ as
\begin{eqnarray}\label{eq:SeqLik}
y_{\windex} \sim N(\kappa_{\windex} \mu_{\windex}, \epsilon^2 \sigma_{\windex}^2) \text{ independently for } \windex \in \Windex_{I}= \Windex\setminus D,\quad \text{ independently of } y_{D} \sim N(\kappa_{D} \mu_{D}, \epsilon^2 V_{D}),
\end{eqnarray}
where $\kappa_{D} = \diag(\kappa_{\windex}, \windex \in D)$ and  $\mu_D = (\mu_{\windex}, \windex \in D)$.

We also assume that the prior process for $\mu$ is Gaussian and under decomposition (\ref{eq:SignalMuDecompB1}) it satisfies
\begin{eqnarray}\label{eq:SeqPrior}
\mu_{\windex} \sim N(0, \lambda_{\windex})  \text{ independently for } \windex \in \Windex.
\end{eqnarray}
 For the prior distribution of   $\mu$ to be proper, we assume $\sum_{\windex\in \Windex} \lambda_{\windex} <\infty$. Denote $\Lambda_D = \diag(\lambda_{\windex}, \windex \in D)$.

The corresponding  posterior distribution $\mu_{\windex} \mid y_{\windex}$ for each $\windex \in \Windex$ is
\begin{eqnarray}\label{eq:SeqPost}
\mu_{\windex} \mid y_{\windex} \sim
N\left(\frac{ y_{\windex} \kappa_{\windex}\lambda_{\windex}}{ \lambda_{\windex} \kappa_{\windex}^2 + \epsilon^2 \sigma^2_{\windex}}, \frac{\sigma^2_{\windex}\lambda_{\windex}}{\epsilon^{-2} \lambda_{\windex} \kappa_{\windex}^2 + \sigma^2_{\windex}}\right), \,\, \windex\in \Windex_I
\end{eqnarray}
independently, and independently of
\begin{eqnarray}\label{eq:SeqPostD}
\mu_{D} \mid y_{D} \sim
 N\left(\left( K_D^T  V_D^{-1} K_D +\epsilon^{2}\Lambda_D^{-1}\right)^{-1} K_D^T V_D^{-1} y_D,   \left(\epsilon^{-2} K_D^T  V_D^{-1} K_D +\Lambda_D^{-1}\right)^{-1}\right).
\end{eqnarray}


Then, if estimates $\hat\mu_{\windex}$ for $\windex \in \Windex$ are available, we obtain the estimate of function $\mu$: $\hat\mu = \sum_{\windex \in \Windex} \hat\mu_{\windex} e_{0,\windex}$. Similarly, a posterior distribution on  $(\mu_{\windex},\, \windex\in \Windex)$ induces a posterior distribution on  function $\mu$.

In most of the paper, we consider this sequence model under the assumption that  $(e_{0,\windex})$ and $(\tilde e_{1,\windex})$ are known. In the case $(e_{0,\windex})$ form an orthonormal eigenbasis of $K K^T$, \citet{KoltchinskiiLounici} propose a way to estimate the corresponding eigenvectors of the covariance matrix in the corresponding finite dimensional model which we discuss in Section~\ref{sec:EstimateEigenFunctions}.



\subsection{Minimax rate of convergence}\label{sec:MinimaxGeneral}

In this section we define the minimax rate of convergence of estimators of $\mu$ under model (\ref{f2model}). This rate is usually used as a benchmark for posterior contraction rates \cite{ghosal_convergence_2000}.  We consider the following smoothness classes for the true unknown function $\mu_0$.
\begin{assumption} \label{smoothGeneral}
Assume that $\mu_0$ belongs to a smoothness class $Q((a_{\windex}), A)$, with $A >0$ and a non-decreasing sequence $(a_{\windex})_{\windex \in \Windex}$, $a_{\windex}\geq 1$:
\begin{equation}\label{eq:DefSeqClassGeneral}
Q((a_\windex), A)=\left\{f = \sum_{\windex \in \Windex} e_{0,\windex} f_{\windex}: \, \sum_{\windex} a_{\windex}^{2} f_{\windex}^2 \leq A^2 \right\}.
\end{equation}
\end{assumption}
Since $\Windex$ is isomorphic to $\mathbb{N}$, the sequence $(a_{\windex})$ is understood to be non-decreasing for some map $\Windex \to \mathbb{N}$.
 For example, a Sobolev class is defined with $e_{0,\windex}$ being the Fourier basis, $\Windex = \mathbb{Z}$ and $a_i = |i|^\beta$,  or in terms of wavelets or vaguelettes as discussed in Section~\ref{sec:SobolevWavelets}. 

Define the risk of estimator $\hat \mu$  of a true function $\mu_0$  in $L^2 $ norm over a set of functions $Q$ by
\begin{equation}\label{eq:minimaxrisk}
R(\hat\mu, Q) = \sup_{\mu_0 \in Q} \mathbb{E}_{\mu_0} ||\hat\mu-\mu_0||^2.
\end{equation}

\begin{definition}\label{def:minimaxrate}
$r_{\epsilon}(Q)$ is the minimax rate of convergence of estimators of $\mu$ under model (\ref{f2model}) over class $Q$  if $\exists 0 < c \leq C <\infty$ that depend only on $Q$, $K$ and $V$ such that $c\leq \inf_{\hat\mu} [r_{\epsilon}(Q)]^{-2} R(\hat\mu, Q) \leq C$,
 where $R(\hat\mu, Q)$ is defined by (\ref{eq:minimaxrisk}).
\end{definition}


The minimax rate of convergence of estimating $\mu$ in $L^2$ norm under model (\ref{f2model}) in sequence space under Assumption~\ref{smoothGeneral} can be derived for given sequences $(\sigma_{\windex})$, $(\kappa_{\windex})$ and $(a_{\windex})$, ${\windex \in \Windex}$,  using Theorem 3 in \citet{belitser_minimax_1995} which is stated in Online supplementary material (Lemma~\ref{lem:BelitserLevit}). 
  The finite-dimensional dependent part of the model in sequence space (\ref{eq:SeqLik}) adds a term of order $\epsilon^2$ to the minimax rate as long as its covariance matrix $V_D$ has finite positive eigenvalues.

\begin{remark}\label{rem:DepMinimax} If the covariance matrix $V_D$ is such that $||V_D||, \, ||V_D^{-1}|| \in (0,\infty)$, then the part of the model (\ref{eq:SeqLik}) with dependent observations $y_D$ is stochastically dominated by an iid model with larger variances $y_\windex \sim N(\kappa_\windex \mu_\windex, \epsilon^2 ||V_D||^2)$, $\windex\in D$, and stochastically dominates an iid model with smaller variances $\epsilon^2 ||V_D^{-1}||^{-2}$.  Since the set $D$ is finite, in both cases these terms affect the minimax rate up to an additional $C \epsilon^2$.
\end{remark}
 The following corollary, for when a fast rate can be obtained, is a straightforward consequence of Lemma~\ref{lem:BelitserLevit}.
\begin{corollary}\label{cor:ParamRate}
Under conditions of Lemma~\ref{lem:BelitserLevit} with $\tilde \sigma_\windex= \sigma_\windex \kappa^{-1}_\windex$,
\begin{enumerate}
\item if $\sum \limits_{\windex \in \Windex} \kappa^{-2}_\windex \sigma_\windex^2 <\infty$ then the minimax rate of convergence is $r_\epsilon(Q)=C\epsilon$ for some $C>0$;
\item if   $\sum \limits_{\windex \in \Windex_N} \kappa^{-2}_\windex \sigma_\windex^2\asymp   \log N$ for large $N$ and $\Windex_N$ is the set of the first $N$ elements of $\Windex$ (ordered so that sequence $(a_\windex)$ is non-decreasing),  then the minimax rate of convergence is $r_\epsilon^2(Q)=C\epsilon^2\log(1/\epsilon)$.
\end{enumerate}
\end{corollary}

\subsection{Posterior contraction rates and their optimality}

Now we study conditions on the prior distribution so that the corresponding posterior distribution $\mu|Y$ contracts to the true value $\mu_0$ at the optimal rate, in the minimax sense. 

\subsubsection{Rate of contraction of posterior distribution}\label{sec:RateGeneral}

Below we present a general result that can be applied for any $a_\windex$, $\lambda_\windex$, $\kappa_\windex$ and $\sigma_\windex$ satisfying stated conditions.


\begin{theorem}\label{th:contractiontheoremGeneral}
Consider the inverse problem (\ref{f2model}) formulated in the sequence space under Assumption~\ref{assumptionB1} with noise $W$ satisfying (\ref{eq:WnoiseDecomp}) and prior distribution (\ref{eq:SeqPrior}). Assume that the true function $\mu_0$ satisfies Assumption~\ref{smoothGeneral}. Assume also that $\min_{\windex \in D}\lambda_{\windex} \geq 1$ and $\left\lVert
K_D^{-1} V_D K_D^{-1}\right\rVert  < \infty$.

Then, for every $M \rightarrow \infty$,
$$
\mathbb{P}\left(\mu: \, ||\mu-\mu_0|| \geq M \rate|Y\right) \stackrel{\mathbb{P}_{\mu_0}}{\to} 0, \text{~ as ~} \epsilon \rightarrow 0,
$$
uniformly over $\mu_0$ in  $Q((a_\windex), A)$ where $\rate$ is given by
$$
\rate  =  \epsilon+\left[\epsilon^2\sum_{\windex \notin \Windex_{\epsilon}}  \sigma^2_{\windex} \kappa_{\windex}^{-2}+    A^2 \sup_{\windex \in \Windex_{\epsilon}} a_{\windex}^{-2} +
 \sum_{\windex \in \Windex_{\epsilon}}  \lambda_{\windex} + \epsilon^4  \max_{\windex \notin \Windex_{\epsilon}} \left[\frac{\sigma^2_{\windex}}{  a_{\windex} \kappa_{\windex}^2\lambda_{\windex} }  \right]^2 \right]^{1/2},
$$
 and $\Windex_{\epsilon}=\{\windex \in \Windex: \, \sigma^{2}_{\windex} \kappa_{\windex}^{-2} > \epsilon^{-2} \lambda_{\windex} \}$. 
\end{theorem}
 The first two terms in the square brackets represent the variance and squared bias respectively, and the remaining terms, which involve prior parameters $\lambda_{\windex}$, represent the bias arising from the prior and the saturation term.


\begin{remark}
The posterior distribution can contract at rate $\epsilon$ if  $(\sigma^2_{\windex})$ is such that conditions (\ref{eq:condDep}) hold and $\sum_{\windex\in \Windex_I} \sigma^2_{\windex} \kappa_{\windex}^{-2}\leq C < \infty$  which is the minimax rate (Corollary~\ref{cor:ParamRate}), under the appropriate choice of prior parameters $(\lambda_{\windex})$. 
\end{remark}

\begin{remark}
There is a phenomenon known as saturation \cite{AgapiouMathe2017} that constrains the posterior contraction rate for an undersmoothing prior. In the above theorem, this is described  by the term $\max_{\windex \notin \Windex_{\epsilon}} \left[\frac{\sigma^2_{\windex}  }{ a_{\windex} \kappa_{\windex}^2\lambda_{\windex} }  \right]$. It will be illustrated in the example below.
\end{remark}

\begin{example}
Motivated by an example in \citet{agapiou_posterior_2013} with covariance operator of the type $( (K K^T)^{-1} + M_x )^{-1}$ where $M_x$ is a nonlinear operator, we consider a particular linear case with $V = (K K^T)^a$, $a\geq 0$. Here $(e_{1,i})$ are the eigenbasis of $K K^T$ and $V$, with eigenvalues $\kappa_i^2$ and $\sigma_i^2 = \kappa_i^{2a}$, respectively.   In this case,
$$
\rate  =\left[\epsilon^2\sum_{i\leq i_{\epsilon}}  \kappa_i^{-2+2a}+  a_{i_{\epsilon}}^{-2} +
 \sum_{i > i_{\epsilon}}  \lambda_i + \epsilon^4  \max_{i\leq i_{\epsilon}} \left[\frac{ 1}{ a_i \kappa_i^{2(1-a)}\lambda_i }  \right]^2 \right]^{1/2},
$$
 where $i_{\epsilon}=\max\left\{i: \, \kappa_i^{-2(1-a)} \leq \epsilon^{-2} \lambda_i \right\}$. In particular, if $a=1$, i.e. if $V = K K^T$, then the posterior contraction rate coincides with the rate of the direct problem. This corresponds to the model of the type $Y = K(\mu + \epsilon W)$ where $W$ is white noise i.e. when the error occurs before operator $K$ is applied.
\end{example}


In the next section we assume that the variances $(\sigma^2_{\windex})$ are unknown, and their plug-in estimator is available.




\subsubsection{Rate of contraction of posterior distribution with unknown variances $(\sigma^2_{\windex})$}\label{sec:RateGeneralPlugIn}

When the $(\sigma^2_{\windex})$ are unknown, their plug-in estimator is used to conduct inference about $\mu$. We investigate how this affects the contraction rate of the posterior distribution of $\mu$. Theorem~\ref{th:contractiontheoremGeneral}  will be used to address this question.

Suppose we have a plug-in estimator $\{\hat \sigma_{\windex}^2\}$ of the error variances $\sigma^2_{\windex}$, $\windex \in \Windex$ under the model  \eqref{eq:SeqLik}. We consider the case where $\hat \sigma_\windex^2$  are consistent estimators of $\sigma^2_\windex$ for all $\windex$ and are independent of $Y$ used to estimate $(\mu_\windex)$. Plugging in an estimator can be thought of as having an informative prior distribution on $\sigma^2_\windex$ that is a point mass at $\hat \sigma^2_\windex$. If $\sigma^2_{\windex}$ are the eigenvalues of $V$, \citet{KoltchinskiiLounici} proposed a way to estimate the eigenfunctions of $V$ which we apply for a particular form  of $V$ under repeated observations in Section~\ref{sec:EstimateEigenFunctions}.

These assumptions are satisfied, for instance, when variances are estimated from a different study; when the data is split into two independent parts: one is used to estimate $(\sigma_\windex)$ and the other one to estimate $(\mu_\windex)$; or when there are repeated observations with Gaussian error, (hence the sample mean and the sample variance are independent),  convolution operators \cite{CavalierHengartner} and statistical inference in econometric problems with instruments \cite{FlorensSimoni}. Another application is an additive error in the operator \cite{HoffmannReiss} that we consider below.

We make the following assumption of consistency of $\hat\sigma_\windex^2$ that holds relative to $\sigma_\windex^2$.
\begin{assumption} \label{AssumeEVVarPlugInRelative}
Assume that the estimated variances $\left(\{\hat\sigma_\windex^2\}_{\windex \in \Windex_N}\right)$  for some $\Windex_N \subset \Windex$ of size $|\Windex_N| = N$  are independent of $(y_{\windex}, \windex \in \Windex)$, and there exists a constant $c_0$ such that $c_0 \epsilon_\sigma\leq 1/2$ and
$$
P\left(|\hat \sigma_\windex^2/\sigma_\windex^2-1|\leq c_0 \epsilon_\sigma, \, i=1,2,\ldots, N\right) \to 1, \text{ as }  \epsilon_\sigma \to 0 \text{ and } N\to \infty.
$$
\end{assumption}
Note that the rate $\epsilon_\sigma$ may depend on $N$ but it does not need to be known.


\begin{theorem}\label{th:contractiontheoremGeneralPlugIn2}
Consider the inverse problem (\ref{f2model}) formulated in the sequence space under Assumption~\ref{assumptionB1}, with noise $W$ satisfying \eqref{eq:WnoiseDecomp} and prior distribution (\ref{eq:SeqPrior}) with $\lambda_{\windex} =0$ for $\windex \notin \Windex_N$. We assume that $\min_{\windex \in D}\lambda_{\windex} \geq 1$ and $\left\lVert K_D^{-1} V_D K_D^{-1}\right\rVert  < \infty$. Assume also that the true function $\mu_0$ satisfies Assumption~\ref{smoothGeneral}.

Suppose that $\sigma_{\windex}$ are unknown, but we observe $(\hat\sigma_\windex, \, \windex \in \Windex_N)$ satisfying Assumption~\ref{AssumeEVVarPlugInRelative}.

Then, the posterior distribution of $\mu$ given $Y$, with plugged in $(\hat\sigma_\windex)$ instead of $\sigma_\windex$, is such that
for every $M \rightarrow \infty$,
$$
\sup_{\mu_0 \in S^\beta(A)}\mathbb{P}\left(\{\mu: \, ||\mu-\mu_0|| \geq M \rate|Y, \hat V, N\}\right) \stackrel{\mathbb{P}_{\mu_0, V}}{\to} 0, \text{~ as ~} \epsilon \rightarrow 0, \epsilon_\sigma \to 0 \text{ and } N\to \infty,
$$
 where $\rate$ is given by
\begin{eqnarray*}
\rate^2  =
  \epsilon^{2} \sum_{ \windex \in \Windex_N \, \& \, \windex \notin \Windex_{\epsilon,-}} \sigma^2_{\windex} \kappa_{\windex}^{-2}
 +  \sum_{  \windex \in \Windex_N \, \& \, \windex \in \Windex_{\epsilon,+}}  \lambda_{\windex} +
 \epsilon^{-2}\sum_{  \windex \in \Windex_N \, \& \, \windex \in \Windex_{\epsilon,-} }\sigma_{\windex}^{-2} \kappa_{\windex}^2\lambda_{\windex}^2
+    \epsilon^4   \max_{  \windex \in \Windex_N \, \& \, \windex \notin \Windex_{\epsilon,+}} \left[\frac{\sigma^2_{\windex} a_{\windex}^{-1}}{  \kappa_{\windex}^2\lambda_{\windex} }  \right]^2
+ \sup_{  \windex \notin \Windex_N \, \& \, \windex \in \Windex_{\epsilon,+}} a_{\windex}^{-2},
\end{eqnarray*}
where $\Windex_{\epsilon,\pm} = \left\{\windex: \,  \sigma^2_\windex/[\lambda_\windex\kappa_\windex^2]> (1 \pm c_0\epsilon_\sigma)^{-1}\epsilon^{-2}\right\}$.

\end{theorem}
In addition to the terms in the rate with known $V$ in Theorem~\ref{th:contractiontheoremGeneral} where $\Windex_{\epsilon}$ replaced by  either $\Windex_{\epsilon,\pm}$ or $\bar{\Windex}_{N}$, there is an additional term (third term in the rate above) that arises from the posterior variance.
 We apply this theorem to the case of repeated observations in Section~\ref{sec:RepeatedObs}, where we estimate $\sigma_{\windex}$ and  the power  for polynomially decreasing $\sigma^2_{\windex}$ (Section~\ref{sec:PolyEstGamma}).  We also discuss the case where the difference between $\hat\sigma_\windex^2$ and $\sigma_\windex^2$ are simultaneously bounded for all $\windex \in \Windex$ (Section~\ref{sec:PlugInAbsBound} in Online supplementary material).

 Next we discuss the problem of error in the forward operator $K$, that can be reformulated as a problem of unknown variances.

\subsubsection{Error in operator}\label{sec:ErrorOperatorGeneral}

We illustrate this on an example with error in forward operator $K$, under model (\ref{eq:ErrorOperator}) introduced in Section~\ref{sec:ErrorOperatorIntro}.  Under the settings of Section~\ref{sec:AssumptionsBasesGen}, this problem can be written in the sequence space as (\ref{eq:SeqLik}) with
$$
 \hat{\kappa}_\windex = \kappa_\windex + \delta \xi_\windex, \quad  \xi_\windex \sim N(0,1)  \text{ \, iid \, for} \quad  \windex \in \Windex_N,
$$
and independently of $y_{\windex}$, for some set $ \Windex_N$ such that  $D \subset \Windex_N \subset \Windex$ of size $|\Windex_N| = N$.
 This model, together with sequence space model (\ref{eq:SeqLik}), can be rewritten as
 \begin{eqnarray}\label{eq:SeqErrorOperator}
\hat y_\windex= y_\windex/\hat{\kappa}_\windex = \mu_\windex + \epsilon \hat\sigma_\windex z_\windex  \quad \text{ independently of }  \quad \hat{\kappa}_\windex = \kappa_j + \delta \xi_\windex, \quad  \windex \in \Windex_N,
 \end{eqnarray}
with $\tilde\sigma_\windex= \sigma_\windex{\kappa}_\windex^{-1}$ and $\hat\sigma^2_\windex=\sigma_\windex^2\hat{\kappa}_\windex^{-2}$ for known $\sigma_\windex$. Thus, this problem can be interpreted as a direct problem with unknown variance, and  we can study its effect on the posterior contraction rate using Theorem~\ref{th:contractiontheoremGeneralPlugIn2}.
 To apply the theorem, we need to verify that Assumption~\ref{AssumeEVVarPlugInRelative}   holds for model (\ref{eq:SeqErrorOperator}) which we do in the following lemma.

\begin{lemma}\label{lem:GeneralKwithError}
Consider model (\ref{eq:SeqErrorOperator}). If for some $t\geq 1$, which may depend on $\epsilon$ and $\delta$,
\begin{equation}\label{cond:ErrorOperator}
A_{N,\delta} := \delta^{-1} [t+\log N]^{-1/2}\min_{\windex \in \Windex_N} \kappa_\windex \to \infty \text{ as } \max(\epsilon,\delta)\to 0,
\end{equation}
then, with probability at least $1-e^{-t/2}$,
$$
|\hat{\sigma}^2_\windex/\tilde\sigma_\windex^2 -1| \leq  
\sqrt{2}A_{N,\delta}^{-1}\left(1 - A_{N,\delta}^{-1}\right)^{-3/2}, \quad \windex \in \Windex_N.
$$
\end{lemma}
Therefore, Assumption~\ref{AssumeEVVarPlugInRelative} holds under condition (\ref{cond:ErrorOperator}), e.g. with $t = 2 \log N$. This condition effectively specifies truncation $N$, i.e. the number of components for which it is possible to estimate $\{\kappa_\windex^{-1}, \, \windex \in \Windex_N \}$ consistently.

\subsubsection{Repeated observations}\label{sec:RepeatedObs}


Now we suppose that we have $m$ independent replicates of the original model (\ref{f2model}) with $\epsilon$ replaced by $\epsilon_0$, and hence, under the Assumptions of Section~\ref{sec:AssumptionsBasesGen}, we can write the corresponding model \eqref{eq:SeqLik} as
\begin{equation} \label{eq: prob model:rep obvs}
Y_{\windex,j} \sim N\left(\kappa_\windex\mu_{\windex},\epsilon_0^2\sigma_{\windex}^2 \right),~~ \windex\in \Windex_I,  \quad Y_{D,j} \sim N\left(K_D \mu_{D},\epsilon_0^2 V_D\right), \quad  j=1,\dots, m,
\end{equation}
independently.   Consequently, for each $\windex \in \Windex_I$, the sample mean ($\bar Y_\windex$) and the sample variance ($s_\windex^2$) are
\begin{eqnarray} \label{eq:probmodel:repobvs:ybar}
\bar Y_\windex &:=& \frac{1}{m} \sum_{j=1}^{m} Y_{\windex,j}  \sim N\left( \kappa_\windex\mu_{\windex},\epsilon_0^2\sigma_{\windex}^2/m \right),\\ s_{\windex}^2 &:=&  \frac{1}{\epsilon_0^2 (m-1)} \sum_{j=1}^{m}(Y_{\windex,j}- \bar Y_{\windex})^2 \sim \frac{\sigma_{\windex}^2}{m-1}\chi^2_{m-1}, \notag
\end{eqnarray}
independently for all $\windex \in \Windex_I$. Also, independently from $\bar Y_\windex$ and $s_{\windex}^2$ for $\windex \in \Windex_I$,
\begin{eqnarray} \label{eq:probmodel:repobvs:ybarD}
\bar Y_D &:=& \frac{1}{m} \sum_{j=1}^{m} Y_{D,j}  \sim N\left( K_D \mu_{D},\epsilon_0^2V_D/m \right),\\
S_{D}^2 &:=&  \frac{1}{\epsilon_0^2 (m-1)} \sum_{j=1}^{m}\left(Y_{D,j}- \bar Y_{D}\right)\left(Y_{D,j}- \bar Y_{D}\right)^T  \sim \frac{1}{m-1} \text{ Wishart}_D( m-1, V_D),\notag
\end{eqnarray}
where a $D$-dimensional Wishart distribution of a random positive definite matrix, parameterised by a positive definite matrix $\Sigma$ and shape parameter $a>0$, has density $f_{D}(S; a, V) =  |S|^{(a-D-1)/2}e^{-\trace(S^T V)/2} 2^{-a D/2} |V|^{-a/2}/\Gamma_D(a/2)$ where $\Gamma_D$ is a multivariate Gamma function.

Now we study when simultaneous asymptotic consistency  of  the following estimator of $(\sigma_\windex^2)_{\windex \in \Windex}$ holds  for large $m$ with high probability, given a set $\Windex_N$ such that $D \subset \Windex_N \subset \Windex$,  $|\Windex_N|=N$ and
\begin{equation}\label{eq:repobs:sigmahat}
\hat{\sigma}_\windex^2 = s_\windex^2 I(\windex \in \Windex_N).
\end{equation}


\begin{proposition} \label{prop: rep obvs M contraction2}
 Assume that we have $m$ independent observations of inverse problem (\ref{f2model}) that can be formulated in the sequence space under Assumption~\ref{assumptionB1}, with noise $W$ satisfying \eqref{eq:WnoiseDecomp}, with repeated observations (\ref{eq: prob model:rep obvs}), with prior distribution (\ref{eq:SeqPrior}), and the true function $\mu_0$ satisfying Assumptions~\ref{smoothGeneral}. 

Consider the estimator of $(\sigma_i^2)$ defined by (\ref{eq:repobs:sigmahat}) with $N$ satisfying $\log N =o(m)$.
 Then, for every $M_0 \rightarrow \infty$,
$$
 \mathbb{P}\left( \mu: \, ||\mu-\mu_0|| \geq M_0 {\rate}_{plugin}|Y, (\hat \sigma_\windex^2)\right)   \stackrel{\mathbb{P}_{\mu_0, V}}{\to} 0, \text{ as } \epsilon \rightarrow 0,
$$
uniformly over $\mu_0 \in Q((a_\windex), A)$, where ${\rate}_{plugin}$ is the rate given in Theorem~\ref{th:contractiontheoremGeneralPlugIn2} with $\epsilon^2=\epsilon_0^2/m$.

If also $\Windex_N \subseteq \Windex \setminus \Windex_\epsilon$ and $\Windex_{\epsilon,\pm} \asymp \Windex_{\epsilon}$ (the latter holds as long as $\epsilon \to 0$ and $m\to \infty$), then ${\rate}_{plugin}= \rate$ stated in Theorem~\ref{th:contractiontheoremGeneral}, i.e. it is not affected by the plug-in estimator.


\end{proposition}
The condition on $N$, $\log N =o(m)$, is required to have simultaneous consistency for $N$ estimators given their precision $m$, and the other condition $\Windex_N \subseteq \Windex \setminus \Windex_\epsilon$ ensures the plugin rate of posterior contraction of $\mu$ coincides with the contraction rate for known $(\sigma_\windex^2)$.

In the next section we illustrate how the results of this section apply to mildly ill-posed problems with correlated Gaussian noise.


 \section{Mildly ill-posed inverse problems with heterogeneous Gaussian noise}\label{sec:PolyEV}

In this section we apply the general results with known and unknown covariance operators of  Gaussian noise to mildly ill-posed inverse problems, under the following assumptions.

\subsection{Assumptions}\label{sec:AssumptionsGeomIP}

In this section, we assume that the setting of Section~\ref{sec:AssumptionsBasesGen} applies, leading to the sequence space problem described in Section~\ref{sec:SequenceSpace}, under mildly ill-posed inverse problems and other sequences decaying polynomially. We set $\Windex = \mathbb{N}$, 
 taking an appropriate mapping, such that sequence $a_i$ in the definition of the smoothness class $Q( (a_i), A)$ is non-increasing.
\begin{assumption} \label{AssumeEVVar}
  Assume that ($\sigma_i^2$) satisfy
 $C_2^{-1}i^{\gamma} \leq \sigma_i \leq C_2i^{\gamma}$
 for some $\gamma \in \mathbb{R}$ and $C_2 \geq 1$.
\end{assumption}

\begin{assumption} \label{AssumeSmoothPoly}
  Assume that ($a_i$)   satisfy
 $C_3^{-1}i^{\beta} \leq a_i \leq C_3 i^{\beta}$
 for some $\beta >0$ and $C_3 \geq 1$.
\end{assumption}
This assumption corresponds to a generalised Sobolev space. If $(e_{0,i})$ is the Fourier basis and $C_3 =1$, this is the definition of a standard Sobolev class. If $(e_{0,i})$ are wavelets then this also defines the standard Sobolev class for the mapping $(j,k)$ to $i = 2^j +k$ and for $C_3 = 2^{\beta}$, given sufficiently regular wavelets.


\begin{assumption} \label{AssumeEpriorV}
Eigenvalues   $(\lambda_i)$  satisfy
  $\lambda_i=\tau_{\epsilon}^2\,i^{-1-2\alpha}$
 for some $\alpha >0$ and $\tau_{\epsilon} > 0$, such that $\epsilon^{-2}\tau_{\epsilon}^2 \rightarrow \infty$ as $\epsilon \to 0$. Denote $\tau=\tau_\epsilon^2$.
\end{assumption}
The latter assumption  implies that  a priori we assume $\mu \in Q\left((i^{\alpha'}), R\right) = S^{\alpha'}(R)$ almost surely, for any $\alpha'<\alpha$ and large enough $R$, for a fixed $\epsilon$.

\begin{assumption} \label{AssumeDepPoly}
$||V_D||, ||V_D^{-1}|| \in (0,\infty)$ and there exist constants $C_{D,1}, \, C_{D,2} >0$ independent of $\epsilon$ such that
\begin{eqnarray*}
\left\lVert(K_D^T  V_D^{-1} K_D +\epsilon^{2}\Lambda_D^{-1})^{-1} \Lambda_D^{-1}\right\rVert \leq C_{D,1} \quad \trace\left( (K_D^T  V_D^{-1} K_D +\epsilon^{2}\Lambda_D^{-1})^{-1} \right)\leq C_{D,2}.
\end{eqnarray*}
\end{assumption}
By Lemma~\ref{lem:PostRateDepTerms}, this assumption ensures that the terms in the posterior contraction rate that correspond to the dependent observations are or order $\epsilon$.



\subsection{Minimax rate of convergence}

The minimax rate  for model $(\ref{f2model})$ with the considered parameter sequences is stated in the following proposition.
\begin{proposition}\label{prop:optimalrateshetero}
Consider the inverse problem (\ref{f2model}) formulated in the sequence space under Assumption~\ref{assumptionB1} with noise $W$ following \eqref{eq:WnoiseDecomp}. Assume that  Assumption~\ref{AssumeEVVar} holds, that matrix $V_D$ has finite positive eigenvalues  and that the true function $\mu_0$ satisfies  Assumptions~\ref{smoothGeneral} and \ref{AssumeSmoothPoly}. 

Then, the minimax rate of convergence of estimating $\mu$ is given by
\[\rate^* :=  \begin{cases}
\epsilon^{\frac{2\beta}{1 + 2\beta + 2(p+\gamma)}}, & \text{if $\gamma > -p-1/2$}.\\
\epsilon (\log |\epsilon|)^{1/2}, & \text{if $\gamma = -p-1/2$}.\\
\epsilon, & \text{if $-\beta/2 -p-1/2<\gamma < -p-1/2$}.
\end{cases} \]
\end{proposition}
It follows from Lemma~\ref{lem:BelitserLevit} together with Remark~\ref{rem:DepMinimax}.

This result implies that for a heterogeneous variance the degree of ill-posedness for $(\ref{f2model})$ changes, in particular to i.e. $\tilde{p} = p+\gamma$ if $p+\gamma >-1/2$. For $p+\gamma =0$, the rate coincides with the minimax rate of the direct problem.
 If $p+\gamma\leq -1/2$, the problem becomes self-regularised, i.e. the  rate of convergence $\epsilon$ can be achieved (up to a log factor in the case $p+\gamma=-1/2$), 
 provided $\gamma  +p+1/2 + \beta/2 >0$. According to \citet{belitser_minimax_1995}, under this constraint the minimax rate of any estimator coincides with the minimax rate of a linear estimator. Since we consider Gaussian errors and Gaussian priors, such constraint is appropriate.

\subsection{Posterior contraction rates}\label{sec:PostContractionPoly}

Having found the minimax rates, we can now discuss the contraction rates achieved by the posterior distribution under the considered Bayesian model. Note these rates also apply when the problem is self-regularised, i.e. $p+\gamma \leq -1/2$.

\begin{theorem}\label{contractiontheorem2}
Consider the inverse problem (\ref{f2model}) formulated in the sequence space under Assumption~\ref{assumptionB1} with noise $W$ satisfying \eqref{eq:WnoiseDecomp} and Assumption~\ref{AssumeEVVar}, with prior distribution (\ref{eq:SeqPrior}) under Assumption~\ref{AssumeEpriorV}, and the true function $\mu_0$ satisfying Assumptions~\ref{smoothGeneral} and \ref{AssumeSmoothPoly}. Let Assumption~\ref{AssumeDepPoly} hold.  

Then, for every $M \rightarrow \infty$, \, $\mathbb{E}_{\mu_0}\mathbb{P}\left(\mu: \, ||\mu-\mu_0|| \geq M \rate|Y\right) \rightarrow 0  \text{ as } \epsilon \rightarrow 0$ \,  uniformly over $\mu_0$ in $S^\beta(A)$  where
\[ \rate:=  \begin{cases}
 (\epsilon^{2} \tau_{\epsilon}^{-2})^{\frac{\beta}{1 + 2\alpha + 2(p+\gamma)} \wedge 1} + \tau_{\epsilon}(\epsilon^{2}\tau_{\epsilon}^{-2})^{\frac{\alpha}{1 + 2\alpha + 2(p+\gamma)}}, & \text{if  $\gamma > -p-1/2$,}\\
 (\epsilon^{2}\tau_{\epsilon}^{-2})^{\frac{\beta}{2\alpha} \wedge 1} + \epsilon[\log (\epsilon^{-1}\tau_{\epsilon})]^{1/2}, & \text{if  $\gamma = -p-1/2$,}\\
(\epsilon^{2}\tau_{\epsilon}^{-2})^{\frac{\beta}{1 + 2\alpha + 2(p+\gamma)} \wedge 1} + \epsilon, & \text{if  $-p -1/2-\alpha <\gamma < -p-1/2$.}
\end{cases} \]
\end{theorem}
The assumption of monotonicity of $\sigma_i^2 /[\lambda_i\kappa_i^2]\asymp i^{1 +2\gamma+2\alpha +2p}$ for large $i$ from Theorem~\ref{th:contractiontheoremGeneral} is reflected in the condition $ -p -1/2-\alpha <\gamma$.
  Recall that parameters $p$ and $\gamma$ are assumed known and given by the problem, as well as the smoothness parameter $\beta$.
Parameters of the prior $\alpha$ and $\tau_{\epsilon}$ can be chosen in some cases so that the posterior contracts at the optimal rate given in Proposition~\ref{prop:optimalrateshetero}.


\begin{corollary}\label{Cor:OptPriorParam}
Let assumptions of Theorem \ref{contractiontheorem2} hold.
 The rate of contraction of the posterior given in Theorem~\ref{contractiontheorem2} matches the minimax rate of convergence, for the following $\alpha$ and $\tau_{\epsilon}= \tau^{1/2}$:
\begin{enumerate}
\item  $\tau_{\epsilon}=const \in (0,\infty)$ and
\[\begin{cases}
\alpha=\beta, & \text{if $\gamma > -\frac{1+2p}{2}$},\\
\alpha\leq\beta, & \text{if $\gamma = -\frac{1+2p}{2}$},\\
\alpha\leq\beta, & \text{ if $-\frac{1+2p}{2} - \alpha<\gamma < -\frac{1+2p}{2}$};
\end{cases} \]

\item for $\tau_{\epsilon}$ depending on $\epsilon$:
\[  \begin{cases}
    \alpha \geq \beta/2 -  (1/2+p+\gamma), \quad \tau_{\epsilon} =C \epsilon^{\frac{2(\beta-\alpha)}{1 + 2\beta + 2(p+\gamma)}}, & \text{if $\gamma > -\frac{1+2p}{2}$},\\
  \epsilon^{-B} \geq \tau_{\epsilon} \geq C\epsilon^{1 -  \max(1/2, \alpha/\beta)} [\log \epsilon^{-1}]^{-
0.5\max(1/2, \alpha/\beta)}, & \text{if $\gamma = -\frac{1+2p}{2}$},\\
  \tau_{\epsilon} \geq C \epsilon^{\min(1/2, 1-\frac{1+2\alpha +2(p+\gamma)}{2\beta})}, & \text{if $-\frac{1+2p}{2} - \alpha<\gamma < -\frac{1+2p}{2}$}.
\end{cases}
\]

\end{enumerate}
\end{corollary}

Observe that when $\gamma > -p-1/2$, the rates obtained are similar to the white noise case, albeit with a different degree of ill-posedness  $\tilde{p}=p+\gamma$. When $p+\gamma<0$, the fastest rate of contraction coincides with the minimax rate of convergence of the direct problem, i.e. the model self-regularises, and it can achieved when we undersmooth a priori.

If $\alpha <  \beta/2 - (1/2 + p+\gamma)$ and $\gamma > -p-1/2$, i.e. if we undersmooth too much a priori, then the minimax rate cannot be achieved for any $\tau$. When $\gamma =0$, this coincides with the findings of \cite{knapik_bayesian_2011}. This is known as saturation \cite{AgapiouMathe2017}.  However, in the self-regularising case $-p-1/2 - \alpha<\gamma \leq -p-1/2$ the optimal rate can be achieved if the appropriate prior scaling is used.

\begin{remark}\label{rem:OptimalCutoff}
Note that  the case $p+\gamma +1/2 >0$, for a given $\alpha \geq \beta/2 -  (1/2+p+\gamma)$, where $\alpha >0$,  the value of $\tau_{\epsilon}$ that leads to the minimax rate is such that the cutoff level $i_{\epsilon } \asymp [\epsilon^{-2}]^{\frac{1}{1 + 2\beta + 2(p+\gamma)}}$ is independent of $\alpha$ and is the cutoff level corresponding to the minimax optimal projection estimator (projecting on the first $i_{\epsilon}$ components). 
\end{remark}
 \citet{FlorensSimoni} consider the case of our setup with $p+\gamma>0$ and  $|\mu_{0,i}| \asymp i^{-\beta-1/2}$.
 In \citet{agapiou_posterior_2013}, the rate of contraction in this setting  is given only for $\beta > \alpha+1/2$ (in our notation) by
\begin{equation}
\epsilon^{\frac{\beta \wedge (p+\gamma+2\alpha+1)}{\delta(1+2\alpha) + 1 + 2p +2\gamma + 2[\beta \wedge (p+\gamma+2\alpha+1)]}},\quad \forall \delta >0,
\end{equation}
(where $\gamma_A = \beta/(\alpha+1/2)$, $\ell_A = p/(2\alpha+1)$, $s_{0,A}=1/(2\alpha+1)$, $\Delta_A = (p+\gamma)/(\alpha+1/2)+1$, with subscript $A$ referring to parameters in \citet{agapiou_posterior_2013}). In particular, the authors' assumption $\Delta_A >2s_{0,A}$ is equivalent to assumption $p+\gamma+\alpha+1/2 >1$ in our notation which is stronger than our assumption $p+\gamma+\alpha+1/2 >0$. In fact, under the latter assumption, it is not possible to achieve the minimax optimal contraction rate for $p+\gamma>-1/2$ with constant $\tau_\epsilon$, which recall is achieved when $\alpha=\beta$. Also, the authors refer to the case $p+\gamma <0$ as self-regularising, with no regularisation being necessary which we do not find in case  $p+\gamma \in (-1/2,0)$; also their rate in this case
is not faster unless $\delta(1/2+\alpha)=-( 1/2 + p +\gamma) -0.5\beta \wedge (p+\gamma+2\alpha+1)$ which is possible only if $\alpha <-1.5  (p+\gamma+1/2)-1/4$ and $\beta< -2( 1/2 + p +\gamma)$.

\begin{example}
Consider case $V = (K K^T)^a$.  For mildly ill-posed inverse problems where $\kappa_i \asymp i^{-p}$, where $p>0$, we have $\gamma = - p a$ and $p+\gamma = p(1-a)$. Take  the prior with $\lambda_i = \tau_{\epsilon}^2 i^{-2\alpha-1}$ and the parameters as specified in Corollary~\ref{Cor:OptPriorParam}, so that the posterior contraction rate coincides with the minimax rate. Then, the rate is $\epsilon \sqrt{\log 1/\epsilon}$, for $a \in (0, 1+0.5/p)$ the rate is $[\epsilon^{2}]^{ \beta/(1 + 2\beta + 2 p(1-a))}$ and for $a > 1+0.5/p$ the rate is $\epsilon$. Note that the model $Y = K(\mu + \epsilon W)$ with white noise $W$ corresponds to $a=1$, where the rate of posterior contraction coincides with the minimax rate of the direct problem.
\end{example}

\subsection{Contraction rate with plugged-in $(\sigma_i)$}

In this section, we investigate how a plug-in estimator of $(\sigma_i)$ affects the rate of contraction 
 for geometric sequences.

\begin{proposition}
\label{prop:plugincontractionratesimproved}
Consider the inverse problem (\ref{f2model}) formulated in the sequence space under Assumption~\ref{assumptionB1} with noise $W$ satisfying \eqref{eq:WnoiseDecomp} and Assumption~\ref{AssumeEVVar} with $\gamma< 0$. Consider the prior distribution (\ref{eq:SeqPrior}) under Assumption~\ref{AssumeEpriorV} and the true function $\mu_0$ satisfying Assumptions~\ref{smoothGeneral} and \ref{AssumeSmoothPoly}. Let Assumption~\ref{AssumeDepPoly} hold.  

Suppose that $(\sigma_i)$ are unknown but their estimators $(\hat\sigma_i)$ are available.
\begin{enumerate}
\item Assume that estimator $(\hat\sigma_i)$  satisfies Assumption \ref{AssumeEVVarPlugIn} (absolute bound).
\begin{enumerate}
\item If $\epsilon_\sigma < C \left[\epsilon^{2}\tau_\epsilon^{-2}  \right]^{-\gamma/(\alpha+1/2+p+\gamma)}$ then the rate of contraction is not affected by using a plug-in estimator of $(\sigma_i)$, i.e. it coincides with the rate given in Theorem~\ref{contractiontheorem2}, up to a constant.

\item If $\epsilon_\sigma \geq C\left[\epsilon^{2}\tau_\epsilon^{-2}  \right]^{-\gamma/(\alpha+1/2+p+\gamma)}$, then
 the contraction rate  of the posterior distribution is given by
\begin{eqnarray*}
{\rate}_{plugin}^2 = \epsilon^{2} \left(\log \epsilon_\sigma^{-1}\right)^{\mathbb{I}\{1+2(p+\gamma)=0\} }
+\tau_{\epsilon}^2 \left[\epsilon^{-1}_{\sigma}\epsilon^{-2}\tau_\epsilon^2  \right]^{-2\alpha/(1+2\alpha+2p)}
+\left[\epsilon^{-1}_{\sigma}\epsilon^{-2}\tau_\epsilon^2  \right]^{-2\beta/(1+2\alpha+2p)}
+ \epsilon^4 \tau_{\epsilon}^{-4} \epsilon_{\sigma}^{\frac{(1+2\alpha+2(p+\gamma)-\beta)_+}{\gamma}}.
\end{eqnarray*}
\end{enumerate}

\item Assume that estimator $(\hat\sigma_i)$  satisfies Assumption \ref{AssumeEVVarPlugInRelative} (relative bound).
The contraction rate  of the posterior distribution is given by
\begin{eqnarray*}
 {\rate}_{plugin}^2  &=&
  \epsilon^{2} \left[\min(N, i_{\epsilon})\right]^{(2p+2\gamma+1)_+} [\log \min(N, i_{\epsilon})]^{I(2p+2\gamma+1=0)}
 +   \tau_\epsilon \left[\tau_\epsilon \epsilon^{2}\right]^{2\alpha/(2\gamma +2p +2\alpha +1)}  I (i_{\epsilon}\leq N)\\
  &+& \min(N,i_{\epsilon})^{-2\beta}
+    \epsilon^4 \tau_\epsilon^2  [\min(N,   i_{\epsilon})]^{2(2p+2\gamma+2\alpha+1-\beta)_+}\\
&+&  \tau_{\epsilon} \epsilon^{-2}  \max\left( N^{-2\alpha -2p-2\gamma}, i_{\epsilon}^{ -2\alpha -2p-2\gamma}\right) I(i_{\epsilon}<N) [\log N]^{I (\alpha +p+\gamma=0)}
\end{eqnarray*}
with $i_{\epsilon,\pm} \asymp i_{\epsilon} = C \left[\tau_\epsilon \epsilon^{-2}\right]^{1/(2\gamma +2p +2\alpha +1)}$.
 If $N \geq C i_{\epsilon}$ for some positive $C>0$, then the rate is the same as if $\sigma_i$'s were known; if $N / i_{\epsilon}\to 0$ then the rate is affected by the plug-in.
\end{enumerate}
\end{proposition}
In the case the bound is absolute, if the error $\epsilon_\sigma$ of estimating $(\sigma_i^2)$ is small enough, then the rate of estimation of $\mu$ is not affected. In the case of relative bound, as long as $N\geq i_{\epsilon}$, the rate of estimation is not affected. 

\subsection{Repeated observations}


\subsubsection{Posterior contraction rate with plugged in sample variances}

Now we assume that we have $m$ repeated observations from model\eqref{f2model} as described in Section~\ref{sec:RepeatedObs}, and that we plug in  $\hat\sigma^2_i = s^2_i$ for $i=1,\ldots, N$. Then, under Assumptions~\ref{AssumeEVVar}, \ref{AssumeEpriorV},  \ref{AssumeSmoothPoly} and  \ref{AssumeDepPoly}, Proposition~\ref{prop:plugincontractionratesimproved} states that
 as long as $\log N= o(m)$ and $N \geq C \left[\tau_\epsilon^2 \epsilon^{-2}\right]^{1/(2\gamma +2p +2\alpha +1)}$, the posterior rate of contraction is the same as if we used the true values of $\sigma_i^2$. Such $N$ exists if $\log \left(\tau_\epsilon^2 \epsilon^{-2}\right) =o(m)$, which is generally a mild constraint.

We have also applied Theorem~\ref{th:GeneralPlugin} to the model with repeated observations using additive error (Assumption~\ref{AssumeEVVarPlugIn}) however in that setting the plug-in effect can lead to a sub-optimal rate for some $\beta$.

We will discuss the choice of $N$ in this setting in practice using the empirical Bayes approach.




\subsubsection{Estimation of the eigenfunctions of $V$}\label{sec:EstimateEigenFunctions}

Here we assume that $e_{1,i} = \phi_i$, i.e. the eigenfunctions of operator $V$, and discuss their estimation using repeated observations. We apply the results of \citet{KoltchinskiiLounici} to evaluate the number of eigenfunctions of $V$ that are possible to estimate consistently, and hence verify whether the posterior distribution can achieve the optimal rate in the minimax sense for repeated observations where the eigenfunctions are unknown. Their conditions assume that the eigenvalues are decreasing, and rely on $trace(V)$ being finite, so for $\sigma^2_i \asymp i^{2\gamma}$  these conditions hold only if $\gamma < -1/2$ which is the case we consider here.

\begin{lemma}\label{lem:EstEigenFunctions}
Recall that $V = \sum_{i=1}^{\infty} \sigma^2_i \phi_i \phi_i^T$ and denote the projection matrices by $P_i = \phi_i \phi_i^T$.
 In the repeated observations model (\ref{eq: prob model:rep obvs}), denote the sample estimator of the covariance matrix by $\hat V$ and the corresponding projection matrices by $\hat P_i$.

Assume that $\sigma^2_i \asymp i^{2\gamma}$ for some $\gamma < -1/2$.
Then, $\hat P_r$ are consistent estimators of  $P_r$ simultaneously for $r=1,\ldots, N$ with high probability if
  $N = o\left(m^{ 1/(1/2-4\gamma)}\right)$.
\end{lemma}
Note that for a fractional Brownian motion with Hurst exponent $H$, $\gamma = 1/2-H\in (-1/2, 1/4)$ hence this theorem does not apply to it (see below for a discussion of the associated known series expansions).
 In other settings, such as a nonparametric regression with a finite-dimensional approximation of $V$ it may be possible to estimate the eigenfunctions for other values of $\gamma$ but this problem is beyond the scope of this paper.

Recall that for the posterior distribution to contract at the optimal rate, we need $N\geq i_\epsilon = \left(\epsilon^{-2}\tau_{\epsilon}^2\right)^{1/(2p+2\gamma+2\alpha+1)}$. 
For $\alpha \geq \max(\beta-1/2, \beta/2 -  (1/2+p+\gamma))$, $p+\gamma+1/2 +\alpha >0$, $\gamma < -1/2$ and $\tau_{\epsilon}$ given in Corollary~\ref{Cor:OptPriorParam}, the condition of the lemma holds if
$$
m \gg \left[\epsilon^{-2}\right]^{(1/2-4\gamma)/(2p+2\gamma+2\beta+1)},
$$
i.e. under this condition estimators of the first $N\geq i_{\epsilon}$ eigenvectors are consistent and the corresponding posterior contraction rate is optimal.


\subsection{Estimated $\gamma$}\label{sec:PolyEstGamma}

Assume  that for large $i$, $\sigma_i = C i^{\gamma}(1+o(1))$. Therefore, we can consider $\sigma_i = C i^{\gamma}$ for all $i\geq i_0$ for some $i_0\geq 1$. Suppose $\gamma$ and $C$ are estimated independently of $(Y_i)_{i=1}^n$, and that with high probability $\hat C/C-1\in [-c_{\max}, c_{\max}]$ and $\hat\gamma - \gamma \in [-z_{\max}, z_{\max}]$ for some   $c_{\max}, \, z_{\max}>0$.

This implies that, with the same high probability,  for all $i=i_0,\ldots, N$,
\begin{eqnarray*}
|\hat\sigma_i^2/\sigma_i^2 -1| &\leq&  \max\left[c_{\max} N^{2 z_{\max}} -1, 1- c_{\max}^{-1} N^{-2 z_{\max}} \right] =: \epsilon_{\sigma}.
 \end{eqnarray*}

If $c_{\max} \to 1$  and $ z_{\max} \log N\to 0$ then the upper bound tends to 0 for all $\gamma\in \mathbb{R}$, i.e. Assumption~\ref{AssumeEVVarPlugInRelative} holds, as long as we have an alternative consistent estimator of $\sigma_i$ for $i < i_0$ or if $i_0=1$. Such bound holds also for the absolute bound (Assumption~\ref{AssumeEVVarPlugIn}) however we found that the conditions for the corresponding plug in  posterior distribution to contract at the optimal rate are stronger, so we do not use it here.

If $c_{\max}$ and $z_{\max}$ are known, then choosing $N \geq c_{\max} \hat{C} \left[\tau_\epsilon \epsilon^{-2}\right]^{1/(2[\hat\gamma-z_{\max}] +2p +2\alpha +1)}$ implies that the rate of contraction of the posterior of $\mu$ is not affected, with high probability. Now we investigate what the expressions for $c_{\max}$ and $z_{\max}$ are under repeated observations.

 \begin{example} In the setting of repeated observations and $\sigma_i^2 = C i^{2\gamma}$ for $i\geq i_0$, we can use $\hat C = s_{i_0}^2 i^{-2\hat\gamma}$ and  for some subset $I_N \subseteq \{i_0,i_0+1,\ldots, N\}$,
 $$\hat\gamma = \frac 1 {2|I_N|} \sum_{i\in I_N}  \frac{\log \left(s_i^2/s_{i_0}^2\right)}{\log (i/i_0)}. $$
Now we derive confidence intervals for $C$ and $\gamma$ when $m$ is large, and hence the expression for relative bound $\epsilon_{\sigma}$.
 The ratio $\frac{s_i^2}{  [i/i_0]^{2\gamma} s_{i_0}^2}$ has distribution $F(m,m)$ hence, approximately for large $m$,
 $ \log \left[s_i^2/( s_{i_0}^2)\right] - 2\gamma\log (i/i_0) \sim N(0, 1/m)$,
 implying that
 $$
 \hat \gamma - \gamma \sim N\left(0, \frac 1 {4 m |I_N|} \sum_{i\in I_N}  \frac 1 {[\log (i/i_0)]^2} \right).
 $$
  This suggests that averaging over a small number of $i$ close to $N$ (or even a single $i=N$) leads to an estimator with the smallest variance of order $[m \log N]^{-1}$, e.g. $I_N = \{N\}$ or   $I_N= \{(N-k) :  N\}$ for a small $k$, for instance $k=5$, provides a more robust estimator.

 For the case $I_N = \{N\}$,
 \begin{equation}\label{eq:hatgammaN}
 \hat\gamma = \frac 1 {2} \frac{\log \left(s_N^2/s_{i_0}^2\right)}{\log (N/i_0)},
\end{equation}
which asymptotically follows distribution $N\left(\gamma, \left[4 m (\log (N/i_0))^2\right]^{-1}\right)$, and hence a $(1-\alpha)100\%$ confidence interval for $\gamma$ is
$$
\left[\hat\gamma-\frac{z_{\alpha/2}}{2\sqrt{m} \log (N/i_0)}, \hat\gamma+\frac{z_{\alpha/2}}{2\sqrt{m} \log (N/i_0)} \right].
$$ Here $z_{a}$ satisfies $  1-\Phi(z_a)=a$, where $\Phi(z) = P(N(0,1) < z)$.

  For large $m$, the distribution of $m s_{i_0}^{2}/\sigma_{i_0}^2 = m \hat C/ C$, $\chi^2_m$, can be approximated by $N(m, 2m)$, and hence for  known $\gamma$ with probability $1-\alpha$,  a $(1-\alpha)100\%$ confidence interval for $C$ is defined by $\sqrt{m/2} (s_{i_0}^2/C i_0^{2\gamma}-1) \in [-z_{\alpha/2}, z_{\alpha/2}]$. 
 Using the above $(1-\alpha)100\%$ confidence interval for $\gamma$ based on (\ref{eq:hatgammaN}) gives the following asymptotic $(1-2\alpha)100\%$ confidence interval for $C$ when $m$ is large:
 $$
\frac{s_{i_0}^2  i_0^{-2\hat\gamma - z_{\alpha/2}/[\sqrt{m} \log (N/i_0)]} }{1+z_{\alpha/2}/\sqrt{m/2}} \leq   C   \leq \frac{s_{i_0}^2 i_0^{-2\hat\gamma+z_{\alpha/2}/[\sqrt{m} \log (N/i_0)]} }{1 -  z_{\alpha/2}/\sqrt{m/2}}.
 $$
 Note that if $i_0=1$, then  a $(1-\alpha)100\%$ confidence interval for $C$ is $\sqrt{m/2} (s_{1}^2/C -1) \in [-z_{\alpha/2}, z_{\alpha/2}]$.




 In particular, this implies that
  with probability $1-2\alpha$, $|\hat C/C - 1| \leq z_{\alpha/2} m^{-1/2}\left[\sqrt{2} + \log(i_0)/\log (N/i_0)\right]=:c_{\max}$ and $z_{\max}= z_{\alpha/2}/\left[\sqrt{m} \log (N/i_0)\right]$ and hence, for large $m$, with probability $1-2\alpha$,
 the relative upper bound for large $m$  is approximately
  $\epsilon_{\sigma} 
 =  m^{-1/2} z_{\alpha/2} [2 + \sqrt{2}]$.
 \end{example}

\subsection{Error in operator }

Now we consider the case when the operator $K$ may be observed with error, and study its effect on the posterior contraction rate for mildly ill-posed inverse  problems (see Section~\ref{sec:ErrorOperatorGeneral} for the setup). 

Given white noise,  Theorem~\ref{th:contractiontheoremGeneralPlugIn2} and Lemma~\ref{lem:GeneralKwithError} imply that for $\mu_0 \in S^\beta(A)$, the  posterior distribution of $\mu$ given data $y_1,\ldots, y_N$ and observed $\hat \kappa_1,\ldots, \hat \kappa_N$ contracts at the optimal rate for $\alpha$ and $\tau$ specified in Corollary~\ref{Cor:OptPriorParam} if  $N\geq i_{\epsilon} \asymp \epsilon^{-2/(1+2\beta+2p)}$ and   $\delta [\log N]^{1/2} N^{p}  =o(1)$ (condition (\ref{cond:ErrorOperator})). There exists $N$ satisfying both conditions if
 $$
 \delta =o\left( \epsilon^{2p/(1+2\beta+2p)} [\log 1/(\epsilon)]^{-1/2} \right),
 $$
 for small $\epsilon$. This holds, for instance, if $\delta \leq \epsilon$.

Another problem is to determine the rate of contraction of the posterior distribution with noisy $\hat{\kappa}_i$ for given values of $\epsilon$ and $\delta$  and compare it to the minimax rate for linear estimators in one dimension, $[\max(\epsilon, \delta)]^{-2\beta/(2\beta+2p+1)}$ (Section~5.2 of  \citet{HoffmannReiss}). For mildly ill-posed problems and white noise,  Theorem~\ref{th:contractiontheoremGeneralPlugIn2} applied to model (\ref{eq:SeqErrorOperator}) with $\tau_\epsilon$ and $\alpha$ specified below in Corollary~\ref{Cor:OptPriorParam} gives the squared rate
\begin{eqnarray*}
\rate^2  =
  \epsilon^{2}  [\min(N, i_{\epsilon})]^{1+2p}
 +  [\epsilon^{-2}]^{-2\beta/(2p+2\beta+1)}
+  \tau_{\epsilon}^{-4}  \epsilon^4   [\min(N,   i_{\epsilon})]^{2(2p+2\alpha+1-\beta)_+}
+ \min(N,i_{\epsilon})^{-2\beta},
\end{eqnarray*}
where $i_{\epsilon} \asymp \left[\epsilon^{-2}\right]^{1/(2p+2\beta+1)}$. Assume that $2p+2\alpha+1-\beta \leq 0$ which is necessary to obtain the optimal contraction  rate with known $K$.
 Taking $N = \delta^{-2/(2p+2\beta+1)}$ satisfies the required condition
$$
\delta N^{p} [\log N]^{1/2} = \delta^{1-2p/(2p+2\beta+1)} [\log (1/\delta)]^{1/2} = o(1), \text{ as } \delta \to 0,
$$
and leads to the optimal contraction rate
$\rate  = [\max(\delta, \epsilon)]^{2\beta/(2p+2\beta+1)}$.

We can easily generalise this to correlated noise with $\sigma_i \asymp i^\gamma$ and known $\gamma$, as it does not affect condition  (\ref{cond:ErrorOperator}).
\begin{lemma}\label{lem:MildIPerrorOp}
Consider the inverse problem (\ref{f2model}) formulated in the sequence space under Assumption~\ref{assumptionB1} with noise $W$ satisfying \eqref{eq:WnoiseDecomp} and Assumption~\ref{AssumeEVVar} with known $\gamma\leq  0$,  and the true function $\mu_0$ satisfying Assumptions~\ref{smoothGeneral} and \ref{AssumeSmoothPoly}. Let Assumption~\ref{AssumeDepPoly} hold.

Assume that operator $K$ is observed with error (\ref{eq:ErrorOperator}), and prior distribution (\ref{eq:SeqPrior}) satisfies Assumption~\ref{AssumeEpriorV} with $\lambda_i =0$ for $i>N$.

Then, for $\mu_0 \in S^{\beta}(A)$,  $\tau_\epsilon$ and $\alpha$ specified in Corollary~\ref{Cor:OptPriorParam} with $2p+2\alpha+1-\beta \leq 0$, the posterior contraction rate of $\mu$ given $(y_1, y_2,\ldots)$ and $(\hat \kappa_i)_{i=1}^N$ with $N$ satisfying (\ref{cond:ErrorOperator}) is given by
\begin{eqnarray*}
\rate  =
  \epsilon   [\min(N, i_{\epsilon})]^{(1/2+p+\gamma)_+} [\log \min(N, i_{\epsilon})]^{0.5 I(1/2+p+\gamma=0)}
 +  \epsilon^{2\beta/(2\beta+(1+2p+2\gamma)_+)}
+ [\min(N,i_{\epsilon})]^{-\beta},
\end{eqnarray*}
where $i_{\epsilon} \asymp [\epsilon^{-2}]^{1/(2\beta+(2p+2\gamma+1)_+)}$. The rate is not affected by the error in the operator if  $\delta =o\left( \epsilon^{2p/( 2\beta+(2p+2\gamma+1)_+)} [\log 1/(\epsilon)]^{-1/2} \right)$ for small $\epsilon$.
Moreover, if $N = \delta^{-2/(2\beta+(2p+2\gamma+1)_+)}$, then
\begin{eqnarray*}
\rate  =  [\max( \epsilon, \delta)]^{2\beta/(2\beta+(2p+2\gamma+1)_+)} \left[\log (1/\max( \epsilon, \delta))\right]^{0.5 I(1/2+p+\gamma=0)}.
\end{eqnarray*}
\end{lemma}
The proof is a straightforward extension of the above argument to the case $\gamma \neq 0$ and is omitted.

\subsection{Empirical Bayes posterior distribution of $\mu$}\label{sec:EBpoly}

As we have seen in Section~\ref{sec:PostContractionPoly}, the posterior distribution of $\mu$ contracts at the optimal rate (in the minimax sense) for a particular set of the prior parameters $\alpha$ and $\tau$ (Corollary~\ref{Cor:OptPriorParam}). In this section  we study the contraction rate of the posterior distribution of $\mu$ with a plugged in estimator of the prior scale $\tau$ that does not rely on knowing $\beta$, true smoothness of $\mu_0$.

Consider the prior Gaussian distribution with $\lambda_i = \tau \lambda_{0,i}$, with empirical Bayes posterior of $\mu$ using the maximum of marginal likelihood estimator of $\tau$ (MMLE) and fixed $\lambda_{0,i}$:
\begin{eqnarray}\label{def:EBtau}
\hat\tau = \arg \max_{\tau >0} p({\bf y} \mid \tau)
 = \arg \min_{\tau >0} \sum_{i=1}^\infty \left[\frac{y_i^2}{\kappa_i^2 \lambda_{0,i}\tau + \epsilon^2 \sigma_i^2} +\log(\kappa_i^2 \lambda_{0,i}\sigma_i^{-2}\tau \epsilon^{-2} +1)\right],
\end{eqnarray}
where $p({\bf y} \mid \tau)=\int p({\bf y} \mid \mu) d P(\mu\mid \tau)$ is the marginal density of ${\bf y}=(y_1,y_2,\ldots)$ with respect to measure $\prod_{i=1}^{\infty} N(0, \epsilon^2 \sigma_i^2)$. Recall that $\tau_{\epsilon}=\tau^{1/2}$.



The following assumption is necessary to show convergence of $\hat\tau$ in probability (but not for optimality of the posterior).
\begin{assumption}\label{Assume:SelfSim}
For $\mu_0 \in S^\beta(A)$, assume that  for any $\beta' \geq \beta$ there exists $N_0\geq 1$ and $c_0>0$ such that for any $N\geq N_0$, $N^{-2(\beta'-\beta)}\sum_{i=1}^N \mu_{0i}^2 i^{2\beta'} \geq c_0$.
\end{assumption}

Now we prove that  the posterior distribution of $\mu$ with plugged in value of $\hat\tau$ defined by (\ref{def:EBtau}) contracts at the optimal rate adaptively, uniformly over $\mu_0 \in S^\beta(A)$ with $0<\beta \leq B_0 <\infty$, in the minimax sense, under some conditions on prior smoothness $\alpha$.
\begin{theorem}\label{th:EBtau1}
Consider the inverse problem (\ref{f2model}) formulated in the sequence space under Assumption~\ref{assumptionB1} with noise $W$ satisfying \eqref{eq:WnoiseDecomp} and Assumption~\ref{AssumeEVVar}, with prior distribution (\ref{eq:SeqPrior}) under Assumption~\ref{AssumeEpriorV} with $\tau = \tau_\epsilon^2$, and the true function $\mu_0$ satisfying Assumptions~\ref{smoothGeneral} and \ref{AssumeSmoothPoly}. 
 Assume that  $1/2+\beta+p+\gamma>0$, and $1/2+\alpha+p+\gamma>0$.
\begin{enumerate}
\item If $\mu_0 \not=0$, then there exist $\underline{\tau}_{\epsilon}$ and $\overline{\tau}_{\epsilon}$ such that $P(\underline{\tau}_{\epsilon} < \hat{\tau}_{\epsilon} < \overline{\tau}_{\epsilon}) \rightarrow 1$ as $\epsilon \to 0$. The posterior distribution of $\mu$ with plugged in $\hat \tau_{\epsilon}$ is consistent.

\item 	If $\mu_0 \not=0$ and $\alpha \geq \max(\beta-0.5, \beta/2 - 1/2 - p -\gamma)$, then the  posterior contraction rate with plugged in $\hat{\tau}_{\epsilon}$ is optimal in the minimax sense.

\item If $\mu_0 \not=0$,  $\alpha \geq \beta-0.5$ and Assumption~\ref{Assume:SelfSim} holds, then
\begin{enumerate}
\item $\underline{\tau}_{\epsilon}$ and $\overline{\tau}_{\epsilon}$ are of the same order,
\item   $ \hat{\tau}_{\epsilon}/\tau^\star_{\epsilon}  \to 1$   in probability as $\epsilon \to 0$ where $\tau^\star_{\epsilon} = \arg \max_{\tau_{\epsilon}>0} \mathbb{E}_{\mu_0} \log p( {\bf y} \mid \tau_{\epsilon}) \asymp [\epsilon^2]^{2(\alpha-\beta)/(1+2\beta+2p+2\gamma)}$.
\end{enumerate}
\item Finally, if $\mu_0 =0$, then $\hat{\tau}_{\epsilon} = \mathcal{O}_P(\epsilon)$ and the plug-in contraction rate is $\rate = \epsilon^2$.
\end{enumerate}

\end{theorem}

We show in the proof that for $\mu_0 \not=0$,
\begin{eqnarray}\label{eq:AsympTau}
\underline{\tau}_\epsilon^2 &\gtrsim& \epsilon^{\frac{1}{(1+\alpha+p+\gamma)}}, \\
	\overline{\tau}_{\epsilon}^2 &\lesssim&
 \begin{cases}
\epsilon^{-4(\alpha-\beta)/(1 + 2p +2\gamma + 2\beta) } (1+o_P(1)), \quad & \text{ if $ \alpha   +1/2 \geq  \beta$},\\
\epsilon^{1/(1+ p +\gamma +\alpha)} (1+o_P(1)), \quad & \text{ if  $\alpha   +1/2 <  \beta$}
\end{cases}
\end{eqnarray}

\begin{corollary}\label{cor:EBtau}
Under the conditions of Theorem~\ref{th:EBtau1}, the posterior distribution of $\mu$ with the plugged in estimator $\hat\tau_\epsilon$ is consistent. The contraction rate is optimal in the minimax sense  adaptively over $S^\beta(A)$ for $0<\beta \leq B_0<\infty$ as long as
\begin{equation}\label{eq:EBalpha}
\alpha \geq \max(B_0-1/2, B_0/2 -1/2-p-\gamma),  \text{ and } \alpha > (- p - \gamma - 1/2)_+.
\end{equation}

\end{corollary}

The theorem states that if the prior smoothness parameter $\alpha$ is chosen large enough, then using the MMLE  $\hat\tau$ allows us to achieve the optimal rate of contraction. If $p+\gamma=0$ and $\mu_0 \in H^\beta(A)$, the results coincide with those of \citet{SzaboEmpiricalBayes}. We have weakened the condition for convergence of $\hat\tau_\epsilon$ in probability and show that it converges to the ``oracle'' value of $\tau_\epsilon$.
 The additional condition $ \alpha \geq  \beta/2 - 1/2 - p - \gamma$ arises only if $p+\gamma <0$ and comes from the saturation term.

Now we consider the case of repeated observations where $\sigma^2_i$ are estimated, and we discuss the choice of $N$ in an adaptive setting, where we plug in the MMLE $\hat \tau$. Recall that the rate is not affected by the plugged in estimator $(\hat \sigma_i^2)$ if $\log(N) = o(m)$ and $ N \geq i_{\epsilon_+}$ (see Proposition~\ref{prop: rep obvs M contraction2}).

\begin{remark}
For a mildly ill-posed inverse problem with geometric sequence $(\sigma_i^2)$ satisfying Assumption~\ref{AssumeEVVar}, for the repeated observations and estimated $\sigma_i^2$,  the rate of contraction of the posterior distribution of $\mu$ with the plugged in $(\hat\sigma_i^2)_{i=1}^N$ is not affected if $\log (1/\epsilon) =o(m)$,   which is generally a mild constraint. In particular, we can choose $N\geq  C  \left[\hat\tau_\epsilon \epsilon^{-2}\right]^{1/(1 + 2p +2\gamma +   2\alpha)}$ for some constant $C>0$,  which leads to a consistent posterior distribution due to the lower bound on $\hat\tau_\epsilon$ given by (\ref{eq:AsympTau}) (with high probability), and for $\mu_0$ satisfying Assumption \ref{Assume:SelfSim} the posterior contraction rate is optimal.
\end{remark}
%


\section{Simulation}\label{sec:SimData}

\subsection{Simulation set up}

We illustrate the theoretical results obtained in Section~\ref{sec:PolyEV} on indirect observations from model (\ref{f2model}) corrupted by the Volterra operator \cite{halmos_hilbert_1974} with  dependent Gaussian noise and conjugate prior (\ref{eq:priorOp}). Here $H_1=H_2=L^2[0,1]$.

The Volterra operator, $K: \, L^2[0,1] \rightarrow L^2[0,1] $ is defined by
\[
K\mu(x):=\int_{0}^{x}\mu(s)ds, \text{~~and~~} K^T\mu(x):=\int_{x}^{1}\mu(s)ds.
\]
The eigenvalues of $KK^T$ and the orthonormal eigenbasis for the range of $K$ are
\[
\kappa_i:=\left[\left(i-\frac{1}{2}\right)^2\pi^2\right]^{-\frac{1}{2}},
\text{~~and~~} \phi_i(x):= \sqrt{2}\sin\left(\left(i-\frac{1}{2}\right)\pi x\right) \text{~~for every~~} i \in \mathbb{N},
\]
where $\kappa_i\asymp i^{-p}$ with $p=1$. The corresponding orthonormal eigenbasis of $K^T K$ is $e_i(x)=\sqrt{2}\cos\left(\left(i-\frac{1}{2}\right)\pi x\right)$.

We will estimate the following approximation of $\mu_0(x)$:
$\mu^N_0(x) := \sum \limits_{i}^N \mu_{{0},i}e_i(x)$,
 where $N$ is the truncation parameter and is large to ensure good approximation: $N \geq  \max\left(50,\epsilon^{-2/(1+2p)}\right)$. We consider a particular function with
 $\mu_{{0},i}:=i^{-3/2}\sin(i)$ which belongs to  $ S^\beta$ with $\beta=1$ (this can be shown using Dirichlet's test, see e.g. \citet{voxman_advanced_1981}).

 We consider Gaussian noise $W$ with truncated covariance operator $V^N = \sum_{i=1}^N \sigma_i^2 \phi_i \phi_i^T$ and $\sigma_i= i^\gamma$.
 Realisations of the data are simulated as follows:
 $Y_i \sim N(\mu_{0,i}\kappa_i,\epsilon^{2}\sigma_i^2)$,
 independently for $i=1,\ldots, N$.

We consider a centered Gaussian prior distribution with  $\lambda_i=\tau_\epsilon^2\,i^{-1-2\alpha}$ and different values of $\alpha$ and $\tau_\epsilon$, including the MMLE $\hat\tau_\epsilon$.

The corresponding posterior distribution is
\[
\hat{\mu}^N |Y \sim N\left(\sum \limits_{i}^N\frac{Y_i\kappa_i\lambda_ie_i }{\lambda_i\kappa_i^2 + \epsilon^{2}\sigma^2_i}, \sum \limits_{i}^N\frac{\epsilon^{2}\sigma^2_i\lambda_i e_i^2 }{\lambda_i\kappa_i^2 + \epsilon^{2}\sigma^2_i}\right).
\]

\begin{figure}
\includegraphics[width=0.98\textwidth, height = 1.6in]{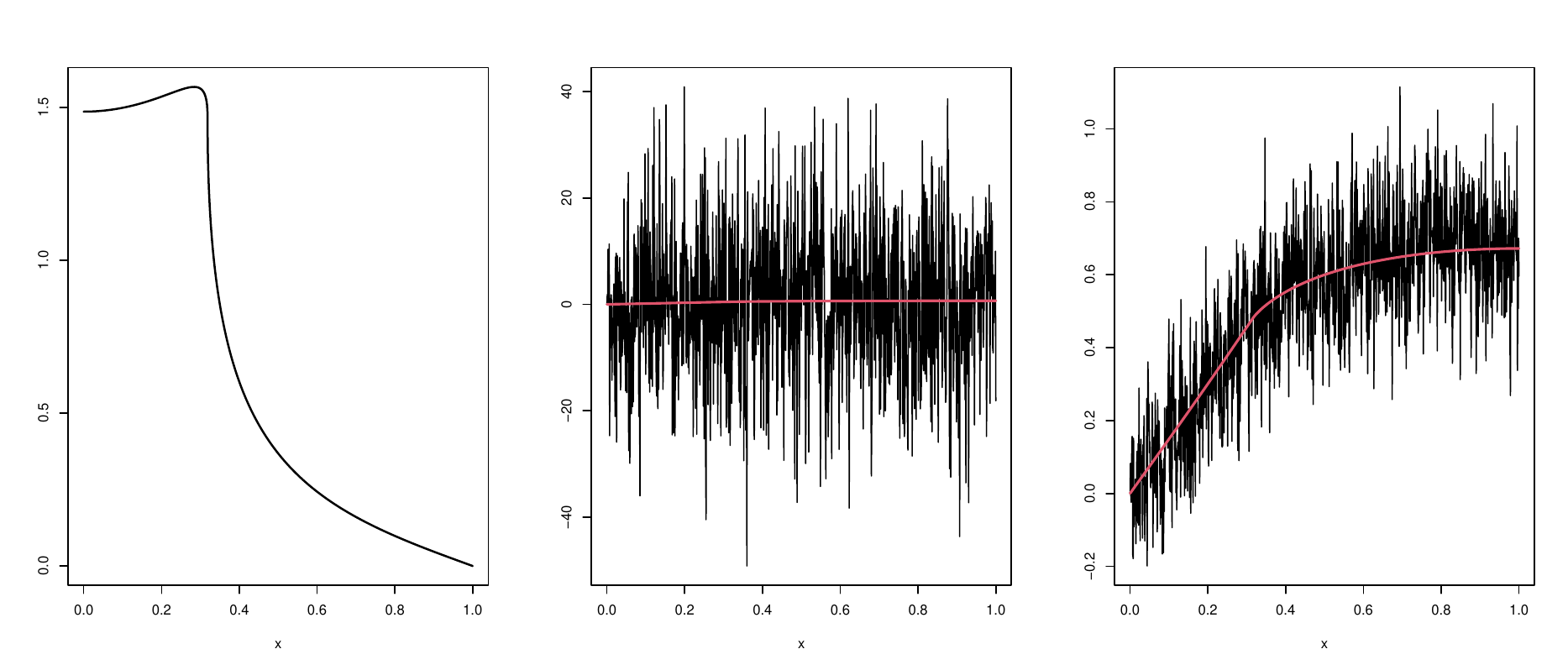}
\caption{Graph of the truncated true function $\mu^N_0(x)$ (left). Noisy data $Y(x)$ with $\epsilon=10^{-2}$  (centre),  $\epsilon=10^{-4}$ (right), with the same indirectly observed noiseless function $Y_0(x) = (K \mu^N_0) (x)$ (red line); $N = 2000$ and $\gamma = 0.5$.} \label{dataplot}
\end{figure}

The (truncated) true function $\mu^N_0(x)$, along with its noisy realisations $Y(x)$ for two different noise levels are shown in Figure~\ref{dataplot}. We can see that with the noise level $\epsilon=10^{-2}$, the observed function $Y$ is very noisy compared to the signal $Y_0=K\mu_0$: the range of observations is approximately $[-44, 43]$ whereas the range of the signal is $[0,0.67]$. For a smaller noise level $\epsilon=10^{-4}$, the range of observations is $[0.2, 0.8]$ which is comparable to the range of $Y_0$.
 \selectlanguage{english}

A key property we want to study here is the variability of the posterior distribution around the true value of the function and posterior coverage, i.e. the posterior probability that the true function lies in the support of the posterior. A common way to investigate the posterior support in Bayesian nonparametric models is to plot a large number of draws from the posterior distribution, where the ``centre'' of the support is displayed using the posterior mean. We plot 100 draws as plotting a larger number in the considered examples does not much change the posterior support. Our main interest is to see how the coverage of the true function by the posterior distribution and its contraction is affected by the choice of prior smoothness $\alpha$, both for a fixed prior scale $\tau$ (non-adaptive case) and for the empirical Bayes $\hat\tau$, as well as by the variance parameters $\epsilon$ and $\gamma$.


\subsection{Non-adaptive posterior distribution of $\mu$ }

\begin{figure}
\centering
\includegraphics[width=1\textwidth,height=4in]{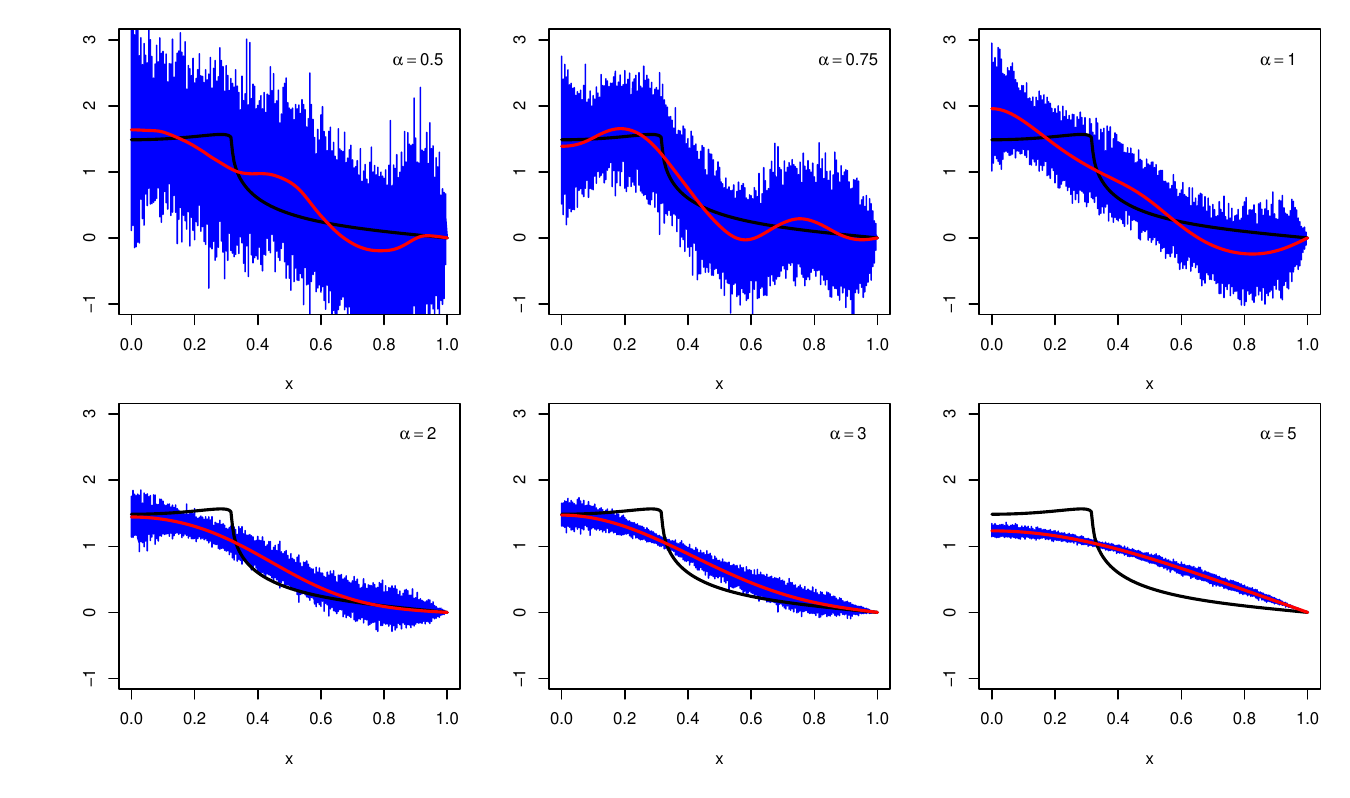}
\caption{{Plots of $\mu^N_0(x)$ (black line) along with the posterior mean (red line), and $100$ draws from the posterior (blue dashes) for $\alpha=(0.5,0.75,1,2,3,5)$ respectively, with $\epsilon=10^{-2}$, $\gamma=0.5$ and $N=2000$ in all cases.}} \label{crediblebandsalpha1}
\end{figure}

Firstly we study how posterior variability and the coverage of the true function varies with prior smoothness $\alpha$ as we fix prior scale $\tau=1$. In this section we also fix $\epsilon=10^{-2}$, $\gamma=0.5$,  $N=2000$, $\tau_\epsilon=1$ and consider how different values of a priori smoothness $\alpha$ affect behaviour of the posterior distribution. We consider the following values of $\alpha=(0.5,0.75,1,2,3,5)$. Draws from the posterior distributions of $\mu^N(x)$ corresponding to these values of $\alpha$ are given in Figure~\ref{crediblebandsalpha1}, for the same realization of $Y(x)$. Individual draws from the corresponding posterior distribution with a priori smoothness $\alpha$ are plotted in Figure~\ref{SingleDrawPostAlpha} in the appendix.  

For  $\alpha \leq \beta =1$, the variability of the posterior is large so that the true function lies inside the credible band. For larger values of $\alpha$, posterior variability around the posterior mean is much smaller, but the bias of the posterior mean increases, so the true function does not lie inside  the credible band. The value of  $\alpha$  (among the considered values) that gives the posterior with the smallest uncertainty while containing the true function appears to be $1$, which is equal to $\beta$ (the smoothness of $\mu_0(x)$), as predicted by theory (Corollary~\ref{Cor:OptPriorParam} with constant $\tau$).


\subsection{Empirical Bayes posterior  distribution of $\mu$}\label{sec:SimulationEB}

In this section we fix $\alpha > \beta$ and apply the Empirical Bayes estimator of $\tau_\epsilon$ defined by (\ref{def:EBtau}), to check if the corresponding empirical Bayes posterior provides reasonable coverage of the true function, and how it is affected by the noise level $\epsilon$ and different values of  $\gamma$. 

\begin{figure}[h]
\centering
{\includegraphics[width = 1\textwidth, height = 1.75in]{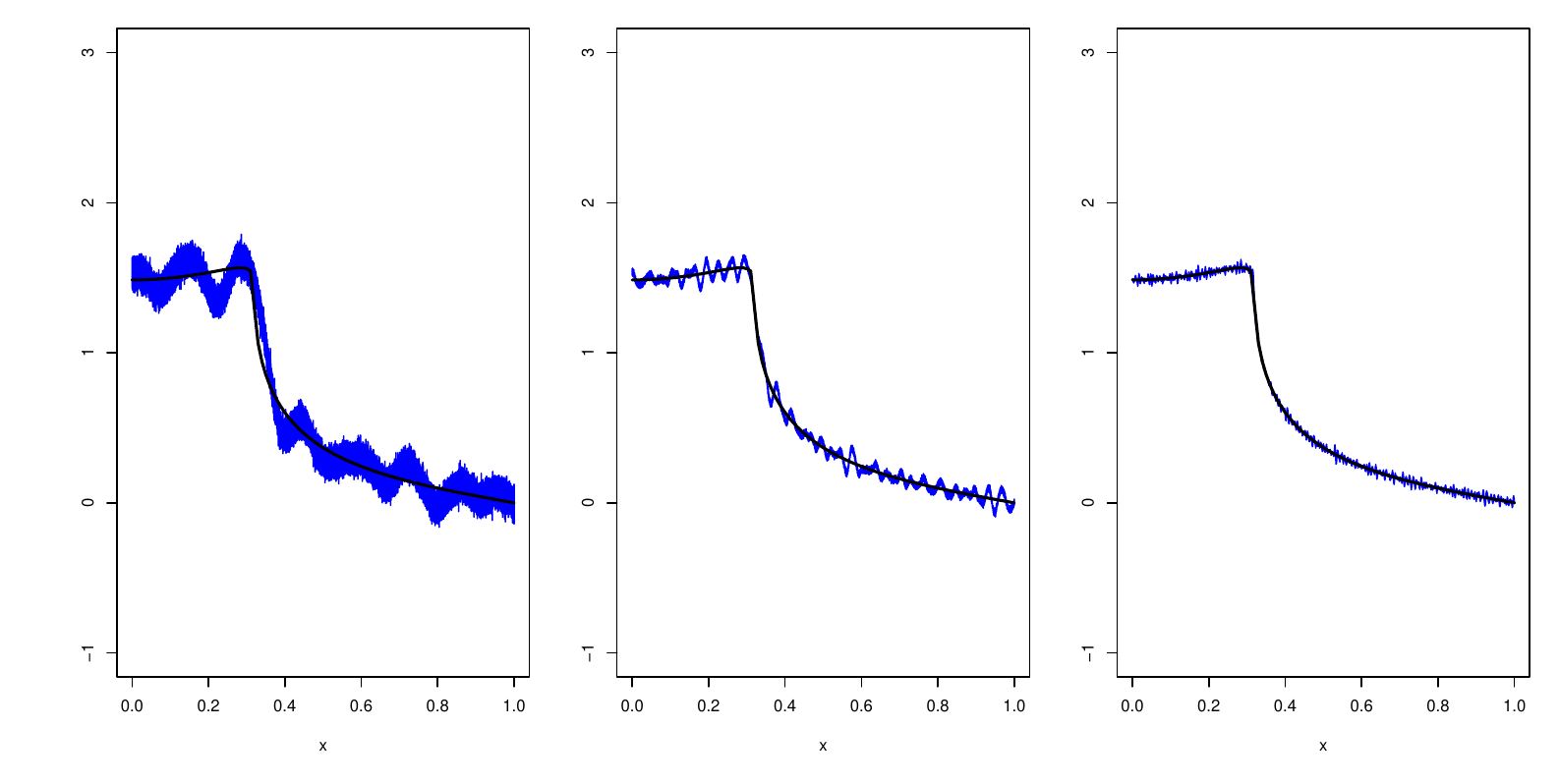}}
\caption{{100 draws from EB posterior with  $\mu_0^N$ (black line), $\epsilon=10^{-4}$ (left), $\epsilon=10^{-6}$ (middle) and $\epsilon=10^{-8}$ (right), $\alpha=1$, $\gamma=0.5$.  }}
\label{fig:EBposterior}
\end{figure}

\begin{figure}[h]
\centering
{\includegraphics[width = 1\textwidth, height = 1.75in]{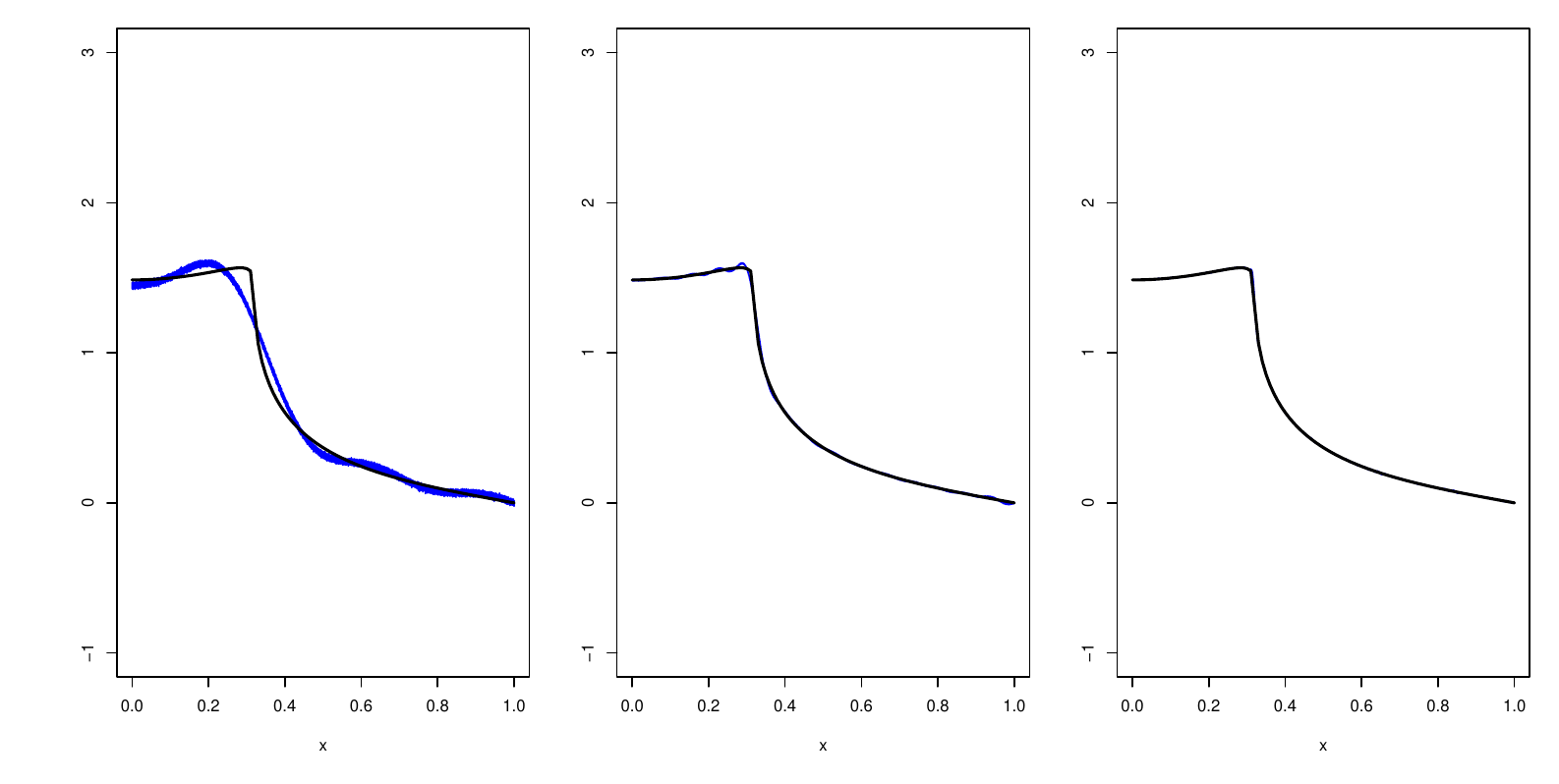}}
\caption{{100 draws from EB posterior with  $\mu_0^N$ (black line), $\epsilon=10^{-4}$ (left), $\epsilon=10^{-6}$ (middle) and $\epsilon=10^{-8}$ (right), $\alpha=5$, $\gamma=0.5$.  }}
\label{fig:EBposteriorAlpha5Gamma05}
\end{figure}

In Figures~\ref{fig:EBposterior} and \ref{fig:EBposteriorAlpha5Gamma05} we can see that with decreasing noise level $\epsilon$ the posterior distribution contracts to the mean, with much smaller bias than for a fixed $\tau$ - the behaviour predicted by the theory. Interestingly, for larger $\alpha$ ($\alpha = 5$ in Figure~\ref{fig:EBposteriorAlpha5Gamma05}), the bias decreases to 0 slower than for smaller $\alpha$ ($\alpha =1$ in Figure~\ref{fig:EBposterior}) however the variability of the posterior decreases faster.

\selectlanguage{english}

Boxplots of values of $\hat\tau $ over 100 simulations for different values of $\alpha$ and different values of $\gamma$ are given in Figure~\ref{fig:EBposteriorHatTauAlpha} in the appendix. 
 In each case, the sampling distribution of $\hat\tau$ concentrates, and the values appear to increase exponentially as functions of $\alpha$. This can be explained by our theory. For self-similar functions and $\alpha >\beta-1/2$, $\hat\tau$ is close to the oracle value $\tau^\star$ (Theorem~\ref{th:EBtau1}) which increases exponentially in $\alpha$ (Corollary~\ref{Cor:OptPriorParam}). 
The values of $\hat\tau$ do not vary much for the considered different values of  $\gamma$ but they do vary with $\alpha$, as expected (Corollary~\ref{Cor:OptPriorParam}).\selectlanguage{english}

We also plotted values of prior variances $\lambda_i$ with plugged in $\hat\tau$ for different values of $\alpha$ and $\gamma$ (Figure~\ref{fig:EBposteriorLambda} in the appendix. 
 Note that the eigenvalues do not much differ for the considered values of $\gamma$, as the only effect is through $\hat\tau$. For each $\gamma$, the values of $\hat\tau$ are such that the values of $\lambda_i = \hat\tau i^{-2\alpha-1}$ are the same at some index $i$ (around $i=50$).  This is expected 
 due to $\lambda_{\hat{i}_{\epsilon}}$  corresponding to the optimal cutoff $\hat{i}_{\epsilon}$ being independent of $\alpha$.


Therefore, the empirical Bayes posterior adapts well to the unknown function and contracts to the true value of the function as the noise level $\epsilon$ vanishes as stated by theory (Theorem~\ref{th:EBtau1}).
 Choosing $\alpha$ larger than $\beta$ and using the Empirical Bayes estimate of $\tau$ does lead to the contraction of the posterior distribution of $\mu$ and good coverage of $\mu_0$. This holds for various values of $\gamma$, including the case of an ill-posed inverse problem ($p+\gamma  >0$), a partially regularised model ($p+\gamma=0$) and a case where the contraction rate is $\epsilon$ ($p+\gamma <-1/2$). The value of $\gamma$ does not have a strong effect on the posterior concentration and on the behaviour of the $\hat\tau$ and the prior eigenvalues $\lambda_i$.




\section{Discussion}

We have considered the inverse problem with Gaussian errors in Hilbert spaces where the covariance operator is not constant but is decomposable in a biorthogonal basis which is an image of a Riesz basis, and both bases span the corresponding Hilbert spaces. We showed that this leads to a sequence space formulation of the inverse problem, and studied its posterior contraction rate, in particular its optimality in the minimax sense, possibly under the error in the forward operator and in the covariance operator. We focused on mildly ill-posed inverse problems with fractional noise, where we also studied optimality of adaptive empirical Bayes posterior distribution.
 We also identified a setting where the posterior distribution can contract at a faster rate, effectively leading to self-regularisation of the inverse problem, which  was also discussed in \citet{JohannesAdaptive2020}.

We extend the general theorem to study the effect of using a plug-in estimator of the variances in sequence space on the posterior contraction rate.  We discuss in which cases the rate of contraction of the posterior distribution is not affected whether the covariance operator is known exactly or observed with error. We consider two types of consistency of the estimator: in terms of the absolute bound on the difference between estimated and true values, and in terms of a relative bound. We study its effect on mildly ill-posed inverse problems, and consider in detail the case of repeated observations. We find that the relative consistency conditions leads to weaker effects of the plug-in on the posterior contraction rate. We have also  applied these results to study an error in operator, e.g. when the eigenfunctions are known but the eigenvalues are estimated, for a fractional noise.

We consider in detail a particular case of the covariance operator whose coefficients in sequence space decrease to 0 at a polynomial rate and illustrate it on the case where the noise is fractional Brownian motion  using fractional wavelets as the bases.  We also derive minimax rates of convergence of estimators of the unknown signal under this model, and show the choice of the prior parameters that leads to this contraction rate of the posterior. We studied the empirical Bayes approach that leads to the corresponding posterior contracting at the optimal rate, under some assumptions on prior smoothness. An alternative approach to adapt in inverse problems uniformly over $S^\beta(A)$ for a range of $\beta$ is a sieve prior proposed by \citet{JohannesAdaptive2020}. Our simulation results confirm the theoretical conclusions.

One can argue that it is  not realistic to assume the knowledge of the covariance operator, and it needs to be estimated in practice. We discuss when estimating of (discretised) eigenfunctions is possible ($\gamma < -1/2$) for repeated observations, and when the number of estimated eigenfunctions is sufficient to achieve the optimal posterior contraction rate of $\mu$.

Another interesting question is estimating the Hurst exponent for fractional noise. While it is possible to estimate $H$, similarly to estimating $\gamma$, using asymptotic expression for coefficients for large indices, fractional wavelets that are used to decompose fractional noise depend on $H$. Therefore, studying the posterior contraction rate with estimated $H$ and using biorthogonal bases that depend on $H$ is a challenging question.



An open question for future research is a joint Bayesian model of the signal and the variance function when the latter is unknown, and its asymptotic behaviour, and the behaviour of the full Bayesian model when the scale parameter is estimated.

\subsection*{Acknowledgments}
Jenovah Rodrigues was supported by The Maxwell Institute Graduate School in Analysis and its Applications, a Centre for Doctoral Training funded by the UK Engineering and Physical Sciences Research Council (grant EP/L016508/01), the Scottish Funding Council, Heriot-Watt University and the University of Edinburgh.

\selectlanguage{english}

\appendix


\section{Minimax rate in sequence space}

\begin{lemma}\label{lem:BelitserLevit} [Theorem 3  \cite{belitser_minimax_1995}]
Consider the problem
$\tilde{Y}_i = \theta_i + \epsilon\tilde{\sigma}_i\xi_i$,   where $\xi_i \sim N(0,1)$ iid for $i=1,2,\ldots$, $\tilde{\sigma}_i \geq 0$ and  $\epsilon > 0$ is  small, and $\theta \in \Theta=  \{(f_i)_{i=1}^{\infty}: \, \sum_i a_i^{2} f_i^2 \leq A^2 \}$.

	Define $c_\epsilon$ to be the solution of the equation
 $\epsilon^2\sum \limits_{i=1}^\infty \tilde{\sigma}_i^2a_i(1-ca_i)_+=c A^2$
 	and $N:=N_\epsilon(\Theta) = \max \{i:a_i \leq c_\epsilon^{-1}\}$. If condition
	\begin{equation}\label{condition}
	\log \epsilon^{-1} \frac{\sum \limits_{i=1}^\infty a_i^2\tilde{\sigma}_i^4(1-c_\epsilon a_i)^2_+}{(\sum \limits_{i=1}^\infty a_i\tilde{\sigma}_i^2(1-c_\epsilon a_i)_+)^2} = o(1),~~ \epsilon \rightarrow 0,
	\end{equation}
	holds, then the minimax rate of convergence of an estimator of $\theta$ in $L^2$ norm over $\Theta$, $r_\epsilon(\Theta)$, satisfies
		$$r_\epsilon^2(\Theta) = \epsilon^2 \sum \limits_{i=1}^N \tilde{\sigma}_i^2(1 -  c_\epsilon a_i(1+o(1))) \,\,\mbox{ as } \epsilon \rightarrow 0.
$$
\end{lemma}

\section{Asymptotic equivalence}

Following \citet{brown_asymptotic_1996}, we show that the considered model (\ref{f2model}) with $f=K \mu$ is equivalent, in Le Cam sense, to a discrete nonparametric regression model
\begin{equation}\label{eq:seqspaceNPR}
Y_i \sim  f_{x_i}  + \xi_i, \quad i=1,\ldots n, \quad x_i = i/n, \quad \xi^{(n)} = (\xi_1,\ldots,\xi_n)^T \sim N(0, \Sigma)
\end{equation}
for $ f_{x_i}$ and $\Sigma_{i,r}$  approximating $f(x_i)$ and $V(x_i, x_r)$, respectively where discretisation operator $D_{x_i}$ is $D_{x_i} = |x_i-x_{i-1}|^{-1} I_{ [x_{i-1},x_i]}$, i.e.
$$
f_{x_i} = \langle f, D_{x_i} \rangle, \quad \Sigma_{ir} =   n^{-1}\langle D_{x_i}, V D_{x_r} \rangle.
$$
For instance, for the white noise model with $V=I$, regular design $x_i-x_{i-1}=1/n$,   $\Sigma = I$.
 This  scheme is typically used in discretisation in inverse problems (Johnstone and Silverman, 1990). 

We consider $f\in \Theta \subset L^2[0,1]$ such that for the considered discretisation operator $D_{x}$,
\begin{equation}\label{eq:AddCondf}
 \lim_{n\to \infty}\sup_{f \in \Theta} n \int_0^1 (f(x) - \sum_{i=1}^n \langle f, D_{x_i}\rangle I(x\in [x_{i-1}, x_i]) )^2 =0.
\end{equation}
This condition is the same as condition (4.5) in \citet{brown_asymptotic_1996} and is a condition for uniform smoothness of functions in $\Theta$ which
  holds e.g. for functions in a H\"older smoothness class $H^{\alpha}(M)$ with $\alpha >1/2$ or with Sobolev smoothness $\beta>1/2$.
 Following the same argument as in Remark 4.6 in \citet{brown_asymptotic_1996}, where this condition is violated, i.e for functions with $\alpha, \beta  \in (0,1/2]$, there is still equivalence for some risk functions, such as with loss $||\hat f - f||^2$.

\begin{theorem}\label{th:asympequiv}
 Consider model (\ref{f2model}) with $f= K\mu \in \Theta$ such that condition (\ref{eq:AddCondf}) holds. Then, model (\ref{f2model}) is asymptotically equivalent, in Le Cam sense, to the model (\ref{eq:seqspaceNPR}) with regular fixed design $(x_i)$ uniformly over $f \in \Theta$ and $V \in {\mathcal V}$ if
$$
\lim_{n\to \infty} \sup_{V\in {\mathcal V}} n^{-1} \left\lVert \left(\langle  nI_{[x_{i-1}, x_i]}, n V  I_{[x_{r-1}, x_r]}\rangle -  \langle D_{x_i}, V  D_{x_r} \rangle \right)_{i,r=1}^n \right\rVert = 0.
$$
\end{theorem}

Consider $\mu_0 \in S^{\beta}(M)$, function $f_0 = K\mu_0 \in S^{\beta+p}(M)$. Then, according to Theorem~\ref{th:asympequiv}, model (\ref{f2model}) is equivalent to discretised model under assumptions (\ref{AssumeEVIP}) and (\ref{AssumeEVVar}) $\beta+p >1/2$, and eigenvectors of $V$ are uniformly bounded and $\gamma < -1/2$.

However, the challenge is usually to find the limiting operator $V$ given matrix $\Sigma$, so in this sense this result is trivial.

 \hspace{2pc}{\it Proof of Theorem~\ref{th:asympequiv}.}

We generally follow the steps in the proof of Theorem 4.1 in \citet{brown_asymptotic_1996}, with their Theorem 3.1 being in the heart of the proof.

For the continuous model (\ref{f2model}), given $f\in \Theta$, define a piece-wise constant function $\bar f_n(x) = D_{x_i} f(x)$.
 Equation (2.2) and Theorem 3.1 in \citet{brown_asymptotic_1996} together with condition (\ref{eq:AddCondf}) imply that Le Cam's distance between the models (\ref{f2model}) with function $f$ and (\ref{f2model}) with function $\bar f_n$ tends to 0 as $n\to \infty$.

For the continuous model
$dY_t =   \bar f_n(t) dt  + \epsilon  dW_t$,  sufficient statistics for unknown $\bar f_n$ are $S_i = n\int_{x_{i-1}}^{x_i}  Y_t dt$, $i=1, \ldots, n$.  $E S_i = \bar f_n(x_i)$ and $Cov(S_i, S_r) = n \langle I_{[x_{i-1}, x_i]}, V  I_{[x_{r-1}, x_r]} \rangle = \epsilon^2  \langle \delta_{\tilde x_i}, V  \delta_{\tilde x_r} \rangle$  for some $ \tilde x_k\in (x_{k-1}, x_k)$.

Since $S_i$ are sufficient statistics for $\bar f_n$ in model (\ref{f2model}), Le Cam distance between these two models over such $\bar f_n$ is 0 (Lemma 3.2 and Theorem 3.1 in \citet{brown_asymptotic_1996}) which  proves the theorem.

\hfill $\square$ \newline


\section{Additional simulation plots}

We continue to investigate performance of the Bayesian model on simulated data, in the set up in Section~\ref{sec:SimData}. 

\subsection{Non-adaptive posterior}

\begin{figure}
\includegraphics[width=1\textwidth,height=5in]{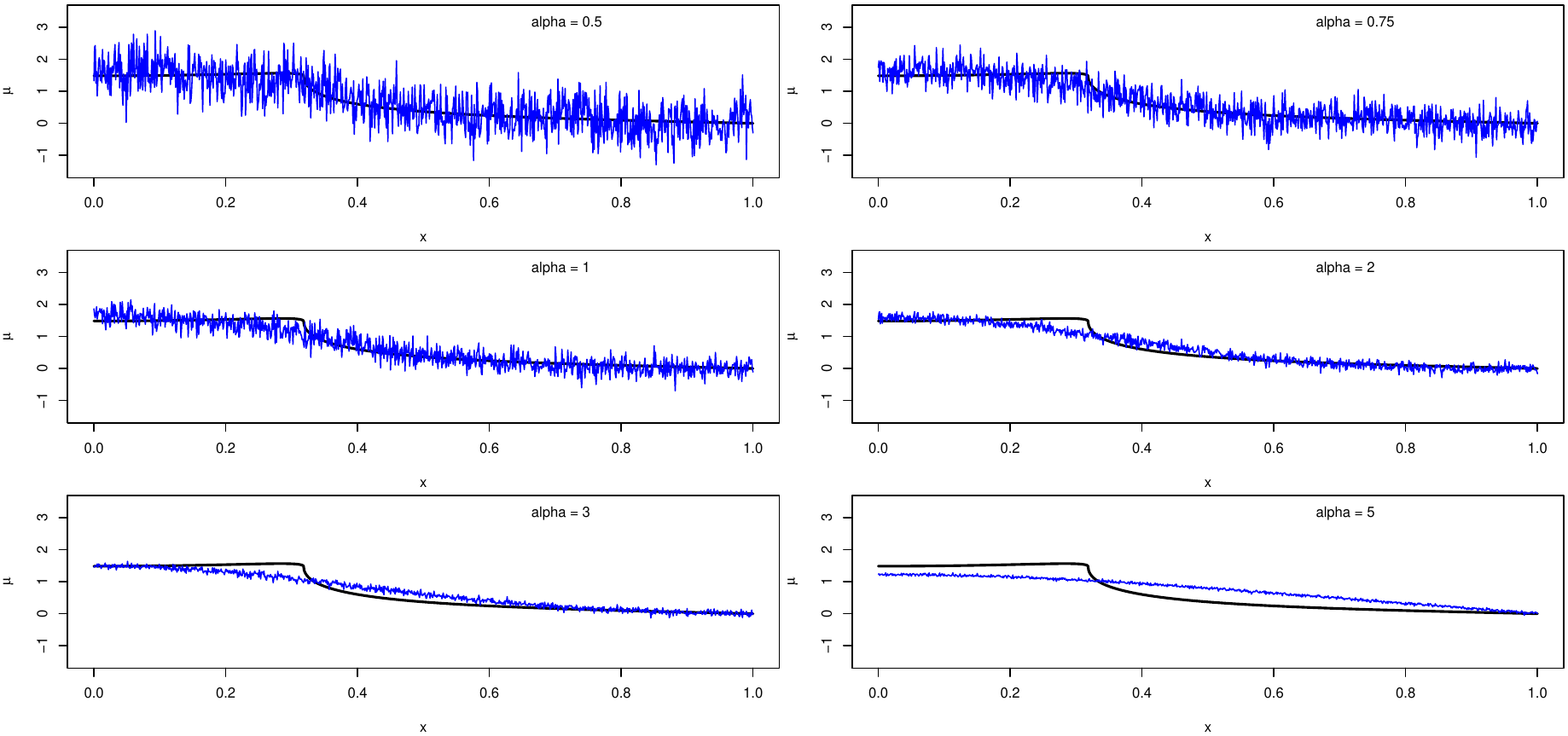}
\caption{{A single draw from the posterior of $\mu^N(x)$ for $\alpha=(0.5,0.75,1,2,3,5)$ respectively, with $\epsilon=10^{-2}$ and $N=2000$ in all cases; $\mu^N_0(x)$ - black line.}} \label{SingleDrawPostAlpha}
\end{figure}

Here we focus on the case of fixed prior smoothness $\alpha$, set $\tau=1$ and plot a single draw from the corresponding posterior distribution for each $\alpha$ (Figure~\ref{SingleDrawPostAlpha}).

\subsection{Empirical Bayes posterior}

Here we study how the   concentration of the empirical Bayes posterior is affected by different values of $\gamma$.

We generate data sets with several values of noise parameter $\gamma$: $\gamma =-2$ (leads to  rate of convergence $\epsilon$), $\gamma=-1$ (corresponds to $p+\gamma=0$ and hence the rate of convergence of a direct problem), $\gamma = 0.5$ (corresponds to the rate of convergence of an ill-posed problem with $\tilde p =p+\gamma= 1.5$).


We start with $\gamma = 0.5$. Level of ill-posedness is $\tilde p = p+\gamma =1.5$. This case is discussed  in the paper. 

\begin{figure}[h]
\centering
{\includegraphics[width = 0.4\textwidth]{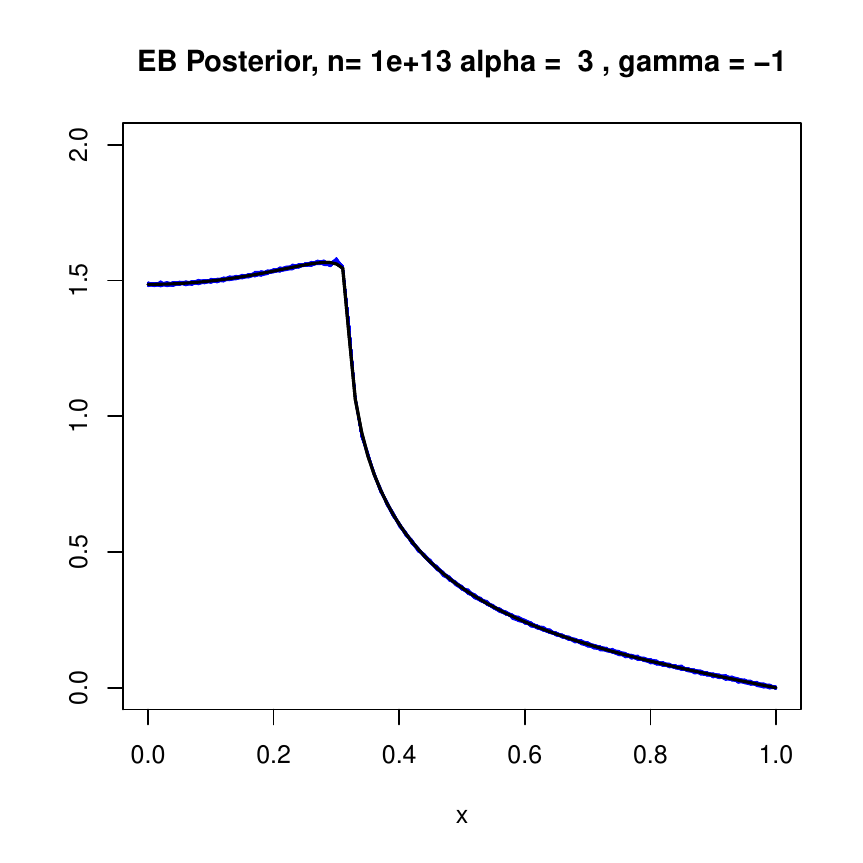}}
{\includegraphics[width = 0.4\textwidth]{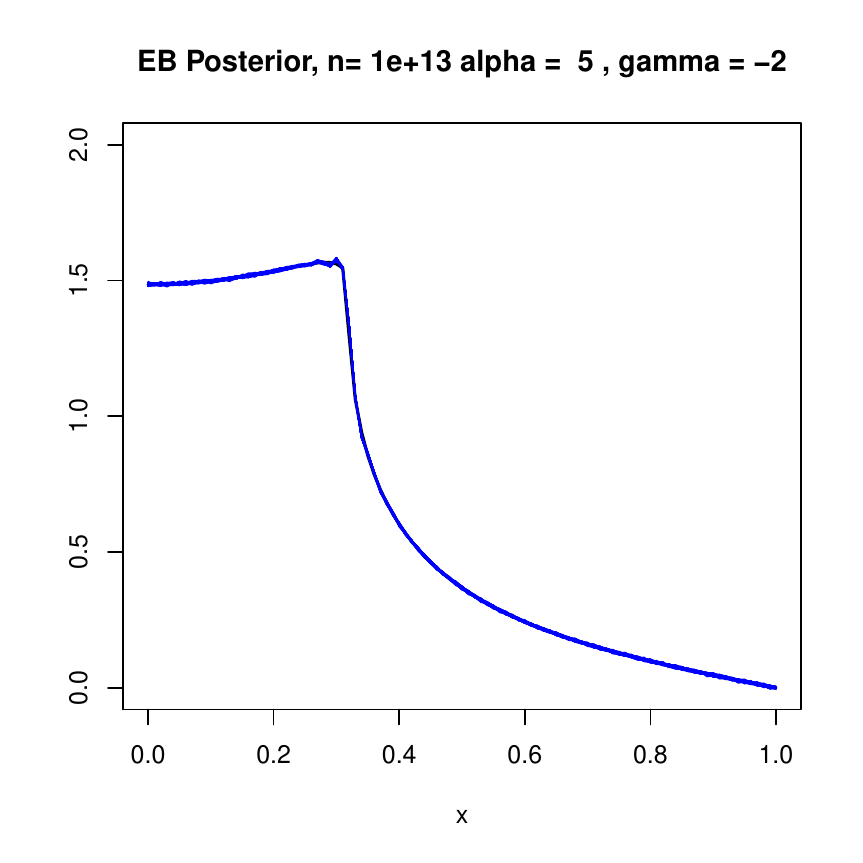}}
\caption{{100 draws from EB posterior (blue), $\alpha = 5$ and $\epsilon=10^{-6.5} \approx 3 \cdot 10^{-7}$;    $\gamma=-1$ (left) and $\gamma=-2$ (right).  Black - truth. }}
\label{fig:EBposteriorAlpha5n10pow13}
\end{figure}

For $\gamma = -1$, the level of ill-posedness is $\tilde p =p+\gamma= 0$, corresponding to the rate of convergence of direct problem. Conditions of Theorem~\ref{th:EBtau1} are satisfied for  $\alpha \geq (B_0-1)/2$.
  For $\gamma = -2$, the rate of contraction is $\epsilon$. Conditions of Theorem~\ref{th:EBtau1} are satisfied for  $\alpha \geq \max((B_0+1)/2,3/2)$.]
  Draws from the EB posterior   for $\gamma = -1$ and $\gamma = -2$ with $\alpha = 5$ are given in Figure~\ref{fig:EBposteriorAlpha5n10pow13}. We can see that the value of $\gamma$ does not have a strong effect on the posterior concentration and on the behaviour of the $\hat\tau$ and the eigenvalues $\lambda_i = \hat\tau i^{-2\alpha-1}$. See discussion in Section~\ref{sec:SimulationEB}. 

\begin{figure}[h]
\centering
{\includegraphics[width = 0.3\textwidth]{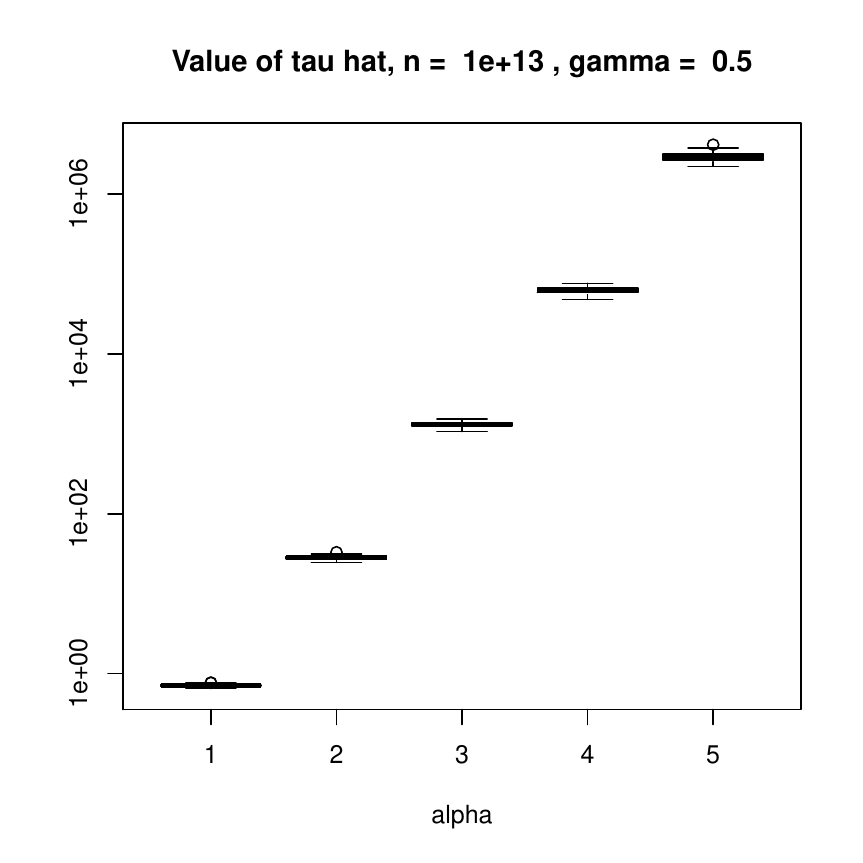}}
{\includegraphics[width = 0.3\textwidth]{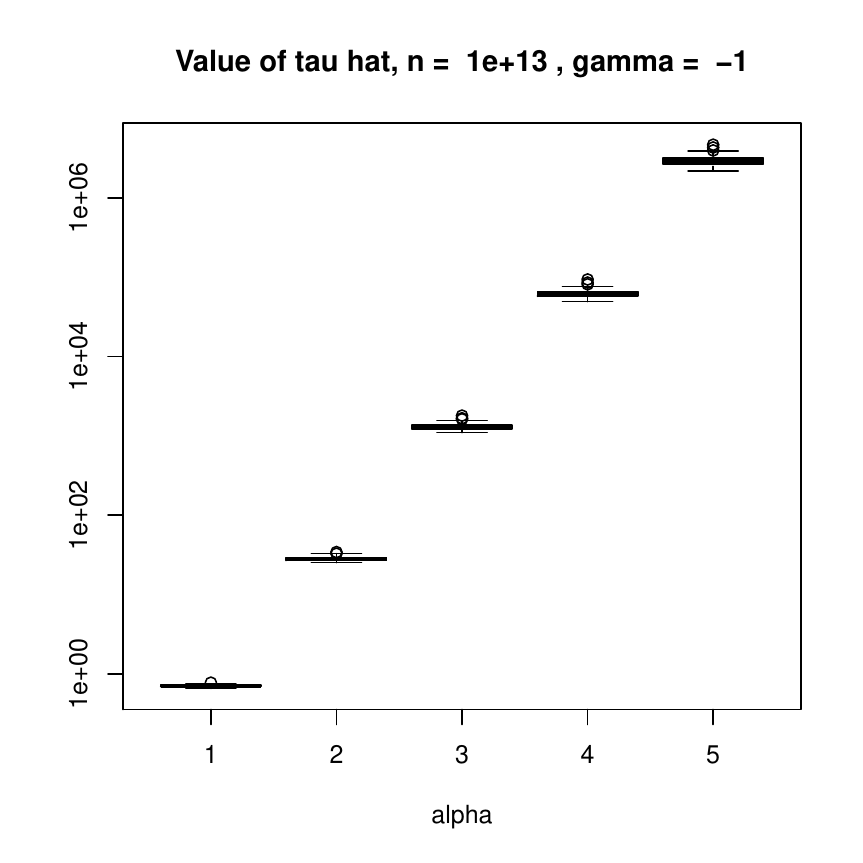}}
{\includegraphics[width = 0.3\textwidth]{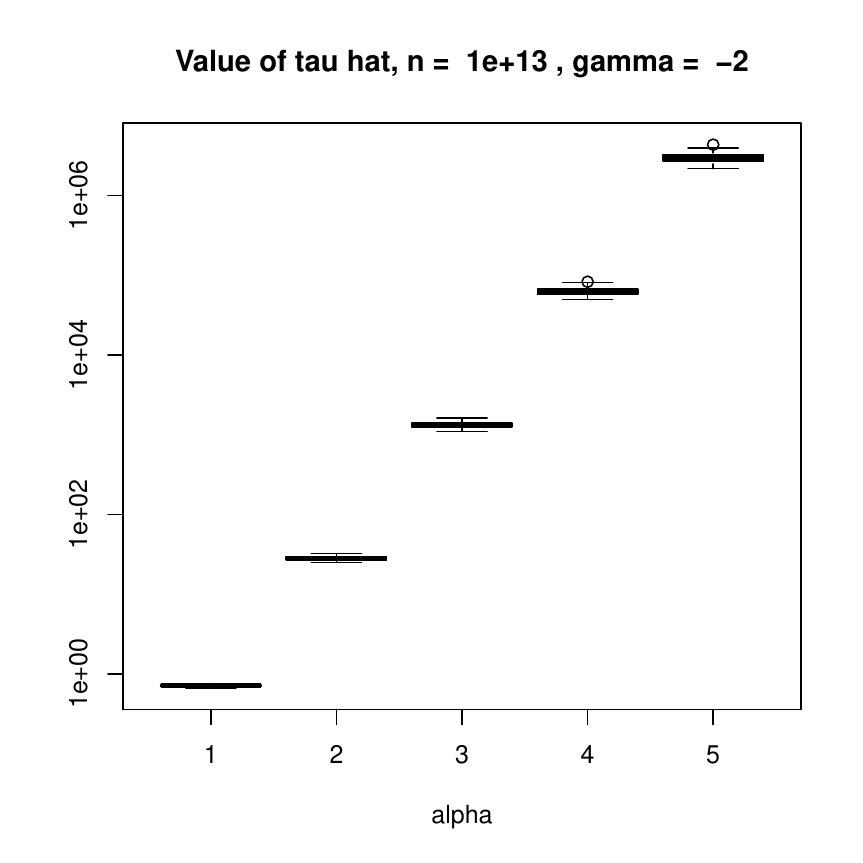}}
\caption{{Boxplots of $\hat\tau$ over 100 draws for $\epsilon=10^{-13/2}$ and various values of $\alpha$; $\gamma=0.5$ (left), $\gamma=-1$ (middle) and $\gamma=-2$ (right).}}
\label{fig:EBposteriorHatTauAlpha}
\end{figure}

\begin{figure}[h]
\centering
{\includegraphics[width = 0.3\textwidth]{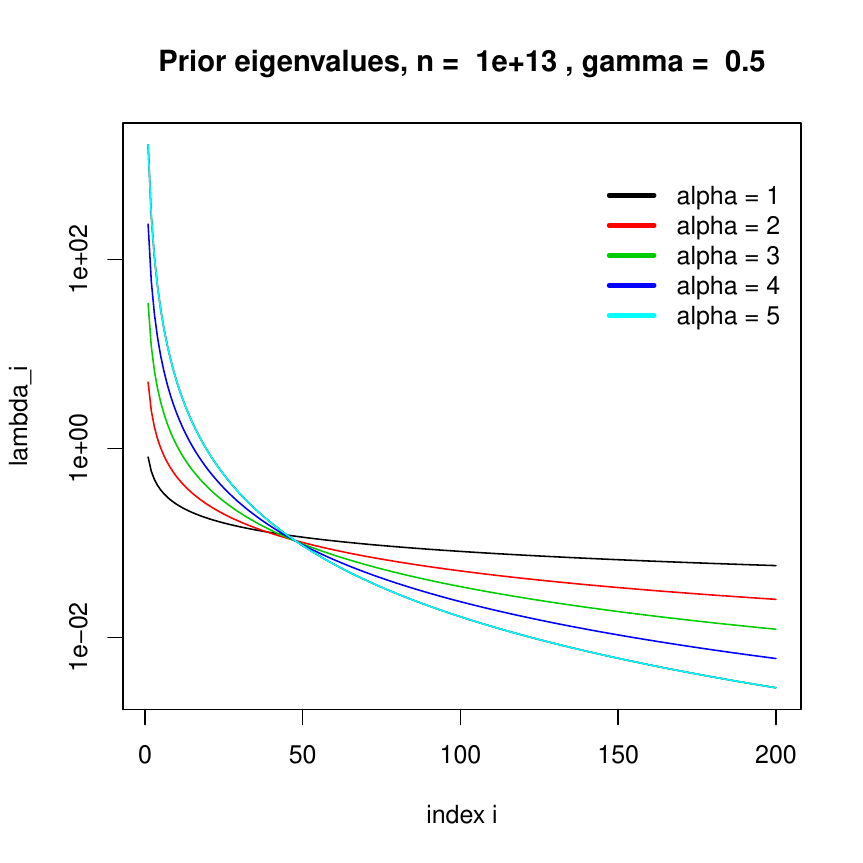}}
{\includegraphics[width = 0.3\textwidth]{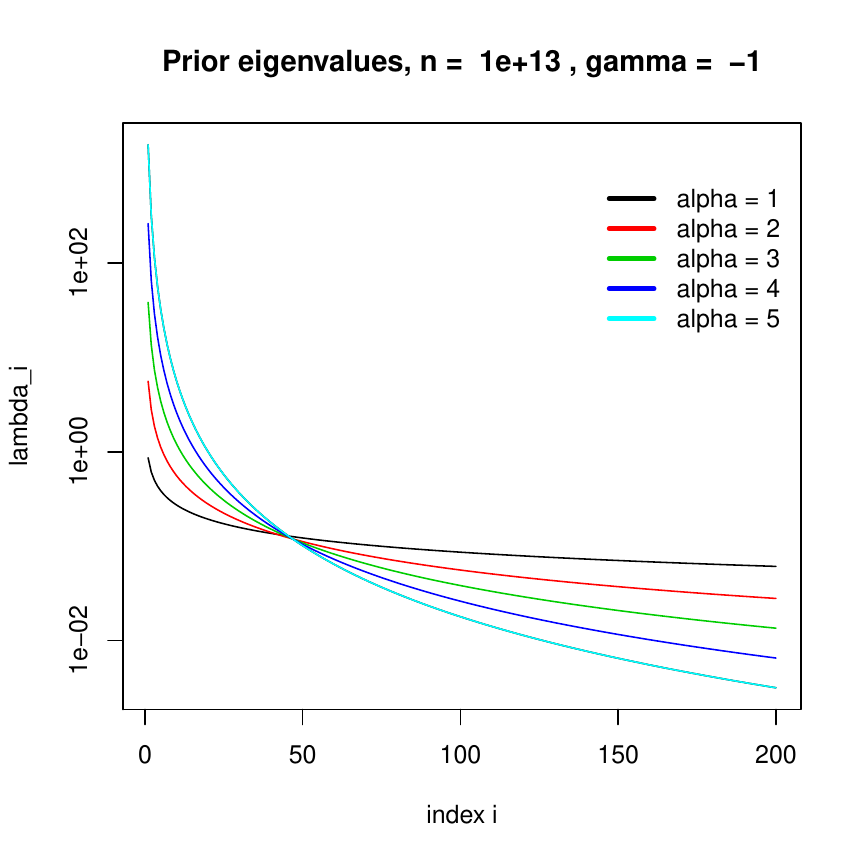}}
{\includegraphics[width = 0.3\textwidth]{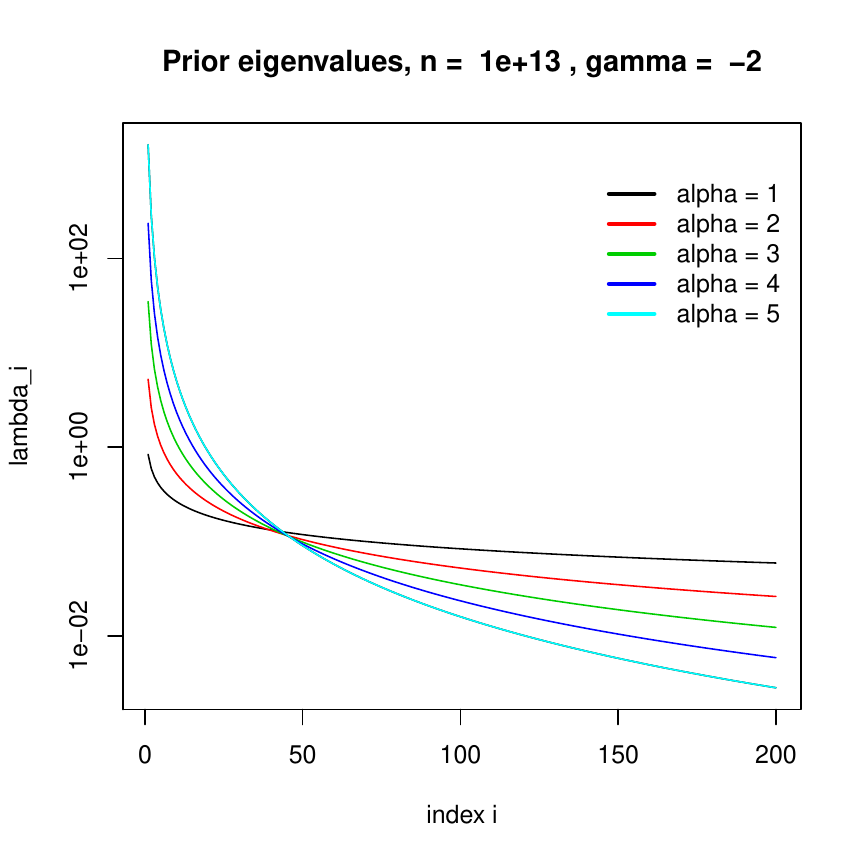}}
\caption{{ Values of $\sqrt{\lambda_i}$ with $\tau_\epsilon = \sqrt{\hat\tau}$,  $\epsilon=10^{-13/2}$; $\gamma=0.5$ (left), $\gamma=-1$ (middle) and $\gamma=-2$ (right) (single draw). }}
\label{fig:EBposteriorLambda}
\end{figure}

%


\section{Proof of the main theorem on contraction of posterior distribution}


\hspace{2pc}{\it Proof of Theorem \ref{th:contractiontheoremGeneral}.}

If we show that $\mathbb{E}_{\mu_0} \mathbb{E} (||\mu -\mu_0||_2^2 \mid Y) \leq C\rate^{2}$ for some $C \in (0,1\infty)$ independent of $\epsilon$  then, by Markov's inequality, for any $M \to \infty$,
$$
\mathbb{P}( ||\mu -\mu_0||_2 \geq M \rate \mid Y) \leq  M^{-2} \rate^{-2} \mathbb{E} (||\mu -\mu_0||_2^2 \mid Y) \leq  C M^{-2} \to 0
$$
hence   $\rate$ is the rate of contraction of the posterior distribution.

We assumed that we can write $\mu =\sum_{\windex \in \Windex} \mu_{\windex} e_{0,\windex}$ where $\{e_{0,\windex}\}$ is a Riesz basis.
Since it is not an orthogonal basis, we cannot use Parseval's identity, however we can use the Riesz basis property to obtain
\begin{eqnarray*}
A \sum_{\windex \in \Windex}  (\mu_i -\mu_{0i})^2  \leq  ||\mu-\mu_0||^2  \leq  B \sum_{\windex \in \Windex}  (\mu_i -\mu_{0i})^2,
\end{eqnarray*}
and therefore we have both lower and upper bound on the $L^2$ norm of $\mu-\mu_0$ in terms of the corresponding $\ell^2$ sequence norm in basis $\{e_{0,\windex}\}$.

Therefore, it is sufficient to study
\begin{eqnarray*}
 \sum_{\windex} \mathbb{E} \left[(\mu_{\windex} -\mu_{0,\windex})^2 \mid Y\right] = \sum_{\windex} \left[\Var(\mu_{\windex}\mid Y) + (\mathbb{E} [\mu_{\windex}  |\mid Y] -\mu_{0,\windex})^2\right].
\end{eqnarray*}
For $\windex \in D$, Lemma~\ref{lem:PostRateDepTerms} implies that under the conditions of the theorem, the expected value of the sum over $D$ with respect to the true distribution of the data is bounded by $C \epsilon^2$.

Mapping $\Windex_I$ to $\mathbb{N}$, we write the sum over $\Windex_I$ indexed equivalently by $\mathbb{N}$. Taking the expected value with respect to the true distribution of the data  and using the explicit form of the posterior distribution
\begin{eqnarray*}
  \mathbb{E}_{\mu_0}[\Var(\mu_i\mid Y)] + \mathbb{E}_{\mu_0}\left[(\mathbb{E} [\mu_i  |\mid Y] -\mu_{0i})^2\right]
&=& \mathbb{E}_{\mu_0}\left[\frac{ Y_i\kappa_i\lambda_i}{ \lambda_i\kappa_i^2 + \epsilon^2 \sigma^2_i} -\mu_{0i}\right]^2 + \frac{\sigma^2_i\lambda_i}{\epsilon^{-2} \lambda_i\kappa_i^2 + \sigma^2_i}\\
&=&  \frac{ \epsilon^2 \sigma_i^2 \kappa_i^2\lambda_i^2}{ [\lambda_i\kappa_i^2 + \epsilon^2 \sigma^2_i]^2}+  \mu_{0i}^2\left[\frac{ \kappa_i^2 \lambda_i}{ \lambda_i\kappa_i^2 + \epsilon^2 \sigma^2_i} -1\right]^2 + \frac{\sigma^2_i\lambda_i}{\epsilon^{-2} \lambda_i\kappa_i^2 + \sigma^2_i}\\
&\asymp&   \frac{\sigma^2_i\lambda_i}{\epsilon^{-2} \lambda_i\kappa_i^2 + \sigma^2_i} + \mu_{0i}^2\left[\frac{ \kappa_i^2 \lambda_i}{ \lambda_i\kappa_i^2 + \epsilon^2 \sigma^2_i} -1\right]^2
\end{eqnarray*}
as the first term is less than the third one.

Recall that $ \sigma^{2}_i/[\lambda_i\kappa_i^2]$ is an increasing sequence.
 Denote 
  $\Windex_{\epsilon} = \{i: \,  \sigma^2_i/[\lambda_i\kappa_i^2] > \epsilon^{-2}\}$.
Then,
\begin{eqnarray*}
  S_1 &=&  \sum_i \frac{\sigma^2_i\lambda_i}{\epsilon^{-2} \lambda_i\kappa_i^2 + \sigma^2_i}  \asymp   \epsilon^{2} \sum_{i\notin \Windex_{\epsilon}} \sigma^2_i \kappa_i^{-2}
 + \sum_{i \in \Windex_{\epsilon}}  \lambda_i,
\end{eqnarray*}
and
\begin{eqnarray*}
 S_2 &=& \sum_i\mu_{0i}^2\left[\frac{ \epsilon^2 \sigma^2_i}{  \kappa_i^2\lambda_i + \epsilon^2 \sigma^2_i}  \right]^2
=||\mu_0||^2_Q  \sum_i  \bar{\mu}_{0,i}^2 \left[\frac{ \epsilon^2 \sigma^2_i a_i^{-1}}{  \kappa_i^2\lambda_i + \epsilon^2 \sigma^2_i}  \right]^2\\
&\asymp& ||\mu_0||^2_{Q}  \epsilon^4  \sum_{i\notin \Windex_{\epsilon}}  \bar{\mu}_{0,i}^2 \left[\frac{\sigma^2_i a_i^{-1}}{  \kappa_i^2\lambda_i }  \right]^2
 +||\mu_0||^2_{Q} \sum_{i\in \Windex_{\epsilon}}  \bar{\mu}_{0,i}^2 a_i^{-2}\\
&\leq&   ||\mu_0||^2_{Q}  \left[ \epsilon^4   \max_{i\notin \Windex_{\epsilon}} \left[\frac{\sigma^2_i a_i^{-1}}{  \kappa_i^2\lambda_i }  \right]^2 + C  \sup_{i \in \Windex_{\epsilon}} a_{i}^{-2} \right]
\end{eqnarray*}
where $\bar{\mu}_{0,i}^2:= \mu_{0,i}^2 i^{2\beta}/||\mu_0||^2_{Q}$ and $||\mu_0||^2_{Q} = \sum_i \mu_{0i}^2 a_i^{2}$.  The lower bound can be proved by taking $\mu_0$ such that $\bar\mu_{0i}=1$  for one of the $i \notin \Windex_{\epsilon}$
and $\bar\mu_{0i}=0$  otherwise to get the first term; to get the second term, denote $i_\epsilon = \arg \sup_{i \in \Windex_{\epsilon}} a_i^{-2}$ (which can be infinite) and set $\bar\mu_{0i}=1$  for  $i=i_\epsilon$  and $\bar\mu_{0i}=0$  for $i\neq i_\epsilon$ and $i \in \Windex_{\epsilon}$.

Hence,
\begin{eqnarray*}
 S_2 &\asymp&||\mu_0||^2_{S^\beta} \left[ \epsilon^4  \max_{i\notin \Windex_{\epsilon}} \left[\frac{\sigma^2_i a_i^{-1}}{  \kappa_i^2\lambda_i }  \right]^2 +  a_{i_\epsilon}^{-2}  \right].
\end{eqnarray*}

Combining these results together, we obtain
\begin{eqnarray*}
\mathbb{E}_{\mu_0} \mathbb{E} (||\mu-\mu_0||^2 \mid Y) &\asymp &
\epsilon^{2} \sum_{i\notin \Windex_{\epsilon}} \sigma^2_i \kappa_i^{-2}
 + \sum_{i \in \Windex_{\epsilon}}  \lambda_i+
   \epsilon^4   \max_{i\notin \Windex_{\epsilon}} \left[\frac{\sigma^2_i a_i^{-1}}{  \kappa_i^2\lambda_i }  \right]^2 +  \sup_{i \in \Windex_{\epsilon}} a_{i}^{-2}.
\end{eqnarray*}

Hence, setting $\rate$ such that
$$
\rate^{-2} \mathbb{E}_{\mu_0} \mathbb{E} (||\mu-\mu_0||^2 \mid Y) \asymp \rate^{-2} \left[\epsilon^2+\epsilon^2\sum_{i\notin \Windex_{\epsilon}}  \sigma^2_i \kappa_i^{-2}
 + \sum_{i \in \Windex_{\epsilon}}  \lambda_i +   \epsilon^4   \max_{i\notin \Windex_{\epsilon}} \left[\frac{\sigma^2_i a_i^{-1}}{  \kappa_i^2\lambda_i }  \right]^2 +  \sup_{i \in \Windex_{\epsilon}} a_{i}^{-2}
\right]=O(1)
$$
as $\epsilon\to 0$  ensures that for every $M \rightarrow \infty$,
$$
\mathbb{E}_{\mu_0}\mathbb{P}(\{\mu: \,||\mu-\mu_0|| \geq M \rate |Y\}) \leq M^{-2} \rate^{-2} \mathbb{E}_{\mu_0} \mathbb{E} (||\mu-\mu_0||^2 \mid Y)\rightarrow 0 \quad \text{ as } \epsilon\to 0.
$$
\hfill $\square$ \newline

Consider the inverse problem (\ref{f2model}) formulated in the sequence space under Assumption~\ref{assumptionB1} with noise $W$ satisfying Assumption~\ref{eq:WnoiseDecomp}, with prior distribution (\ref{eq:SeqPrior}). Assume that the true function $\mu_0$ satisfies Assumption~\ref{smoothGeneral}. Assume also that $\min_{\windex \in D}\lambda_{\windex} \geq 1$ and $\left\lVert K_D^{-1} V_D K_D^{-1}\right\rVert  < \infty$.

\begin{lemma}\label{lem:PostRateDepTerms} For the inverse problem (\ref{f2model}) formulated in the sequence space under Assumption~\ref{assumptionB1} with noise $W$ satisfying Assumption~\ref{eq:WnoiseDecomp}, with prior distribution (\ref{eq:SeqPrior}),
\begin{eqnarray*}
\mathbb{E}_{\mu_0}\sum_{\windex \in D} \mathbb{E} [(\mu_\windex - \mu_{0,\windex})^2 \mid y_D] \leq\epsilon^{2} || \mu_{0, D}||^2\, \left\lVert (K_D^T  V_D^{-1} K_D +\epsilon^{2}\Lambda_D^{-1})^{-1} \Lambda_D^{-1}\right\rVert +2\epsilon^2 \trace\left( (K_D^T  V_D^{-1} K_D +\epsilon^{2}\Lambda_D^{-1})^{-1} \right).
\end{eqnarray*}
Therefore, for $\mu_0 \in Q((a_{\windex}),A)$, this term is bounded by $C \epsilon^2$ with $C =  2C_{D,2} +A^2 C_{D,1}$, as long as there exist constants $C_{D,1}, \, C_{D,2} >0$ independent of $\epsilon$ such that
\begin{eqnarray}\label{eq:condDep}
\left\lVert(K_D^T  V_D^{-1} K_D +\epsilon^{2}\Lambda_D^{-1})^{-1} \Lambda_D^{-1}\right\rVert \leq C_{D,1} \quad \trace\left( (K_D^T  V_D^{-1} K_D +\epsilon^{2}\Lambda_D^{-1})^{-1} \right)\leq C_{D,2}.
\end{eqnarray}
\end{lemma}

\begin{proof}[Proof of Lemma \ref{lem:PostRateDepTerms}]
 The bias of the posterior mean is (omitting indices $_D$ in $K$, $V$ and $\Lambda$ for simplicity)
 \begin{eqnarray*}
\mathbb{E} (\mu_\windex|y_D) - \mu_{0,\windex} &=& (K^T  V^{-1} K +\epsilon^{2}\Lambda^{-1})^{-1} K^T V^{-1} [K \mu_0 +\epsilon V^{1/2} \xi_D] - \mu_0  \\
&=& [(K^T  V^{-1} K +\epsilon^{2}\Lambda^{-1})^{-1} K^T V^{-1} K - I] \mu_0 +\epsilon (K^T  V^{-1} K +\epsilon^{2}\Lambda^{-1})^{-1} K^T V^{-1} V^{1/2} \xi_D\\
&=& -\epsilon^{2}(K^T  V^{-1} K +\epsilon^{2}\Lambda^{-1})^{-1} \Lambda^{-1} \mu_0 +\epsilon (K^T  V^{-1} K +\epsilon^{2}\Lambda^{-1})^{-1} K^T V^{-1} V^{1/2} \xi_D,\end{eqnarray*}
hence, for $\xi_D \sim N(0, I_{|D|})$,
 \begin{eqnarray*}
\mathbb{E}_{\mu_0} ||\mathbb{E} (\mu_\windex|y_D) - \mu_{0,\windex}||^2 &=& (K^T  V^{-1} K +\epsilon^{2}\Lambda^{-1})^{-1} K^T V^{-1} [K \mu_0 +\epsilon V^{1/2} \xi_D] - \mu_0  \\
&=& \epsilon^{4}\left\lVert(K^T  V^{-1} K +\epsilon^{2}\Lambda^{-1})^{-1} \Lambda^{-1} \mu_0\right\rVert^2 \\
&& +\epsilon^2 \trace\left( (K^T  V^{-1} K +\epsilon^{2}\Lambda^{-1})^{-1} K^T V^{-1}  K^T (K^T  V^{-1} K +\epsilon^{2}\Lambda^{-1})^{-1}\right)\\
&\leq&
\epsilon^{2} || \mu_{0, D}||^2\, ||(K^T  V^{-1} K +\epsilon^{2}\Lambda^{-1})^{-1} \Lambda^{-1}|| +\epsilon^2 \trace\left( (K^T  V^{-1} K +\epsilon^{2}\Lambda^{-1})^{-1} \right).
 \end{eqnarray*}
\end{proof}


\begin{remark}\label{rem:CondDepSuff}
Conditions (\ref{eq:condDep}) hold if  $\min_{\windex \in D}\lambda_{\windex} \geq 1$ and $||K_D^{-1} V_D K_D^{-1}|| \leq  C_{KVK} < \infty$, with $C_{D,1} =   C_{KVK}$ and  $C_{D,2} = |D| \, C_{KVK}$.
\end{remark}

\begin{proof}[Proof of Remark \ref{rem:CondDepSuff}]
Conditions (\ref{eq:condDep}) hold if  $||\Lambda_D^{-1}|| \leq 1$ and $||K_D^{-1} V_D K_D^{-1}|| \leq  C_{KVK} < \infty$ which imply
$$
||(K_D^T  V_D^{-1} K_D +\epsilon^{2}\Lambda_D^{-1})^{-1} \Lambda_D^{-1}|| \leq ||(K_D^T  V_D^{-1} K_D )^{-1} ||\leq  C_{KVK}
$$
and
$$
\trace\left( (K_D^T  V_D^{-1} K_D +\epsilon^{2}\Lambda_D^{-1})^{-1}  \right) =
\trace\left( K_D^{-1}  V_D K_D^{-1} (I +\epsilon^{2}\Lambda_D^{-1} K_D^{-1}  V_D K_D^{-1})^{-1}  \right)
\leq \trace\left( K_D^{-1}  V_D K_D^{-1} \right) \leq  |D| \, C_{KVK}.
$$
\end{proof}

\begin{lemma} \label{lem:PostRateDepTermsFBM}
 Under assumptions of Proposition~\ref{cor:vaguelettes}, conditions (\ref{eq:condDep}) of Lemma \ref{lem:PostRateDepTerms} hold for the scaling coefficients of fractional noise.
\end{lemma}

\begin{proof}[Proof of Lemma \ref{lem:PostRateDepTermsFBM}]
Recall that for fractional noise under the compactly supported wavelets, set $D$ is finite, and covariance matrix $V_D$ has variances $\sigma_{-d}^2$ on the diagonal and $V_{D; k,m} = \sigma_{-d}^2 \rho^{(-d)}_{|k-m|}$ where
$$
\sigma_{-d}^2  = \sum_{\ell=0}^{\infty} \left[\gamma_\ell^{(-d)}\right]^2 < \infty\,\, \text{ and } \,\, \rho^{(-d)}_{m} = \sigma_{-d}^{-2} \sum_{\ell=0}^{\infty} \gamma_\ell^{(-d)} \gamma_{\ell+|m|}^{(-d)} <\infty.
$$
Condition  (\ref{eq:MildIPwavelets}) implies that all diagonal elements of $K_D$ are non-zero and finite, and under the assumptions of the proposition, all diagonal elements of $\Lambda_D$ are $\geq 1$ and finite. Therefore,
$$
\trace\left( (K_D^T  V_D^{-1} K_D +\epsilon^{2}\Lambda_D^{-1})^{-1} \right) \leq  \trace\left( (K_D^{-1}  V_D K_D^{-1} \right) \lesssim \trace(V_D) \lesssim |D| \sigma^2_d <\infty,
$$
and
$$
||(K_D^T  V_D^{-1} K_D +\epsilon^{2}\Lambda_D^{-1})^{-1} \Lambda_D^{-1}||  \leq ||K_D^{-1}  V_D K_D^{-1} || \leq \trace\left( (K_D^{-1}  V_D K_D^{-1} \right) <\infty.
$$
Therefore, conditions (\ref{eq:condDep}) hold.
\end{proof}


\section{Proofs for plugin estimators}


\subsection{General case with absolute bound}\label{sec:PlugInAbsBound}

We use an isomorphism $\Windex \to \mathbb{N}$ and  index sequences by $i \in \mathbb{N}$.

\begin{assumption} \label{AssumeEVVarPlugIn}
Assume that the estimated  eigenvalues ($\{\hat\sigma_i^2\}$) of operator $V$  are independent of $Y$, and there exists a constant $c_0$ such that
$$
P(|\hat \sigma_i^2 -  \sigma_i^2|\leq c_0 \epsilon_\sigma, \, i=1,2,\ldots) \to 1 \text{ as }  \epsilon_\sigma \to 0.
$$
\end{assumption}
This type of assumption is useful if for estimating $\mu$, it is important to consider the full sequence $i=1,2,\ldots$ without truncation. We included this case for completeness, as we find that in all cases we consider the relative error yields better results.

If $\sigma_i\to 0$, an estimator of $\sigma_i^2$ can be truncated at some small value $\tilde \epsilon_{\sigma}^2$:
 $\tilde \sigma_i^2 = \max(\tilde \epsilon_{\sigma}, \hat\sigma_i^2), \quad i=1,2,\ldots$.
 If the error rate $c_0 \epsilon_{\sigma}$ is known, we can use $\tilde \epsilon_{\sigma}=c_0 \epsilon_{\sigma}$.

Plugging-in the estimator of $\sigma_i^2$, the posterior distribution of $\mu_i$ (\ref{eq:SeqPost}) becomes
$$
p(\mu_i \mid Y_i, \tilde \sigma_i^2)  \sim N\left(\frac{ Y_i\kappa_i\lambda_i}{ \lambda_i\kappa_i^2 + \epsilon^2 \tilde \sigma^2_i}, \frac{\tilde \sigma^2_i\lambda_i}{\epsilon^{-2} \lambda_i\kappa_i^2 + \tilde\sigma^2_i}\right)\, \quad i=1,2,\ldots, \text{ independently}.
$$

In particular, Theorem~\ref{th:contractiontheoremGeneral}    implies that when $\sigma_i\geq c_1 >0$ for all $i$, then using a consistent plug-in estimator (i.e. when plugged in values are $\hat\sigma_i^2 = \sigma_i^2 + o(1)$), the effect on the contraction rates of the posterior distribution is a larger constant on the definition of the summation set. When sequence $(\sigma_i)$ decreases to 0 as $i$ increases then the error of estimation of $V$ can affect the rate, with the effect depending on the speed of decay of $(\sigma_i)$.

We investigate how the contraction rates are affected as $\epsilon_\sigma \rightarrow 0$ when $\tilde \epsilon_{\sigma}=c_0 \epsilon_{\sigma}$.

\begin{theorem}\label{th:GeneralPlugin}
Consider the inverse problem (\ref{f2model}) formulated in the sequence space under Assumption~\ref{assumptionB1} with noise $W$ satisfying Assumption~\ref{eq:WnoiseDecomp}, with prior distribution (\ref{eq:SeqPrior}). Assume that the true function $\mu_0$ satisfies Assumption~\ref{smoothGeneral}. Assume also that $\min_{\windex \in D}\lambda_{\windex} \geq 1$ and $||K_D^{-1} V_D K_D^{-1}||  < \infty$.

 Consider an estimator of  $\{\sigma_i^2\}$ satisfying assumption (\ref{AssumeEVVarPlugIn}),  and use a plugin estimator $\tilde \sigma_i^2 = \max(c_0 \epsilon_{\sigma}, \hat\sigma_i^2)$.

Then, for every $M \rightarrow \infty$,
$$
 \mathbb{P}(\{\mu: \, ||\mu-\mu_0|| \geq M {\rate}_{plugin}|Y, \tilde \sigma_i\})   \stackrel{\mathbb{P}_{\mu_0, V}}{\to} \text{ as } \epsilon \rightarrow 0
$$
uniformly over $\mu_0$ in  $Q((a_i), A)$ where ${\rate}_{plugin}$ is given by
\begin{eqnarray}\label{eq:PlugInRate}
{\rate}_{plugin}^2  =
    \epsilon^{2}   \sum_{i \in \bar I_{\epsilon}(2) \cup \bar I_\sigma(2)} \sigma_i^2  \kappa_i^{-2}
 +  \epsilon^2 c_0 \epsilon_\sigma\sum_{i\in   I_\sigma(2)  \cap \bar I_{\sigma\epsilon}(1/3)}   \kappa_i^{-2} \\
+     \sum_{i\in I_{\epsilon}(2/3)}  \lambda_i  +  \sum_{i\in   I_\sigma(2)\cap  I_{\sigma\epsilon}(1/3)}   \lambda_i \notag\\
+    \epsilon^4
 \max\left\{ \max_{i\in \bar I_{\epsilon}(1)} \left[\frac{ \sigma_i^2  }{ a_i \kappa_i^2\lambda_i }  \right]^2, \max_{i\in \bar I_{\sigma \epsilon}(1)} \left[\frac{ \epsilon_\sigma i^{-\beta}}{  \kappa_i^2\lambda_i }  \right]^2\right\}
 +  \max_{i\in I_{\epsilon}(1)\cap I_{\sigma \epsilon}(1)} a_i^{-2} \notag
\end{eqnarray}
where
\begin{eqnarray}
I_{\epsilon}(a) = \{i: \,  \sigma_i^2/[\lambda_i\kappa_i^2]>  a \epsilon^{-2}\},\,\,
I_\sigma(a) = \{i: \,  \sigma^2_i < a c_0 \epsilon_\sigma\},\,\,
I_{\sigma \epsilon}(a) = \{i: \, c_0 \epsilon_\sigma /[\lambda_i \kappa_i^2] > a \epsilon^{-2} \}.
\end{eqnarray}
In particular, if $I_{\sigma}(2/3) \subseteq I_{\epsilon}(2)$  then the rate consists of the same terms as the contraction rate for the known $V$, with slightly larger constant in $I_{\epsilon}$ in some terms:
\begin{eqnarray}\label{eq:PlugInRateSimple}
{\rate}_{plugin}^2 =
  \epsilon^{2}   \sum_{i \in \bar I_{\epsilon}(2)} \sigma_i^2  \kappa_i^{-2}
+     \sum_{i\in I_{\epsilon}(2)}  \lambda_i
+    \epsilon^4   \max_{i\in \bar I_{\epsilon}(1)} \left[\frac{ \sigma_i^2 i^{-\beta}}{  \kappa_i^2\lambda_i }  \right]^2
 +  \max_{i\in I_{\epsilon}(1)} a_i^{-2}.
\end{eqnarray}

\end{theorem}


Note that set $I_\sigma(a)$ corresponds to the values $\sigma_i^2$ that are not estimated well. Sets $\bar I_{\epsilon}(a)$ and $\bar I_{\sigma \epsilon}(a)$ represent the observations that contributes strongly to estimation of $\mu$, where $\sigma_i^2$ are estimated well and not well, respectively. In the simpler version (\ref{eq:PlugInRateSimple}), the first term is the leading part of the posterior variance, the second term corresponds to the bias caused by the prior, and the last term is the bias that comes from the true function. The third term represents saturation and is discussed in detail in \citet{AgapiouMathe2017}. Note that this rate is similar to the rate of the case where $\sigma_i^2$ are known except the value of $a$ in the sets $I_{\epsilon}(a)$ here is not always $a=1$ as it was in Theorem~\ref{th:contractiontheoremGeneral}. In the rate with all terms \eqref{eq:PlugInRate},  each of these terms consists of two parts: where the variance is estimated well and where it is not. 

\hspace{2pc}{\it Proof  of Theorem \ref{th:GeneralPlugin}.}
The proof is following the lines of the proof of Theorem~\ref{th:contractiontheoremGeneral}. 

Define $\Omega_{\sigma}= \{  |\hat \sigma_i^2 -  \sigma_i^2|\leq c_0 \epsilon_\sigma, \, i=1,2,\ldots  \}$. Under assumption (\ref{AssumeEVVarPlugIn}),
$P(\Omega_{\sigma} ) \to 1$ as    $\epsilon_\sigma \to 0$. Then, on $\Omega_\sigma$, $\hat\sigma^2_i \leq \sigma_i^2 + c_0 \epsilon_\sigma$ and
$$
\tilde \sigma_i^2 = \max(\hat\sigma^2_i, a_\sigma) \geq \max(\sigma_i^2 - c_0 \epsilon_\sigma,  a_\sigma).
$$
Consider the case $a_\sigma = c_0 \epsilon_\sigma$.
 Then,
$$
\mathbb{E}_{\mu_0, V} \mathbb{P}( ||\mu -\mu_0||_2 \geq M \rate \mid Y, \hat V) \leq \mathbb{E}_{\mu_0, V} \mathbb{P}( ||\mu -\mu_0||_2 \geq M \rate \mid Y, \hat V) I(\Omega_{\sigma}) + P_{\mu_0, V}(\Omega_{\sigma})
$$
and the latter term goes to 0 by assumption (\ref{AssumeEVVarPlugIn}). For the former term,
 Markov's inequality implies:
$$
\mathbb{P}( ||\mu -\mu_0||_2 \geq M \rate \mid Y, \hat V) \leq  M^{-2} \rate^{-2} \mathbb{E} (||\mu -\mu_0||_2^2 \mid Y, \hat V).
$$

Using the Riesz basis assumption, as in the proof of Theorem~\ref{th:contractiontheoremGeneral}, we have that
\begin{eqnarray*}
\mathbb{E} (||\mu-\mu_0||^2\mid Y, \hat V) &\asymp& \sum_i \mathbb{E} [(\mu_i -\mu_{0i})^2 \mid Y, \hat V] = \sum_i \left[\Var(\mu_i\mid Y, \hat V) + \left(\mathbb{E} [\mu_i  \mid Y, \hat V] -\mu_{0i}\right)^2\right].
\end{eqnarray*}

Taking expected value with respect to the true distribution of $Y$ (under the assumption that it is independent of $(\hat\sigma_i)$) and using the explicit form of the posterior distribution, we have on $\Omega_\sigma$
\begin{eqnarray*}
&&  \mathbb{E}_{\mu_0}[\Var(\mu_i\mid Y, \tilde V ) ] + \mathbb{E}_{\mu_0}\left[(\mathbb{E} [\mu_i  |\mid Y, \tilde V] -\mu_{0i})^2  \right]
=  \frac{ \epsilon^2 \sigma_i^2 \kappa_i^2\lambda_i^2}{ [\lambda_i\kappa_i^2 + \epsilon^2 \tilde\sigma^2_i]^2}+  \mu_{0i}^2\left[\frac{ \epsilon^2 \tilde\sigma^2_i}{ \lambda_i\kappa_i^2 + \epsilon^2 \tilde\sigma^2_i}\right]^2 + \frac{\tilde\sigma^2_i\lambda_i}{\epsilon^{-2} \lambda_i\kappa_i^2 + \tilde\sigma^2_i}\\
&& \leq  \frac{ \epsilon^2 \lambda_i^2 \kappa_i^2 \sigma_i^2}{ [\lambda_i\kappa_i^2 +  \epsilon^{2} \max(\sigma_i^2 - c_0 \epsilon_\sigma,c_0 \epsilon_\sigma)]^2 }  + \frac{[\sigma_i^2 +c_0 \epsilon_\sigma]\lambda_i}{\epsilon^{-2} \lambda_i\kappa_i^2 + [\sigma_i^2 + c_0 \epsilon_\sigma]} +  \mu_{0i}^2\left[\frac{ [\sigma_i^2 + c_0 \epsilon_\sigma]}{ \lambda_i\kappa_i^2/\epsilon^2 + [\sigma_i^2 + c_0 \epsilon_\sigma]} \right]^2.
\end{eqnarray*}

Denote $I_{\epsilon}(a) = \{i: \,  \sigma_i^2/[\lambda_i\kappa_i^2]>  a \epsilon^{-2}\}$ and $I_\sigma (a) = \{i: \,  \sigma^2_i < a c_0 \epsilon_\sigma\}$.
 Define also $i_{\epsilon}(a) = \min\{i: \, \sigma_i^2 /[\lambda_i\kappa_i^2]> a\epsilon^{-2}\}$.
 Note that if $ \sigma^{2}_i/[\lambda_i\kappa_i^2]$ increases then $I_{\epsilon}(a) =\{i \geq i_{\epsilon}(a)\}$.
 And if $\sigma_i$ decreases then we can write  $  I_\sigma(a)= \{i \geq i_{\sigma}(a)\}$ where $i_{\sigma}(a) = \max\{i: \, \sigma^2_i \leq a c_0 \epsilon_\sigma\}$. 
  Note that  $I_{\epsilon}(a_1) \subseteq I_{\epsilon}(a_2)$ and  $I_{\sigma}(a_1) \subseteq I_{\sigma}(a_2)$ if $a_1 > a_2$.
Then,
\begin{eqnarray*}
 S_1 &=&
  \sum_i   \frac{ \epsilon^2 \lambda_i^2 \kappa_i^2 \sigma_i^2}{ [\lambda_i\kappa_i^2 +  \epsilon^{2} \max(\sigma_i^2 - c_0 \epsilon_\sigma,c_0 \epsilon_\sigma)]^2 }\\
  &\leq& \epsilon^2  \sum_{ i\in \bar I_{\epsilon}(1)}     \sigma_i^2/\kappa_i^2
+  \epsilon_\sigma^{-2}  \epsilon^{-2}  \sum_{  i\in  I_{\epsilon}(1) \cap I_{\sigma}(2)}   \lambda_i^2 \kappa_i^2 \sigma_i^2
  +  \sum_{  i\in  I_{\epsilon}(1) \cap \bar I_{\sigma}(2)}  \frac{ \epsilon^2 \lambda_i^2 \kappa_i^2 \sigma_i^2}{ [\lambda_i\kappa_i^2 +  \epsilon^{2} \sigma_i^2/2]^2 }.
  \end{eqnarray*}
The first   and the last terms are  the same as in the case of no plug-in, up to a constant in the indices.
In the polynomial case, the second term is
\begin{eqnarray*}
  \epsilon_\sigma^{-2}  \tau_\epsilon^4 \epsilon^{-2}  \sum_{ i> \max(  i_{\epsilon}(1) , i_{\sigma})}  i^{-4\alpha-2 -2 p +2\gamma}
  \leq   \epsilon_\sigma^{-2}  \tau_\epsilon^4 \epsilon^{-2} [\max(  i_{\epsilon}(1) , i_{\sigma})]^{-4\alpha-1 -2 p +2\gamma} \\
  =  \epsilon_\sigma^{-2}  \tau_\epsilon^4 \epsilon^{-2} [ \max[(\tau_\epsilon^2 \epsilon^{-2})^{1/(2\tilde p + 2\alpha+1)}, \epsilon_\sigma^{1/(2\gamma)} ]]^{-4\alpha-1 -2 p +2\gamma}
  \end{eqnarray*}
If $ (\tau_\epsilon^2 \epsilon^{-2})^{1/(2\tilde p + 2\alpha)}> \epsilon_\sigma^{1/(2\gamma)} $ then the upper bound is
$$
  \epsilon_\sigma^{-2}  \tau_\epsilon^4 \epsilon^{-2} (\tau_\epsilon^2 \epsilon^{-2})^{-(4\alpha+1 +2 p -2\gamma)/(2\tilde p + 2\alpha+1)}
= \epsilon_\sigma^{-2}  \tau_\epsilon^2 (\tau_\epsilon^2 \epsilon^{-2})^{-(2\alpha  -2\gamma)/(2\tilde p + 2\alpha+1)}
= \epsilon_\sigma^{-2}  \epsilon^{2} (\tau_\epsilon^2 \epsilon^{-2})^{(1+2p)/(2\tilde p + 2\alpha+1)}
$$
which tends to 0 if $ \epsilon^{2[(2\alpha +6\gamma)] }  <  \tau_\epsilon^{-2(1+2p-4\gamma) }$.


If $ (\tau_\epsilon^2 \epsilon^{-2})^{1/(2\tilde p + 2\alpha)}< \epsilon_\sigma^{1/(2\gamma)} $ then the upper bound is
\begin{eqnarray*}
\epsilon_\sigma^{-2}  \tau_\epsilon^4 \epsilon^{-2} \epsilon_\sigma^{(-4\alpha-1 -2 p +2\gamma)/(2\gamma)} =
 \tau_\epsilon^4 \epsilon^{-2} \epsilon_\sigma^{(-4\alpha-1 -2 p -2\gamma)/(2\gamma)}
  \end{eqnarray*}
 which tends to 0 if $  \tau_\epsilon^4 \epsilon^{-2} \epsilon_\sigma^{(-4\alpha-1 -2 p -2\gamma)/(2\gamma)} =o(1)$.
 It holds if
 $ (\tau_\epsilon^2 \epsilon^{-2})^{1/(2\tilde p + 2\alpha)} < (\tau_\epsilon^4 \epsilon^{-2})^{1/(4\alpha+1 +2 \tilde p)}  =o(1) \epsilon_\sigma^{1/(2\gamma)}$.

The next term is
 \begin{eqnarray*}
S_2  &\leq &
   \sum_{i\in   I_\sigma(2)}  \frac{ 3\epsilon^2 \lambda_i  c_0 \epsilon_\sigma}{ \lambda_i\kappa_i^2 + 3\epsilon^2  c_0 \epsilon_\sigma  }
 +  \sum_{i\in \bar I_\sigma(2)}  \frac{ 1.5\epsilon^2 \lambda_i  \sigma_i^2  }{ \lambda_i\kappa_i^2 +  1.5\epsilon^2  \sigma_i^2  } \\
 &\leq &
 \sum_{i\in   I_\sigma(2)  \, \& \, i\in  I_{\sigma\epsilon}(1/3)}   \lambda_i
+  3\epsilon^2 c_0 \epsilon_\sigma\sum_{i\in   I_\sigma(2)  \, \& \, i\in  \bar I_{\sigma\epsilon}(1/3)}   \kappa_i^{-2}
   +  1.5 \epsilon^2\sum_{i\in \bar I_\sigma(2)  \, \& \, i\in \bar I_{\epsilon}(2/3)}     \sigma_i^2    \kappa_i^{-2}\\
  && +  \sum_{i\in \bar I_\sigma(2)  \, \& \, i\in I_{\epsilon}(2/3)}  \lambda_i  .
\end{eqnarray*}
where $I_{\sigma \epsilon}(a) = \{i: \, c_0 \epsilon_\sigma /[\lambda_i \kappa_i^2] > a \epsilon^{-2} \}$.

Combining these terms, we obtain
 \begin{eqnarray*}
S_1+S_2  &\leq & 2 \epsilon^{2}   \sum_{i \in \bar I_{\epsilon}(2) \cup I_\sigma(2)} \sigma_i^2 \kappa_i^{-2}
 +  3\epsilon^2 c_0 \epsilon_\sigma\sum_{i\in   I_\sigma(2)  \, \& \, i\in  \bar I_{\sigma\epsilon}(1/3)}   \kappa_i^{-2} \\
  &+&   \sum_{i\in I_{\epsilon}(2/3)}  \lambda_i  +  \sum_{i\in   I_\sigma(2)  \, \& \, i\in  I_{\sigma\epsilon}(1/3)}   \lambda_i
\end{eqnarray*}

The remaining term is
\begin{eqnarray*}
 S_3 &=&
||\mu_0||^2_{S^\beta} \sum_i  \bar{\mu}_{0,i}^2 \left[\frac{ \epsilon^2 [\sigma_i^2 + c_0 \epsilon_\sigma] i^{-\beta}}{  \kappa_i^2\lambda_i + \epsilon^2 [\sigma_i^2 + c_0 \epsilon_\sigma]}  \right]^2\\
&\leq& ||\mu_0||^2_{S^\beta}  \epsilon^4  \sum_{i\in \bar I_{\epsilon}(1) \cup \bar I_{\sigma \epsilon}(1)}  \bar{\mu}_{0,i}^2 \left[\frac{[\sigma_i^2 + c_0 \epsilon_\sigma]  i^{-\beta}}{  \kappa_i^2\lambda_i }  \right]^2
 +||\mu_0||^2_{S^\beta} \sum_{i\in I_{\epsilon}(1)\cap I_{\sigma \epsilon}(1)}  \bar{\mu}_{0,i}^2 i^{-2\beta}\\
&\leq&   ||\mu_0||^2_{S^\beta}  \left[ \epsilon^4   \max_{i\in \bar I_{\epsilon}(1) \cup \bar I_{\sigma \epsilon}(1)} \left[\frac{[\sigma_i^2 + c_0 \epsilon_\sigma] i^{-\beta}}{  \kappa_i^2\lambda_i }  \right]^2 + C \max_{i\in I_{\epsilon}(1)\cap I_{\sigma \epsilon}(1)} i^{-2\beta} \right]
\end{eqnarray*}
where $\bar{\mu}_{0,i}^2:= \mu_{0,i}^2 i^{2\beta}/||\mu_0||^2_{S^\beta}$ and $||\mu_0||^2_{S^\beta} = \sum_i \mu_{0i}^2 i^{2\beta}$.


We can rewrite the first term as follows:
\begin{eqnarray*}
\max_{i\in \bar I_{\epsilon}(1) \cup \bar I_{\sigma \epsilon}(1)} \left[\frac{[\sigma_i^2 + c_0 \epsilon_\sigma] i^{-\beta}}{  \kappa_i^2\lambda_i }  \right]^2
=\max\left\{ \max_{i\in \bar I_{\epsilon}(1)} \left[\frac{2\sigma_i^2 i^{-\beta}}{  \kappa_i^2\lambda_i }  \right]^2, \max_{i\in \bar I_{\sigma \epsilon}(1)} \left[\frac{2c_0 \epsilon_\sigma i^{-\beta}}{  \kappa_i^2\lambda_i }  \right]^2\right\}.
\end{eqnarray*}

Combining these results together, we obtain that on $\Omega_{\sigma}$,
\begin{eqnarray*}
C \mathbb{E}_{\mu_0} \mathbb{E} (||\mu-\mu_0||^2 \mid Y, \hat V) \lesssim
  \epsilon^{2}   \sum_{i \in \bar I_{\epsilon}(2) \cup \bar I_\sigma(2)} \sigma_i^2  \kappa_i^{-2}
 +  \epsilon^2 c_0 \epsilon_\sigma\sum_{i\in   I_\sigma(2)  \cap \bar I_{\sigma\epsilon}(1/3)}   \kappa_i^{-2} \\
+     \sum_{i\in I_{\epsilon}(2/3)}  \lambda_i  +  \sum_{i\in   I_\sigma(2)\cap  I_{\sigma\epsilon}(1/3)}   \lambda_i\\
+    \epsilon^4
 \max\left\{ \max_{i\in \bar I_{\epsilon}(1)} \left[\frac{ \sigma_i^2 i^{-\beta}}{  \kappa_i^2\lambda_i }  \right]^2, \max_{i\in \bar I_{\sigma \epsilon}(1)} \left[\frac{ \epsilon_\sigma i^{-\beta}}{  \kappa_i^2\lambda_i }  \right]^2\right\}
 +  \max_{i\in I_{\epsilon}(1)\cap I_{\sigma \epsilon}(1)} i^{-2\beta}
\end{eqnarray*}
and hence this can be taken as $\rate^2$, up to a constant. This bound is uniform in $\mu_0\in S^\beta(A)$, and proves the rate stated in the theorem.

Note that
\begin{eqnarray*}
I_\epsilon(a_1) \cap I_\sigma(a_2) &= \{i: \sigma_i^2/[\lambda_i \kappa_i^2]\geq a_1 \epsilon^{-2}  \} \cap \{i: \sigma_i^2 < a_2 c_0 \epsilon_\sigma \}  \\
&\subseteq  \{i: c_0 \epsilon_\sigma/[\lambda_i \kappa_i^2] \geq  a_1/a_2 \epsilon^{-2}  \}  = I_{\sigma \epsilon}(a_1/a_2),
\end{eqnarray*}
which also implies that $\bar I_{\sigma \epsilon}(a_1/a_2) \subseteq \bar I_\epsilon(a_1) \cup \bar I_\sigma(a_2)$, and hence
$$
I_\sigma(2)  \cap \bar I_{\sigma\epsilon}(1/3) \subseteq  \bar I_\epsilon(2/3)  \cap I_\sigma(2).
 \text{ and }
I_\epsilon(2)  \cap \bar I_{\sigma\epsilon}(1/3) \subseteq   I_\epsilon(2) \cap \bar I_\sigma(2/3).
$$

If $I_{\sigma}(2/3) \subseteq I_{\epsilon}(2)$  then $\bar I_{\epsilon}(2) \cap \bar I_\sigma(2/3) = \bar I_{\epsilon}(2)$ and $\bar I_{\epsilon}(2) \cup I_\sigma(2/3) = \emptyset$, and the contraction rate can be written as
\begin{eqnarray*}
C \mathbb{E}_{\mu_0} \mathbb{E} (||\mu-\mu_0||^2 \mid Y, \hat V) \leq
  \epsilon^{2}   \sum_{i \in \bar I_{\epsilon}(2)} \sigma_i^2  \kappa_i^{-2}
+     \sum_{i\in I_{\epsilon}(2)}  \lambda_i
+    \epsilon^4   \max_{i\in \bar I_{\epsilon}(1)} \left[\frac{ \sigma_i^2 i^{-\beta}}{  \kappa_i^2\lambda_i }  \right]^2
 +  \max_{i\in I_{\epsilon}(1)} i^{-2\beta}
\end{eqnarray*}
which is almost the upper bound in the case $V$ is known, except the summation in the first two sums is over $I_{\epsilon}(2)$ rather than $I_{\epsilon}(1)$. Here  the  sum of $\epsilon^2 c_0 \epsilon_\sigma \kappa_i^{-2}$ over $I_\epsilon(2)  \cap \bar I_{\sigma\epsilon}(1/3)$ is bounded by the sum of $\epsilon^2 \lambda_i $ (up to a constant). The theorem is proved.

\hfill $\square$ \newline


\hspace{2pc}{\it Proof  of Theorem \ref{th:contractiontheoremGeneralPlugIn2}.}

The proof is following the lines of the proof of Theorem~\ref{th:contractiontheoremGeneral}. 
 Define $\Omega_{\sigma}= \{  |\hat \sigma_i^2/\sigma_i^2-1|\leq c_0 \epsilon_\sigma, \, i=1,2,\ldots, N  \}$. Under assumption (\ref{AssumeEVVarPlugIn}),
$P(\Omega_{\sigma} ) \to 1$ as    $\epsilon_\sigma \to 0$. Then, on $\Omega_\sigma$,
$$
\mathbb{E}_{\mu_0, V} \mathbb{P}( ||\mu -\mu_0||_2 \geq M \rate \mid Y, \hat V) \leq \mathbb{E}_{\mu_0, V} \mathbb{P}( ||\mu -\mu_0||_2 \geq M \rate \mid Y, \hat V) I(\Omega_{\sigma}) + P_{\mu_0, V}(\Omega_{\sigma})
$$
and the latter term goes to 0 by Assumption \ref{AssumeEVVarPlugInRelative}. For the former term,
 Markov's inequality implies:
$$
\mathbb{P}( ||\mu -\mu_0||_2 \geq M \rate \mid Y, \hat V) \leq  M^{-2} \rate^{-2} \mathbb{E} (||\mu -\mu_0||_2^2 \mid Y, \hat V).
$$
Using the Riesz basis property, we have that
\begin{eqnarray*}
\mathbb{E} (||\mu-\mu_0||^2\mid Y, \hat V) &\asymp& \sum_i \mathbb{E} [(\mu_i -\mu_{0i})^2 \mid Y, \hat V] = \sum_i \left[\Var(\mu_i\mid Y, \hat V) + (\mathbb{E} [\mu_i  \mid Y, \hat V] -\mu_{0i})^2\right].
\end{eqnarray*}
Taking expected value with respect to the true distribution of $Y$ (under the assumption that it is independent of $(\hat\sigma_i)$) and using the explicit form of the posterior distribution, we have on $\Omega_\sigma$,
\begin{eqnarray*}
&&  \mathbb{E}_{\mu_0}[\Var(\mu_i\mid Y, \hat V ) ] + \mathbb{E}_{\mu_0}[(\mathbb{E} [\mu_i  |\mid Y, \hat V] -\mu_{0i})^2  ]
=  \frac{ \epsilon^2 \sigma_i^2 \kappa_i^2\lambda_i^2}{ [\lambda_i\kappa_i^2 + \epsilon^2 \hat\sigma^2_i]^2}+  \mu_{0i}^2\left[\frac{ \epsilon^2 \hat\sigma^2_i}{ \lambda_i\kappa_i^2 + \epsilon^2 \hat\sigma^2_i}\right]^2 + \frac{\hat\sigma^2_i\lambda_i}{\epsilon^{-2} \lambda_i\kappa_i^2 + \hat\sigma^2_i}\\
& \leq& \frac{ \epsilon^2 \sigma_i^2 \kappa_i^2\lambda_i^2}{ [\lambda_i\kappa_i^2 + \epsilon^2 \sigma^2_i (1- c_0\epsilon_\sigma)]^2}+  \mu_{0i}^2\left[\frac{ \epsilon^2  \sigma^2_i (1+c_0\epsilon_\sigma)}{ \lambda_i\kappa_i^2 + \epsilon^2  \sigma^2_i (1+c_0\epsilon_\sigma)}\right]^2 + \frac{ \sigma^2_i\lambda_i (1+c_0\epsilon_\sigma)}{\epsilon^{-2} \lambda_i\kappa_i^2 +  \sigma^2_i(1+c_0\epsilon_\sigma)}.
\end{eqnarray*}

Denote $i_{\epsilon}(a) = \max\{i: \,  \sigma^2_i/[\lambda_i\kappa_i^2]\leq a\epsilon^{-2}\}$ and recall that we assume $c_0\epsilon_\sigma\leq 1/2$. Then,
similarly to the proof of Theorem~\ref{th:contractiontheoremGeneral},
\begin{eqnarray*}
  S_1 &=&  \sum_{i\leq N} \frac{\sigma^2_i\lambda_i(1+c_0\epsilon_\sigma)}{\epsilon^{-2} \lambda_i\kappa_i^2 + \sigma^2_i(1+c_0\epsilon_\sigma)}
    \asymp   \epsilon^{2} \sum_{i\leq   i_{\epsilon+}} \sigma^2_i \kappa_i^{-2}
 + (1+c_0\epsilon_\sigma)\sum_{  i_{\epsilon+} < i \leq N}  \lambda_i,
\end{eqnarray*}
where $i_{\epsilon+} = i_{\epsilon}(1/(1+c_0\epsilon_\sigma))$, and
\begin{eqnarray*}
 S_2 &=& \sup_{\mu_0 \in S^\beta(A)}\sum_{i\leq N}\mu_{0i}^2\left[\frac{ \epsilon^2 \sigma^2_i(1+c_0\epsilon_\sigma)}{  \kappa_i^2\lambda_i + \epsilon^2 \sigma^2_i(1+c_0\epsilon_\sigma)}  \right]^2
 \asymp  A^2  \left[ \epsilon^4  (1+c_0\epsilon_\sigma)^2 \max_{i\leq \min(N,  i_{\epsilon+})} \left[\frac{\sigma^2_i i^{-\beta}}{  \kappa_i^2\lambda_i }  \right]^2 + C i_\epsilon(2/3)^{-2\beta} I( N >   i_{\epsilon+})\right].
\end{eqnarray*}
The first term is
\begin{eqnarray*}
S_3 &=& \sum_{i\leq N}  \frac{ \epsilon^2 \sigma_i^2 \kappa_i^2\lambda_i^2}{ [\lambda_i\kappa_i^2 + \epsilon^2 \sigma^2_i (1- c_0\epsilon_\sigma)]^2}
\asymp \epsilon^2\sum_{i\leq \min(N,  i_{\epsilon-} )}    \sigma_i^2 \kappa_i^{-2}  + \epsilon^{-2}(1- c_0\epsilon_\sigma)^{-2}\sum_{ i_{\epsilon-} < i \leq N}\sigma_i^{-2} \kappa_i^2\lambda_i^2\\
&\lesssim &  \epsilon^2\sum_{i\leq \min(N,  i_{\epsilon-} )}  \sigma_i^2 \kappa_i^{-2} + (1- c_0\epsilon_\sigma)^{-1}\sum_{ i_{\epsilon-} < i \leq N}\lambda_i
\end{eqnarray*}
where $ i_{\epsilon-} = i_{\epsilon}( 1/(1-c_0\epsilon_\sigma))$ and for $i>i_{\epsilon-}$, $\sigma_i^{-2} \kappa_i^2\lambda_i < (1- c_0\epsilon_\sigma) \epsilon^{2} $. Note that the last term is similar to the second term in $S_1$. We also have the truncation bias term
\begin{eqnarray*}
S_4 &=& \sup_{\mu_0 \in S^\beta(A)} \sum_{i>N} \mu_{0,i}^2  \asymp   N^{-2\beta}.
\end{eqnarray*}
Noting that $i_{\epsilon,-} \geq i_{\epsilon,+}$ and combining the terms proves the theorem.

\hfill $\square$ \newline

\subsection{Error in operator}

\hspace{2pc}{\it Proof of Lemma~\ref{lem:GeneralKwithError}.}

Using the definitions of $\hat{\sigma^2}_j$ and $\tilde\sigma_j^2$, we have
$$
|\hat{\sigma^2}_j - \tilde\sigma_j^2| = |\hat{\kappa}_j^{-2} - {\kappa}_j^{-2}| \leq \sqrt{2} \delta |\xi_j| \hat{\kappa}_j^{-3/2} {\kappa}_j^{-3/2}
$$
and
$$
|\hat{\sigma^2}_j/ \tilde\sigma_j^2 -1| = |\hat{\kappa}_j^{-2}{\kappa}_j^2 -1| \leq \sqrt{2} \delta |\xi_j| \hat{\kappa}_j^{-3/2} {\kappa}_j^{1/2}
$$
Provided that $\delta < C_{n}^{-1} \min_{j\leq N} \kappa_j $, on $\Omega_j = \{|\xi_j| \leq C_n\}$ with probability $\geq 2\Phi(C_n)$,
$$
|\xi_j|/(\kappa_j + \delta \xi_j)^{3/2} \leq 
 \sqrt{2} C_n/(\kappa_j - \delta C_n)^{3/2}.
$$
This holds simultaneously for $j=1,\ldots, N$ on $\Omega = \cup_{j=1}^N \Omega_j$ with $P(\Omega) \geq 1- N(1-2\Phi(C_n))$. For instance, if $C_n = \sqrt{2\log n+ 2 \log N}$, then $P(\Omega) = 1-1/n(1+o(1))$.

If also $\delta \leq  A^{-1} C_{n}^{-1} \min_{j\leq N} \kappa_j$ for some $A>1$, then,  on $\Omega$,
$$
|\hat{\sigma^2}_j/ \tilde\sigma_j^2 -1| \leq  \sqrt{2} (1 - A^{-1})^{-3/2} \delta C_n/\kappa_j  \leq \sqrt{2}A^{-1}/(1 - A^{-1})^{3/2}.
$$
The lemma is proved.
 \hfill $\square$ \newline


\section{Proofs for geometric spectra}

\hspace{2pc}{\it Proof  of Proposition \ref{prop:plugincontractionratesimproved} }

We prove the first part, for absolute bound. The expression for the rate for the relative bound follows easily from Theorem~\ref{th:contractiontheoremGeneralPlugIn2}.

Use Theorem~\ref{th:GeneralPlugin} to derive the rate with $\kappa_i\asymp i^{-p }$, $\sigma_i \asymp i^{\gamma }$, $\lambda_i \asymp \tau_\epsilon^2 i^{-2\alpha +1}$. First we investigate the sets
\begin{eqnarray*}
I_{\epsilon}(a) &=& \{i: \,  \sigma_i^2/[\lambda_i\kappa_i^2]>  a \epsilon^{-2}\} = \{i: \,  i  >  (a \tau_\epsilon^2 \epsilon^{-2})^{1/(2\tilde p +2\alpha+1)}\},\\
I_\sigma(a) &=& \{i: \,  \sigma^2_i < a c_0 \epsilon_\sigma\} = \{i: \,   i >( a c_0 \epsilon_\sigma)^{1/(2\gamma)}\},\\
I_{\sigma \epsilon}(a) &=& \{i: \, c_0 \epsilon_\sigma /[\lambda_i \kappa_i^2] > a \epsilon^{-2} \} = \{i: \, i   > [a \tau_\epsilon^2 \epsilon^{-2} \epsilon_\sigma^{-1}/c_0]^{1/(2 p +2\alpha+1)} \}.
\end{eqnarray*}
Denote
\begin{eqnarray*}
i_\epsilon (a) &=& (a \tau_\epsilon^2 \epsilon^{-2})^{1/(2\tilde p +2\alpha+1)},\\
i_\sigma (a) &=& ( a c_0 \epsilon_\sigma)^{1/(2\gamma)},\\
i_{\sigma \epsilon} (a) &=& [a \tau_\epsilon^2 \epsilon^{-2} \epsilon_\sigma^{-1}/c_0]^{1/(2 p +2\alpha+1)}.
\end{eqnarray*}
Then, the squared contraction rate given by Theorem~\ref{th:GeneralPlugin} is
\begin{eqnarray*}
	{\rate}_{plugin}^2  = \epsilon^{2}   \sum_{i \leq \max(i_\epsilon (2),  i_\sigma(2))} i^{2\tilde p}
 +  \epsilon^2 c_0 \epsilon_\sigma\sum_{i> i_\sigma(2)  \, \&\, i\leq  i_{\sigma\epsilon}(1/3)}  i^{2p} \\
+     \sum_{i > i_{\epsilon}(2/3)} \tau_\epsilon^2 i^{-2\alpha-1}  +  \sum_{i >  \max(i_\sigma(2), i_{\sigma\epsilon}(1/3))}  \tau_\epsilon^2 i^{-2\alpha-1}\\
+    \epsilon^4
 \max\left\{ \max_{i\leq i_{\epsilon}(1)} \left[\frac{ \sigma_i^2 i^{-\beta}}{  \kappa_i^2\lambda_i }  \right]^2, \max_{i\leq i_{\sigma \epsilon}(1)} \left[\frac{ \epsilon_\sigma i^{-\beta}}{  \kappa_i^2\lambda_i }  \right]^2\right\}
 +  \max_{i > \max(i_{\epsilon}(1), i_{\sigma \epsilon}(1))} i^{-2\beta}\\
\asymp \epsilon^{2}  \max(i_\epsilon (2),  i_\sigma(2))^{(2\tilde p+1)_+}
 +  \epsilon^2 c_0 \epsilon_\sigma   (i_{\sigma\epsilon}(1/3)^{2p+1} - i_\sigma(2)^{2p+1} )|I(i_{\sigma\epsilon}(1/3) > i_\sigma(2))  \\
+    \tau_\epsilon^2 [i_{\epsilon}(2/3)]^{-2\alpha}  +  \tau_\epsilon^2 [ \max(i_\sigma(2), i_{\sigma\epsilon}(1/3))]^{-2\alpha}\\
+    \epsilon^4
 \max\left\{   [i_{\epsilon}(1) ]^{2(2\tilde p +2\alpha +1 -\beta)_+}  ,   \epsilon_\sigma^2 [i_{\sigma \epsilon}(1)]^{(2 p +2\alpha +1-\beta)_+}  \right\}
 +  [\max(i_{\epsilon}(1), i_{\sigma \epsilon}(1))]^{-2\beta}\\
  \asymp \epsilon^{2}  \max((2 \tau_\epsilon^2 \epsilon^{-2})^{(2\tilde p+1)_+/(2\tilde p +2\alpha+1)},  (\epsilon_\sigma)^{(2\tilde p+1)_+/(2\gamma)}) \\
 +  \epsilon^2   \epsilon_\sigma  \left( [\tau_\epsilon^2 \epsilon^{-2} \epsilon_\sigma^{-1}]^{(2p+1)/(2 p +2\alpha+1)}  - [( \epsilon_\sigma)^{(2p+1)/(2\gamma)}]  \right)_+  \\
+    \tau_\epsilon^2 (\tau_\epsilon^2 \epsilon^{-2})^{-2\alpha/(2\tilde p +2\alpha+1)} +  \tau_\epsilon^2 [ \max([ \tau_\epsilon^2 \epsilon^{-2} \epsilon_\sigma^{-1}]^{1/(2 p +2\alpha+1)}, (  \epsilon_\sigma)^{1/(2\gamma)})]^{-2\alpha}\\
+    \epsilon^4
 \max\left\{   [(\tau_\epsilon^2 \epsilon^{-2})^{1/(2\tilde p +2\alpha+1)}]^{2(2\tilde p +2\alpha +1 -\beta)_+}  ,   \epsilon_\sigma^2 [[a \tau_\epsilon^2 \epsilon^{-2} \epsilon_\sigma^{-1}/c_0]^{1/(2 p +2\alpha+1)}]^{(2 p +2\alpha +1-\beta)_+}  \right\}\\
 +  [\max( (\tau_\epsilon^2 \epsilon^{-2})^{1/(2\tilde p +2\alpha+1)}, [ \tau_\epsilon^2 \epsilon^{-2} \epsilon_\sigma^{-1}]^{1/(2 p +2\alpha+1)})]^{-2\beta}.
\end{eqnarray*}

First consider the case $i_{\epsilon}(1) \leq i_{\sigma}(1)$, i.e. if
$\epsilon_\sigma \leq (\tau_\epsilon^2 \epsilon^{-2})^{2\gamma/(2\tilde p +2\alpha+1)}$. Below we consider $a=1$ and omit this argument.
 In particular, this implies that $i_\epsilon \leq i_{\sigma \epsilon} \leq i_\sigma$. Then,
 \begin{eqnarray*}
	{\rate}_{plugin}^2
\asymp \epsilon^{2}  \epsilon_\sigma^{(2\tilde p+1)_+/(2\gamma)}
+    \tau_\epsilon^2 (\tau_\epsilon^2 \epsilon^{-2})^{-2\alpha/(2\tilde p +2\alpha+1)}   +  \tau_\epsilon^2 [ \epsilon_\sigma]^{-2\alpha/(2\gamma)}
+  [\tau_\epsilon^2 \epsilon^{-2} \epsilon_\sigma^{-1}]^{-2\beta/(2 p +2\alpha+1)}
\\
+    \epsilon^4 \max\left\{    ( \tau_\epsilon^2 \epsilon^{-2})^{2(2\tilde p +2\alpha +1 -\beta)_+/(2\tilde p +2\alpha+1)} ,   \epsilon_\sigma^2 [\tau_\epsilon^2 \epsilon^{-2} \epsilon_\sigma^{-1}]^{(2 p +2\alpha +1-\beta)_+/(2 p +2\alpha+1)}   \right\}\\
\leq
 \epsilon^{2}   (\tau_\epsilon^2 \epsilon^{-2})^{(2\tilde p+1)_+/(2\tilde p +2\alpha+1)}
+    \tau_\epsilon^2 (\tau_\epsilon^2 \epsilon^{-2})^{-2\alpha/(2\tilde p +2\alpha+1)}
+  (\tau_\epsilon^2 \epsilon^{-2})^{-2\beta/(2\tilde p +2\alpha+1)}\\
+    \epsilon^4    ( \tau_\epsilon^2 \epsilon^{-2})^{2(2\tilde p +2\alpha +1 -\beta)_+/(2\tilde p +2\alpha+1)}\\
\leq    (\tau_\epsilon^2)^{(2\tilde p+1)_+/(2\tilde p +2\alpha+1)} (\epsilon^{2})^{1-(2\tilde p+1)_+/(2\tilde p +2\alpha+1)}
+  (\tau_\epsilon^2 \epsilon^{-2})^{-2\beta/(2\tilde p +2\alpha+1)}\\
+    \epsilon^4    ( \tau_\epsilon^2 \epsilon^{-2})^{2(2\tilde p +2\alpha +1 -\beta)_+/(2\tilde p +2\alpha+1)}.
\end{eqnarray*}
This coincides with the rate of contraction of $\mu$ when $(\sigma_i)$ are known.

Now we consider large $\epsilon_{\sigma}$ i.e. $i_\epsilon(1) \geq i_\sigma(1)$.
%

Observe, $i_\epsilon \geq i_\sigma$ implies
\begin{eqnarray*}
	I_\sigma  \cap  \bar I_{\epsilon} = \{i: i_\sigma \leq i < i_\epsilon \}\\
	\bar I_\sigma  \cap  \bar I_{\epsilon} = \{i: i<i_\sigma,~ i < i_\epsilon \}=\bar I_\sigma \\
			\bar{I}_\sigma  \cap  I_{\epsilon} = \{i: i_\epsilon \leq i < i_\sigma \} = \emptyset \\
\end{eqnarray*}

Thus,
\begin{eqnarray*}
	\epsilon^{-2}  \sum_{i\in \bar I_\sigma  \cap I_{\epsilon}} \lambda_i^2 \kappa_i^2 \sigma_i^{-2}
	= 0
\end{eqnarray*}

\begin{eqnarray*}
	\epsilon^{2}\sum_{i \in \bar I_{\epsilon} \cap \bar I_\sigma} \sigma_i^2 \kappa_i^{-2} = \epsilon^{2}\sum_{i \leq (i_\epsilon \wedge i_\sigma)} i^{2(p+\gamma)} \asymp \epsilon^{2} i_\sigma^{(1+2(p+\gamma))_+} (\log  i_\sigma)^{\mathbb{I}\{1+2(p+\gamma)=0\} }
\end{eqnarray*}

\begin{eqnarray*}
	\epsilon^2 c_0 \epsilon_\sigma\sum_{i\in  I_\sigma  \cap  \bar I_{\epsilon}}    \kappa_i^{-2} = \epsilon^2 c_0 \epsilon_\sigma\sum_{i_\sigma \leq i< i_\epsilon}    i^{2p} = \epsilon^2 c_0 \epsilon_\sigma i_\epsilon^{1+2p}
\end{eqnarray*}

\begin{eqnarray*}
	\sum_{i \in I_{\epsilon} }  \lambda_i =\tau_{\epsilon}^2 \sum \limits_{ i_{\epsilon}\leq i}i^{-1-2\alpha} \asymp \tau_{\epsilon}^2 i_{\epsilon}^{-2\alpha}
\end{eqnarray*}
%
%

Specifically,
\begin{eqnarray*}
 \epsilon^4 \max\left[    \max_{i\in \bar I_{\epsilon}\cap I_{\sigma}} \left[\frac{c_0^2 \epsilon_\sigma^2 i^{-2\beta}}{  \kappa_i^4\lambda_i^2 }  \right], \max_{i\in \bar I_{\epsilon}\cap \bar I_{\sigma}} \left[\frac{\sigma_i^4 i^{-2\beta}}{  \kappa_i^4\lambda_i^2 } \right]     \right]\\
 =  \epsilon^4 \tau_{\epsilon}^{-4}  \max\left[   \epsilon_\sigma^2   \max_{i\in \bar I_{\epsilon}\cap I_{\sigma}} \left[c_0^2 i^{2(1+2\alpha+2p-\beta)} \right], \max_{i\in \bar I_{\sigma}} \left[i^{2(1+2\alpha+2(p+\gamma)-\beta)} \right]     \right],\\
\end{eqnarray*}
where
\begin{eqnarray*}
\epsilon^4 \tau_{\epsilon}^{-4} \epsilon_\sigma^2  \max_{i\in \bar I_{\epsilon}\cap I_{\sigma}} \left[c_0^2 i^{2(1+2\alpha+2p-\beta)} \right] =  \epsilon^4  \tau_{\epsilon}^{-4} \epsilon_\sigma^2  \left[c_0^2 i_\epsilon^{2(1+2\alpha+2p-\beta)} \vee i_\sigma^{2(1+2\alpha+2p-\beta)} \right]\\
  = c_0^2 i_\epsilon^{-2\beta} \vee \epsilon^4  \tau_{\epsilon}^{-4} i_\sigma^{2(1+2\alpha+2(p+\gamma)-\beta)},
\end{eqnarray*}
and
\begin{eqnarray*}
\epsilon^4 \tau_{\epsilon}^{-4} \max_{i\in \bar I_{\sigma}} \left[i^{2(1+2\alpha+2(p+\gamma)-\beta)} \right] = \epsilon^4 \tau_{\epsilon}^{-4} \vee \epsilon^4  \tau_{\epsilon}^{-4} i_\sigma^{2(1+2\alpha+2(p+\gamma)-\beta)}.
\end{eqnarray*}

%

Consequently, we obtain the following rates using Theorem \ref{th:GeneralPlugin}:
\begin{eqnarray*} 
	{\rate}_{plugin}^2  = \epsilon^{2} i_\sigma^{(1+2(p+\gamma))_+} (\log  i_\sigma)^{\mathbb{I}\{1+2(p+\gamma)=0\} }
	+ \epsilon^2 c_0 \epsilon_\sigma i_\epsilon^{1+2p}
	+\tau_{\epsilon}^2 i_{\epsilon}^{-2\alpha}
	+i_\epsilon^{-2\beta}\\
	+\left[\epsilon^4 \tau_{\epsilon}^{-4} \vee \epsilon^4  \tau_{\epsilon}^{-4} i_\sigma^{2(1+2\alpha+2(p+\gamma)-\beta)} \right]
\end{eqnarray*}

Note, using the definition of $i_\epsilon$,
\begin{eqnarray*}
	\epsilon^2 c_0 \epsilon_\sigma i^{2p+1}_\epsilon= \epsilon^2c_0\epsilon_\sigma  i_\epsilon^{2p+1+2\alpha} i_\epsilon^{-2\alpha} = \epsilon^2c_0\epsilon_\sigma(\epsilon^{-2} \tau^2_\epsilon\epsilon_\sigma^{-1})i_\epsilon^{-2\alpha}   = c_0\tau^2_\epsilon i_\epsilon^{-2\alpha}.\\
	\epsilon^2i_{\sigma}^{(2(\gamma+p)+1)} = \epsilon^2c_0\epsilon_\sigma i_\sigma^{2p+1} \leq \epsilon^2c_0\epsilon_\sigma i_\epsilon^{2p+1}
	= c_0\tau^2_\epsilon i_\epsilon^{-2\alpha}.
\end{eqnarray*}

Therefore,
\begin{eqnarray*} 
	{\rate}_{plugin}^2  =  \epsilon^{2} (\log  i_\sigma)^{\mathbb{I}\{1+2(p+\gamma)=0\} }
	+\tau_{\epsilon}^2 i_{\epsilon}^{-2\alpha}
	+i_\epsilon^{-2\beta}
	+\left[\epsilon^4 \tau_{\epsilon}^{-4} \vee \epsilon^4  \tau_{\epsilon}^{-4} i_\sigma^{2(1+2\alpha+2(p+\gamma)-\beta)} \right]\\ 
	= \epsilon^{2} (\log  \epsilon_\sigma^{-1})^{\mathbb{I}\{1+2(p+\gamma)=0\} }\\
	+\tau_{\epsilon}^2 [\epsilon^{-1}_{\sigma}\epsilon^{-2}\tau_\epsilon^2  ]^{-2\alpha/(1+2\alpha+2p)}
	+[\epsilon^{-1}_{\sigma}\epsilon^{-2}\tau_\epsilon^2  ]^{-2\beta/(1+2\alpha+2p)}\\
	+\left[\epsilon^4 \tau_{\epsilon}^{-4} \vee \epsilon^4 \tau_{\epsilon}^{-4} \epsilon_{\sigma}^{\frac{1+2\alpha+2(p+\gamma)-\beta}{\gamma}}\right]
\end{eqnarray*}

For the second part, the terms are essentially the same as in the main theorem, with additional term which behaves as
 a) $-2\alpha -2p-2\gamma>0$: $\tau_{\epsilon} \epsilon^{-2}   N^{ -2\alpha -2p-2\gamma}$,
 b) $-2\alpha -2p-2\gamma=0$:  $\tau_{\epsilon} \epsilon^{-2}   \log (N- i_{\epsilon})$,
 c) $-2\alpha -2p-2\gamma<0$:  $\tau_{\epsilon} \epsilon^{-2} i_{\epsilon}^{ -2\alpha -2p-2\gamma}$.

Putting it together, we have
$$
\tau_{\epsilon} \epsilon^{-2}  \max( N^{ -2\alpha -2p-2\gamma}, i_{\epsilon}^{ -2\alpha -2p-2\gamma}) I(i_{\epsilon}<N) [\log N]^{I (\alpha +p+\gamma=0)}
$$
\hfill $\square$ \newline

\section{Vaguelette-vaguelette } 

\hspace{2pc}{\it Proof of  Corollary~\ref{cor:vaguelettes}. }

Recall that for fractional noise, $\epsilon = n^{-(1-H)}$ and $\sigma_{\windex}^2= C 2^{-j(2H-1)}(1+o(1)$ for large $j$. Prior variance is $\lambda_{\windex} = \tau 2^{- 2  \alpha j  }$.  For the vaguelette-vaguelette approach, assumptions on the scaling coefficients of  Theorems~\ref{th:contractiontheoremGeneral} (and \ref{contractiontheorem2}) hold due to Lemma~\ref{lem:PostRateDepTermsFBM}.

Therefore, applying Theorem~\ref{th:contractiontheoremGeneral} for operator $K$ with $\kappa_{jk} \asymp 2^{-p j}$,
for $\mu_0 \in S^\beta(A)=Q((a_\windex), A)$ in one dimension with $\windex = (j,k)$ and $a_\windex = 2^{j \beta}$, the rate of contraction of the posterior distribution is,
\begin{eqnarray*}
\rate 
  &\asymp&  \epsilon + \left[\epsilon^2\sum_{\windex \in \Windex_0} 2^{j(2p+1-2H)} +  \sup_{\windex \in \Windex_1} 2^{-2 j \beta} +
 \sum_{\windex \in \Windex_1}  \tau 2^{-j (2  \alpha+1)} + \epsilon^4  \max_{\windex \in \Windex_0} \left[2^{j(2p+1-2H-\beta+2\alpha+1)} \right]^2 \right]^{1/2}.
\end{eqnarray*}
 where
 $$
 \Windex_\epsilon = \{\windex \in \Windex_I: \, \sigma^{2}_\windex \beta_\windex^{-2} > \epsilon^{-2} \lambda_\windex \}
 = \{ (j,k): \,   2^{j(2p+1-2H+2\alpha+1)} > n^{2(1-H)} \} =\left\{ (j,k): \,   j >  \frac{(1-H)}{(p+1-H+\alpha)} \log_2 n \right\} = \left\{ (j,k): \,   j >  J_H \right\}
 $$
 where $J_H=\frac{(1-H)}{(p+1-H+\alpha)} \log_2 n$.

Therefore,
\begin{eqnarray*}
\rate 
  &\asymp&  \left[n^{-2(1-H)}\sum_{j \leq  J_H} 2^{j(2p+2-2H)} +  \sup_{j >  J_H} 2^{- 2 j \beta} +
 \sum_{j >  J_H}  \tau 2^{- 2  \alpha j} + \tau^{-2} n^{-4(1-H)} \max_{j \leq  J_H} \left[2^{j(2p+1-2H-\beta+2\alpha+1)} \right]^2 \right]^{1/2}\\
  &\asymp&  \left[n^{-2(1-H)}  2^{J_H(2p+2-2H)} +    2^{- 2 J_H \beta} +
   \tau 2^{- 2  \alpha J_H} + \tau^{-2} n^{-4(1-H)} \left[2^{J_H(2p+2-2H-\beta+2\alpha)_+} \right]^2 \right]^{1/2}\\
  &\asymp&  n^{ -\alpha(1-H)/(p+1-H+\alpha)} +    n^{-   \beta (1-H)/(p+1-H+\alpha)} +
   \tau n^{- \alpha (1-H)/(p+1-H+\alpha)} + \tau^{-1} \left[  n^{-2\min( p+1-H+\alpha,  0.5\beta)(1-H)/(p+1-H+\alpha)} \right].
\end{eqnarray*}
For $\tau = const >0$ independent of $n$, the rate is
\begin{eqnarray*}
\rate &\asymp&    n^{-    \min(\beta, \alpha) (1-H)/(p+1-H+\alpha)} +
   n^{-\min( 2(p+1-H+\alpha),   \beta)(1-H)/(p+1-H+\alpha)}.
\end{eqnarray*}

\hfill $\square$ \newline

\section{Proofs, repeated observations}

This proposition is a general tool for identifying the rate $\epsilon_\sigma$.
\begin{proposition} \label{prop: rep obvs consistencyFull}
Assume probability model (\ref{eq: prob model:rep obvs}), and consider the truncated estimator of $(\sigma_i^2)_{i\geq 1}$ defined by (\ref{eq:repobs:sigmahat}).
 Fix $c_0 \in \mathbb{R}^+$, and define
\begin{equation}\label{def:Msigma}
M_\sigma = M_\sigma(\epsilon_\sigma) = \inf_{M \in \mathbb{N}}\{\sigma_i^2\leq c_0 \epsilon_\sigma, \, \forall i>M\}.
\end{equation}
 Then,  $\forall M \geq M_\sigma$ and  $\epsilon_\sigma$  such that $c_0 \epsilon_\sigma/c_\sigma \leq 1/2$,
\[
P(|\hat \sigma_i^2 -  \sigma_i^2|\leq c_0 \epsilon_\sigma, ~ \forall i \geq 1) \geq 1-  2M  e^{- (m-1) (c_0 \epsilon_\sigma/c_\sigma)^2/6}.
\]
In particular, $P(|\hat \sigma_i^2 -  \sigma_i^2|\leq c_0 \epsilon_\sigma, ~i=1,\dots, M) \rightarrow 1$
if $m\epsilon_\sigma^2 \rightarrow\infty$ and $\frac{\log M}{m\epsilon_\sigma^2} \rightarrow 0$ as $m\to \infty$.
\end{proposition}


\hspace{2pc}{\it Proof of Proposition \ref{prop: rep obvs consistencyFull}.}

We use the following lemma from \citet{laurent_adaptive_2000}.
\begin{lemma}\label{lemma: chi squared r.v. bounds}
	Let $(Y_1,\dots,Y_D)$ be i.i.d Gaussian variables, with mean $0$ and variance $1$. Let $a_1,\dots,a_D$ be non-negative. We set
	\[
	|a|_\infty = \sup_{i=1,\dots,D} |a_i|, \mbox{~~and~~} |a|^2_2 = \sum_{i=1}^{D} a_i^2.
	\]
		Let $Z= \sum_{i=1}^{D} a_i(Y_i^2-1)$. 	
	Then, the following inequalities hold for any positive $x$:
	\begin{eqnarray*}
	P(Z \geq 2|a|_2\sqrt{x} + 2 |a|_\infty x) \leq e^{-x}\\
	P(Z \leq - 2|a|_2\sqrt{x} ) \leq e^{-x}	
	\end{eqnarray*}
	\end{lemma}

Consequently, setting $D=m-1$, and $a_i=1$ for all $i$, implies
\begin{eqnarray*}
&Z= \sum_{i=1}^{m-1} (Y_i^2-1), \mbox{~~and~} \\
&P(Z \geq 2|a|_2\sqrt{x} + 2 |a|_\infty x) = P(Z \geq 2\sqrt{m-1}\sqrt{x} + 2 x) \leq e^{-x}\\
&	P(Z \leq - 2|a|_2\sqrt{x} ) = 	P(Z \leq - 2\sqrt{m-1}\sqrt{x} ) \leq e^{-x}	
\end{eqnarray*}

Furthermore, for some $C>0$,
\begin{eqnarray} 
\nonumber
2\sqrt{m-1}\sqrt{x} + 2 x = C &\iff \sqrt{x} = -\frac{\sqrt{m-1}}{2}+\frac{\sqrt{m-1+2C}}{2}\\ \label{eq:exact constant lb}
&\Longrightarrow x = \frac{(m-1)+C}{2} -\frac{1}{2} \sqrt{(m-1)[m-1+2C]}\\ \nonumber
2\sqrt{m-1}\sqrt{x} = C &\iff \sqrt{x} =\frac{C}{2\sqrt{m-1}}\\ \label{eq:exact constant ub}
&\Longrightarrow x = \frac{C^2}{4(m-1)}
\end{eqnarray}

Hence, for a fixed $i$,
\begin{eqnarray*}
P(|\hat \sigma_i^2 -  \sigma_i^2|\geq  c_0 \epsilon_\sigma) &= P(\frac{m-1}{\sigma_i^2} \hat \sigma_i^2 -  (m-1) \geq  \frac{m-1}{\sigma_i^2} c_0 \epsilon_\sigma)\\
&\qquad + P(\frac{m-1}{\sigma_i^2} \hat \sigma_i^2 -  (m-1) \leq  -\frac{m-1}{\sigma_i^2} c_0 \epsilon_\sigma),\\
&= P(Z \geq \frac{m-1}{\sigma_i^2} c_0 \epsilon_\sigma) + P(Z \leq -\frac{m-1}{\sigma_i^2} c_0 \epsilon_\sigma).
\end{eqnarray*}

Note,
\begin{eqnarray*}
P(|\hat \sigma_i^2 -  \sigma_i^2|\leq c_0 \epsilon_\sigma, ~i=1,\dots, M)
&= 1- P(\exists i \in \{1,\dots, M \}: |\hat \sigma_i^2 -  \sigma_i^2|\geq c_0 \epsilon_\sigma) \\
&\geq 1- \sum_{i=1}^M P(|\hat \sigma_i^2 -  \sigma_i^2|\geq c_0 \epsilon_\sigma)
\end{eqnarray*}
Using Equations (\ref{eq:exact constant lb}) and (\ref{eq:exact constant ub}), Lemma \ref{lemma: chi squared r.v. bounds} implies
\begin{eqnarray*}
\sum_{i=1}^M P(|\hat \sigma_i^2 -  \sigma_i^2|\geq c_0 \epsilon_\sigma)
&\leq \sum_{i=1}^M (e^{-x_{1,i}} + e^{-x_{2,i}}),
\end{eqnarray*}
where, for $i=1,\dots, M$,
\begin{eqnarray} \label{eq:massart tail exponent}
x_{1,i} &= \frac{m-1}{2}(1+\frac{c_0\epsilon_\sigma}{\sigma_i^2}) - \frac{m-1}{2}\sqrt{(1+2\frac{c_0\epsilon_\sigma}{\sigma_i^2})}\\ \label{eq:massart head exponent}
x_{2,i} &= \frac{m-1}{4}(\frac{c_0\epsilon_\sigma}{\sigma_i^2})^2.
\end{eqnarray}

Note that for $x\geq 0$, $1+x- \sqrt{1+2x} \geq 0$ and
\begin{eqnarray}\label{eq:lowerbound}
1+x- \sqrt{1+2x} = \frac{x^2}{1+x+ \sqrt{1+2x}} \geq \frac{x^2}{2(1+x)}.
\end{eqnarray}

Now, we show that ${x}_{1,i} \leq x_{2,i}$:
\begin{eqnarray*}
\frac{m-1}{2} \frac{c_0 \epsilon_\sigma}{\sigma_i^2} +\frac{m-1}{2} - \frac{(m-1)}{2}\sqrt{1+2\frac{c_0 \epsilon_\sigma}{\sigma_i^2}} &\leq& \frac{m-1}{4}(\frac{c_0\epsilon_\sigma}{\sigma_i^2})^2\\
 \iff y^2/2-y-1+\sqrt{1+2y} &\geq& 0
\end{eqnarray*}
Denoting $g(y)= y^2/2-y-1+\sqrt{1+2y}$, we have
\begin{eqnarray*}
g(0)&=0,~~\mbox{and}\\
g'(y) &=y-1+\frac{1}{\sqrt{1+2y}} >0 \iff  y>1 \quad or \quad y\leq 1 \quad and \quad  2 y^2(3/2 - y) >  0
\end{eqnarray*}
which does hold. Hence,  $g(y) \geq 0$, for all $y \geq 0$, and therefore ${x}_{1,i} \leq x_{2,i}$.

Thus,
\begin{eqnarray*}
\sum_{i=1}^M P(|\hat \sigma_i^2 -  \sigma_i^2|\geq c_0 \epsilon_\sigma)
&\leq \sum_{i=1}^M (e^{- {x}_{1,i}} + e^{-x_{2,i}}) \leq \sum_{i=1}^M (2e^{- {x}_{1,i}}) \leq 2M e^{-\min \limits_{i=1, \dots M} {x}_{1,i}}\\
& \leq 2M  e^{- (m-1) (c_0 \epsilon_\sigma/C_2)^2/6}
\end{eqnarray*}
since
$$
\min \limits_{i=1, \dots M}  x_{1,i} =  x_{1,1} = \frac{m-1}{2} \left[\frac{c_0 \epsilon_\sigma}{C_2} +1 - \sqrt{1+2\frac{c_0 \epsilon_\sigma}{C_2}}\right]\geq \frac{m-1}{4} \frac{(c_0 \epsilon_\sigma/C_2)^2}{1+c_0 \epsilon_\sigma/C_2} \geq \frac{m-1}{6} (c_0 \epsilon_\sigma/C_2)^2
$$
using inequality (\ref{eq:lowerbound}) and small $\epsilon_\sigma$ (eg such that $c_0 \epsilon_\sigma/C_2 \leq 1/2$).

For $i >M$, to have
\begin{eqnarray*}
P(|\hat \sigma_i^2 -  \sigma_i^2|\leq c_0 \epsilon_\sigma, ~i>M) = P(\sigma_i^2\leq c_0 \epsilon_\sigma, ~i>M)
= 1,
\end{eqnarray*}
we need $M > M_\sigma = \inf_{M}\{\sigma_i^2\leq c_0 \epsilon_\sigma, \, \forall i>M\}$.

Note, the parameters of interest are $m, \epsilon_\sigma$ and $M$. Thus,
\[
e^{-\frac{(m-1)}{6}(c_0 \epsilon_\sigma/C_2)^2 + \log M} \rightarrow 0 \iff 	\frac{(m-1)}{6}(c_0 \epsilon_\sigma/C_2)^2 - \log M \rightarrow \infty
\]

However, since
\[
\frac{(m-1)}{6} (c_0 \epsilon_\sigma/C_2)^2 - \log M \geq Cm\epsilon_\sigma^2 - \log M
\]
for some $C>0$, it suffices to show
\[
Cm\epsilon_\sigma^2 - \log M = Cm\epsilon_\sigma^2 (1- \frac{\log M}{Cm\epsilon_\sigma^2})  \rightarrow \infty \iff m\epsilon_\sigma^2 \rightarrow\infty,~~ \mbox{and}~~ \frac{\log M}{m\epsilon_\sigma^2} \rightarrow 0.
\]

\hfill $\square$ \newline

\hspace{2pc}{\it Proof of Proposition \ref{prop: rep obvs M contraction2}.}
The relative bound holds with $|\hat\sigma_i^2/\sigma_i^2 -1 |\leq 2x/m+\sqrt{  x/m}$ with probability $\geq 1 - 2e^{-x}$ \cite{laurent_adaptive_2000}. Hence, simultaneously for $i=1,\ldots, N$, this upper bound holds with probability at least $2 N e^{-x}$. Taking $x = 2\log N$ implies that with probability at least $1 - 2/N$, $|\hat\sigma_i^2/\sigma_i^2 -1 |\leq 4 m^{-1}\log N +\sqrt{ 2 m^{-1} \log N}$ for all $i=1\ldots, N$, and this upper bound tends to 0 if $\log N =o(m)$. According to Theorem~\ref{th:contractiontheoremGeneralPlugIn2}, the rate of contraction of the posterior distribution is not affected if $ N \geq i_{\rate,+}$. This proves the proposition.
 \hfill $\square$ \newline


\hspace{2pc}{\it Proof of Lemma~\ref{lem:EstEigenFunctions}.}

We use notation of \citet{KoltchinskiiLounici}. Our aim is to verify their Theorem 3 for $r=1,\ldots, N$ with $t = t_0 \log N$ to obtain a simultaneous upper bound on $||P_r- \hat P_r||^2$ for $r\leq N$  with probability at least $1-e^{-t_0}$ in terms of characteristics of $\Sigma=V$, and to find conditions when this bound tends to 0. Recall that $\Sigma = \sum_{i=1}^{\infty} \sigma^2_i \phi_i \phi_i^T$ and $P_r = \phi_r \phi_r^T$.

For $\sigma_{i} \asymp i^{\gamma}$ with $\gamma < -1/2$,
$$
||\Sigma||_{\infty} =\sup_{i=1, 2, \ldots, } \sigma_i^2 = \sigma_1^2 \in (0,\infty)
$$
and hence
 $$
 r(\Sigma) =  \frac{\trace{\Sigma}}{||\Sigma||_{\infty}} \asymp \sum_{i=1}^{\infty} i^{2\gamma} \asymp 1.
 $$
Hence, $r(\Sigma)=o(m)$ for $\gamma < -1/2$ for any $N$. Note that for $\gamma \in [-1/2, 0]$, $\trace{\Sigma}=\infty$.

The $r$th spectral gap is
$$
g_r = \sigma^2_{r} - \sigma^2_{r+1} \asymp  r^{2\gamma}[ 1-(1+1/r)^{2\gamma}] \gtrsim  r^{2\gamma}[1-2^{2\gamma}] \gtrsim  r^{2\gamma}, \quad r=1,\ldots, N,
$$
and $\bar g_r = \min(g_r, g_{r-1}) = g_r$ for $r\geq 2$ and  $\bar g_1 = g_1$, hence $\bar g_r = g_r$. Here $m_r$, the multiplicity of $\sigma_r$, is 1.

Using Theorem 3 of \citet{KoltchinskiiLounici},  for $\gamma < - 1/2$ we obtain
\begin{eqnarray*}
E || P_r - \hat P_r||_2^2 &\lesssim&  \frac{m_r ||\Sigma||_{\infty}^4}{\bar g_r^4} \max([r(\Sigma)/m]^2,  [r(\Sigma)/m]^4) + r(\Sigma)/m\\
&\lesssim&  r^{-8\gamma} m^{-2} + 1/m.
\end{eqnarray*}

Using this bound and Markov's inequality, we can derive
 $$
 P ( || P_r - \hat P_r||_2^2 > \varepsilon ) \leq \varepsilon^{-1} E || P_r - \hat P_r||_2^2
 $$
and hence the probability of such events simultaneously for all $r=1,\ldots, N$ is
\begin{eqnarray*}
P ( || P_r - \hat P_r||_2^2 > \varepsilon ,\quad r=1,\ldots, N) &\leq&
 \sum_{r=1}^N P ( || P_r - \hat P_r||_2^2 > \varepsilon) \leq  \varepsilon^{-1} \sum_{r=1}^N E || P_r - \hat P_r||_2^2\\
 &\leq&  C_{\gamma} \varepsilon^{-1} N^{ 1-8\gamma }  m^{-2}   + N/m.
\end{eqnarray*}


The first term in the upper bound tends to 0 if   $N = o\left(m^{1/(1/2-4\gamma)}\right)$,
and the second term goes to 0 if $N = o\left(m \right)$.

Therefore, it is possible to estimate $P_r$ simultaneously for $r=1,\ldots, N$ with high probability if $\gamma < -1/2$ and $N = o\left(m^{1/(1/2-4\gamma)}\right)$.
 \hfill $\square$ \newline


\hspace{2pc}{\it Proof of Corollary \ref{cor:frBgeom0}.}

If $ \tau = \tau_n^2 < n^{2(\alpha+p)}$, assuming that $1-2H +2p+2\alpha>0$,  the posterior contraction rate
\begin{gather*}
   n^{2H-2} \sum_{i\leq   i_{\epsilon-}}  i^{2p+1-2H}
 +  \sum_{  i_{\epsilon+} < i \leq n}  \tau_n^2 i^{-1-2\alpha} +  n^{2(1-H)}  \tau_n^4  \sum_{   i_{\epsilon-} < i \leq n} i^{-3+2H -2p-4\alpha}
+     n^{4(H-1)} \tau_n^{-4} \max_{i\leq i_{\epsilon+}}  i^{2(2-2H-\beta+2p  +2\alpha)}
+   i_{\epsilon+}^{-2\beta}\\
\lesssim  n^{-2(1-H)}  [\tau_n^2 n^{2(1-H)}]^{(p+1-H)/(1-H+\alpha+p)}
 +   \tau_n^2   [\tau_n^2 n^{2(1-H)}]^{-\alpha/(1-H+\alpha+p)}
   +    \tau_n^2 [\tau_n^2 n^{2(1-H)}]^{ - \alpha/(1-H+\alpha+p)}\\
+     n^{4(H-1)} \tau_n^{-4} [\tau_n^2 n^{2(1-H)}]^{2(1-H-0.5\beta+p  +\alpha)_+/(1-H+\alpha+p)}
+  [\tau_n^2 n^{2(1-H)}]^{-\beta/(1-H+\alpha+p)}\\
\lesssim   \tau_n^2  [\tau_n^2 n^{2(1-H)}]^{-\alpha/(1-H+\alpha+p)}
+  [\tau_n^2 n^{2(1-H)}]^{-\min( 2,   \beta/(1-H+\alpha+p))}
+  [\tau_n^2 n^{2(1-H)}]^{-\beta/(1-H+\alpha+p)}.
 \end{gather*}
 Hence, for $\beta\leq 2(1-H+\alpha+p)$, the $\tau_n$ minimising this expression is $ \tau_n^2   =  [ n^{2(1-H)}]^{(\alpha-\beta)/(1-H+\beta+p)}$, and the squared  posterior contraction rate is $[ n^{2(1-H)}]^{-\beta/(1-H+\beta+p)}$ which coincides with the minimax rate.

\hfill $\square$ \newline


\section{Proofs, empirical Bayes approach}

\hspace{2pc}{\it Proof of Theorem \ref{th:EBtau1}. }

Recall that we consider an inverse problem with heterogeneous variance that can be written in the sequence space as
\begin{equation}\label{eq:model ass.}
Y_i|\mu_i \sim N(\kappa_i\mu_i, \epsilon^2\sigma_i^2),~~~and~~ \mu_i \sim N(0,\lambda_i), \text{ independently },
\end{equation}
where $\kappa_i\asymp i^{-p}$, $\lambda_i = \tau i^{-1-2\alpha}$, $\sigma_i \asymp i^{\gamma}$, with noise level $\epsilon \to 0$, $\tau= \tau_\epsilon^2$.

Using the substitution $\nu^{1+2\alpha_s} = \epsilon^{-2}\tau$ similarly to \citet{SzaboEmpiricalBayes}, the log likelihood for $\nu$ (e.g. with respect to $\prod_{i=1}^{\infty} N(0, \epsilon^2 \sigma_i^2)$) under the considered model can be written as
\begin{equation*}
\ell(\nu) = -\frac{1}{2}\sum_{i=1}^\infty [\log(1+   \nu^{1+2\alpha_s} i^{-1-2\alpha_s}) + \frac{\epsilon^{-2} Y_i^2 \sigma_i^{-2} i^{1+2\alpha_s}}{ i^{1+2\alpha_s}+  \nu^{1+2\alpha_s} } ],
\end{equation*}
where $\alpha_s=\alpha + p +\gamma$,  which coincides with the likelihood in the direct white noise case considered by \citet{SzaboEmpiricalBayes} with $\alpha_s$ instead of $\alpha$, $\epsilon^{-2}$ instead of $n$ and $X^2_i = Y_i^2/\sigma_i^2$.

In particular, the derivative map studied in  \citet{SzaboEmpiricalBayes} can be written as
\[
\mathbb{M} (\nu) = -\ell'(\nu) = \frac{1+2\alpha_s}{2}\left[
\epsilon^{-2}\sum_{i=1}^\infty \frac{\nu^{2\alpha_s}i^{1+2\alpha_s} }{ [ i^{1+2\alpha_s}+ \nu^{1+2\alpha_s} ]^2 } X_i^2
- \sum_{i=1}^\infty \frac{\nu^{2\alpha_s} }{  i^{1+2\alpha_s}+ \nu^{1+2\alpha_s}  }
  \right]
\]
 where $E \mathbb{M}(\nu) = \frac{1+2\alpha_s}{2}[\epsilon^{-2} h(\nu) - C_{\alpha, \nu}]$ where
\[
h(\nu) = \sum_{i=1}^\infty \frac{ i^{1+2\alpha_s}  \nu^{2\alpha_s} \eta^2_{0,i}  }{ [ i^{1+2\alpha_s}+ \nu^{1+2\alpha_s} ]^2 },
\]
with $\eta^2_{0,i} = \kappa_i^2 \mu_{0,i}^2/  \sigma_i^{2}$ and
$$
C_{\alpha_s, \nu}  = \sum_{i=1}^\infty \frac{\nu^{-1} }{  (i/\nu)^{1+2\alpha_s}+ 1}.
$$
By Lemma A.1 in \citet{SzaboEmpiricalBayes}, $C_{\alpha, \nu} \leq C_{\alpha} = \int_0^{\infty} (1+x^{1+2\alpha})^{-2} dx$ for all $\nu>0$, and if $\nu \to \infty$, $C_{\alpha, \nu} \to C_{\alpha}$.

 If $\mu_0 \in H^{\beta}(A)$ then $\eta_0 = \sum_i \eta_{0,i} e_i \in H^{\beta_s}(A)$ with $\beta_s = \beta+p+\gamma$. Hence, as long as $\beta_s=\beta+p+\gamma >0 $ and $\alpha_s=\alpha + p +\gamma>0$,  the statement and the proof  of Theorem 2.1  \citet{SzaboEmpiricalBayes} applies with $\alpha_s$ and $\beta_s$ instead of $\alpha$ and $\beta$.

In the considered paper, we study $\mu_0 \in S^{\beta}(A)$, so we adapt the proof of   Theorem 2.1 of \citet{SzaboEmpiricalBayes} to this case. The change of the functional space only affects the expression of $h(\nu)$ which is stated in Proposition~\ref{prop: bounds for h(nu) and overline(nu)}.

We need to modify the proof of Theorem 2.2 slightly for the considered inverse problem with heterogeneous noise and for $\mu_0 \in S^{\beta}(A)$ which we do in Proposition~\ref{prop: posterior risk: hetero}.



For positive constants $l < L$, define
\begin{eqnarray*}
	\overline{\nu} = \sup \{ \nu>0 :\epsilon^{-2}h(\nu) \geq l \}, \\
	\underline{\nu} = \sup \{ \nu>0 :\epsilon^{-2}h(\nu) \geq L \} .
\end{eqnarray*}
This implies that for $\nu \geq \overline{\nu}$, $h(\nu) \leq l \epsilon^2$, and $\nu \geq \underline{\nu}$, $h(\nu) \leq L \epsilon^2$.

\begin{proposition}\label{prop: bounds for h(nu) and overline(nu)}
	Let the assumptions  \eqref{eq:model ass.} hold, and assume $1+2\alpha_s>0$. 
Then, for $\mu_0 \in S^{\beta}$,  and $\nu \geq 1$,
\begin{align*}
h(\nu) \lesssim ||\mu_0||^2_{S^\beta}
\begin{cases}
\nu^{-(1+ 2\beta_s)}	, &\text{ for } \alpha_s +1/2\geq \beta_s,\\
\nu^{-2(1+\alpha_s)}	, &\text{ for }\alpha_s+1/2 < \beta_s,
\end{cases}
\end{align*}	
which imply $\underline{\nu} \gtrsim (\frac{\epsilon^{-2} ||\mu_0||^2_{S^\beta}}{L} )^{\frac{1}{2(1+\alpha_s)}}$ and
\begin{align}\label{eq:UpprBoundOverlineNu}
\overline{\nu} \lesssim
\begin{cases}
(\frac{\epsilon^{-2} ||\mu_0||^2_{S^\beta}}{l} )^{\frac{1}{1+ 2 \beta_s }}	,  &\text{ for } \alpha_s +1/2\geq \beta_s,\\
(\frac{\epsilon^{-2} ||\mu_0||^2_{S^\beta}}{l} )^{\frac{1}{2(1+\alpha_s)}}	,&\text{ for } \alpha_s +1/2 < \beta_s,
\end{cases}
\end{align}

If  also  Assumption~\ref{Assume:SelfSim} holds then $h(\nu) \gtrsim c_0 \nu^{-1-2\beta_s}$ and hence $\underline{\nu}  \gtrsim [L \epsilon^2]^{-1/(1+2\beta_s)}$.
\end{proposition}
Note that for $\alpha_s +1/2\geq \beta_s$,
$$
[\epsilon^2]^{-1/(1+2\beta_s)} \lesssim  \underline{\nu} \leq \overline{\nu} \lesssim [\epsilon^2]^{-1/(1+2\beta_s)},
$$
i.e. $\underline{\nu}$ and $ \overline{\nu}$ are of the same order.

Using Proposition~\ref{prop: bounds for h(nu) and overline(nu)} instead of the corresponding upper bounds on $h(\nu)$ and on $\overline{\nu}$ in the proof of Theorem 2.1 in \citet{SzaboEmpiricalBayes}, we obtain that $P_{\mu_0}(\underline{\nu} \leq \nu \leq \overline{\nu}) \to 1$ as $\epsilon \to 0$.

Now we show that $\hat\nu$ converges in probability.
\begin{proposition}\label{prop:convEBtau}
	Let the assumptions  \eqref{eq:model ass.} hold, and assume $\alpha + p +\gamma+1/2>0$. 
Then, for $\mu_0 \in S^{\beta}$ satisfying Assumption~\ref{Assume:SelfSim}  and $\alpha +0.5 \geq \beta$,  and $\epsilon$ small enough,
$$
\hat\nu/\nu^\star \to 1 \text{ in probability}
$$
where $\nu^\star = \arg \max_{\nu} E \ell(\nu) \asymp [\epsilon^2]^{-1/(1+2\beta_s)}$.
\end{proposition}

Now we study the posterior risk.

\begin{proposition}\label{prop: posterior risk: hetero}
Consider model \eqref{eq:model ass.}, with  parametrisation $\nu^{1+2\alpha_s} = \epsilon^{-2}\tau$ and $1+2\alpha_s>0$ where $\alpha_s = \alpha+p+\gamma$. Assume $\mu_0 \in S^\beta(A)$.

 If $\mu_0 \neq 0$, then,  for any $\nu \in [\underline{\nu}, \overline{\nu}]$,
	\begin{eqnarray}\label{eq:contraction rates: parametrized}
	\mathbb{E}_{\mu_0} \mathbb{E} \sum_{i=1}^{\infty}\left[ (\mu_i-\mu_{0i})^2\mid Y\right]
	& \lesssim& 	 [\epsilon^2]^{1 - 2(0.5+p+\gamma)_+/[2\beta + 1+2(p+\gamma)]}  [\ln(1/\epsilon)]^{I(0.5+p+\gamma=0)}\\
&&+ [\epsilon^{2}]^{2\min(\beta, 1+2\alpha_s)/(2\alpha_s+2)}  ||\mu_0||_{S^\beta}^2.
 \notag
	\end{eqnarray}
As $\epsilon \to 0$, the upper bound tends to 0, hence the posterior distribution is consistent.


In particular, if $ \alpha +1/2\geq \beta$ and $ \alpha \geq  \beta/2 - 1/2 - (p+\gamma)$, then
	\begin{align}\label{eq:contraction rates1}
	\mathbb{E}_{\mu_0}  \sum_{i=1}^{\infty}\mathbb{E}\left[ (\mu_i-\mu_{0i})^2\mid Y, \hat \nu\right]
	& \lesssim& [\epsilon^{2}]^{2\beta/(2\beta+(1+2\gamma+2p)_+)} + \epsilon^2 [\ln(1/\epsilon)]^{I(0.5+p+\gamma=0)},
	\end{align}
i.e. the posterior distribution contracts at the optimal rate, in a minimax sense.

If $\mu_0 =0$, then $\mathbb{E}_{\mu_0}  \sum_{i=1}^{\infty}\mathbb{E}\left[ (\mu_i-\mu_{0i})^2\mid Y, \hat \nu\right] \lesssim \epsilon^2$.
\end{proposition}

\hspace{2pc}{\it Proof of Proposition \ref{prop: bounds for h(nu) and overline(nu)}.}

	Denote $i_\nu :=\max \{i: i \leq \nu \} $. Then, for $1+2\alpha_s>0$ (or equivalently $\gamma > -1/2 -(p+\alpha)$),
\begin{align*}
i^{1+2\alpha_s}+ \nu^{1+2\alpha_s} \asymp
	\begin{cases}
	\nu^{1+2\alpha_s}, &\text{ for } i \leq i_\nu, \\
	i^{1+2\alpha_s}, &\text{ for } i > i_\nu.
	\end{cases}
\end{align*}
	
	Hence,
	\begin{align*}
h(\nu)
&=\sum_{i=1}^\infty \frac{ i^{1+2\alpha_s}  \nu^{2\alpha_s}  }{ [ i^{1+2\alpha_s}+ \nu^{1+2\alpha_s} ]^2 }  \eta^2_{0,i}
 \asymp \sum_{i \leq i_\nu} \frac{ i^{1+2\alpha_s}  \nu^{2\alpha_s}  }{ [ \nu^{1+2\alpha_s} ]^2 }  \eta^2_{0,i}
+
\sum_{i > i_\nu} \frac{ i^{1+2\alpha_s}  \nu^{2\alpha_s}  }{ [ i^{1+2\alpha_s} ]^2 }  \eta^2_{0,i}\\
& \asymp  \nu^{-2-2\alpha_s} \sum_{i \leq i_\nu}  i^{1+2\alpha_s}  i^{-2p-2\gamma} \mu^2_{0,i}
+
 \nu^{2\alpha_s}  \sum_{i > i_\nu}  i^{-1-2\alpha_s}   i^{-2p-2\gamma} \mu^2_{0,i}.
	\end{align*}
Using $\mu_0 \in S^\beta(A)$, we have
\begin{align*}
h(\nu)
& \lesssim   \nu^{-2-2\alpha_s} \sum_{i \leq i_\nu}  i^{1+2(\alpha_s- \beta_s)}   \mu^2_{0,i} i^{2\beta}
+
\nu^{2\alpha_s} \sum_{i > i_\nu}  i^{-1-2(\alpha_s + \beta_s)}  \mu^2_{0,i} i^{2\beta}\\
& \lesssim \nu^{-2-2\alpha_s + 2(\alpha_s+0.5- \beta_s)_+}  ||\mu_0||^2_{S^\beta}+
\nu^{  -(1+2 \beta_s )} ||\mu_0||^2_{S^\beta}\\
& \lesssim \nu^{-2-2\alpha_s + 2(\alpha+0.5- \beta)_+}  ||\mu_0||^2_{S^\beta},
\end{align*}
where $i_\nu \asymp \nu$ for $\nu \geq 1$. There exists  $\mu_0 \in S^\beta(A)$  such that the equality holds, up to a constant, similarly to the proof of Theorem 1.

Note, when   $\alpha_s +0.5 \leq \beta_s$ and $\nu \geq 1$, $h(\nu) \lesssim \nu^{-2-2\alpha_s}  ||\mu_0||^2_{S^\beta}$.
When $\alpha_s+0.5 = \beta_s $, both terms are of the same order. Hence, for  $\nu \geq 1$,
\begin{align*}
h(\nu) \lesssim ||\mu_0||^2_{S^\beta}
\begin{cases}
\nu^{-(1+  2\beta_s )}	, &\text{ for } \alpha_s +1/2 \geq   \beta_s,\\
\nu^{-2(1+\alpha_s)}	, &\text{ for } \alpha_s +1/2 <  \beta_s.
\end{cases}
\end{align*}	

Consequently, for $\alpha_s +1/2 \geq   \beta_s$ and  $\nu \geq 1$,
\[
l \leq \epsilon^{-2} h(\nu) \lesssim ||\mu_0||^2_{S^\beta}\epsilon^{-2} \nu^{-(1+  2\beta_s )} \iff \nu \lesssim (\frac{\epsilon^{-2} ||\mu_0||^2_{S^\beta}}{l} )^{\frac{1}{(1+ 2\beta_s)}}.
\]

Doing the same for $\alpha_s +1/2 < \beta_s$ implies,
\begin{align*}
\overline{\nu} \lesssim
\begin{cases}
(\frac{\epsilon^{-2} }{l} )^{\frac{1}{1+ 2 \beta_s }}	,  &\text{ for } \alpha_s +0.5\geq  \beta_s,\\
(\frac{\epsilon^{-2} }{l} )^{\frac{1}{2(1+\alpha_s)}}	,&\text{ for }\alpha_s +0.5 <  \beta_s.
\end{cases}
\end{align*}

  If also Assumption~\ref{Assume:SelfSim} holds, then for $\nu \geq 1$ and $1+2\alpha_s-2\beta_s \geq 0$,
\begin{align*}
h(\nu)& \gtrsim c_0 \nu^{-2-2\alpha_s}    i_\nu^{1+2\alpha_s-2\beta_s}  \asymp   \nu^{-1-2\beta_s}.
\end{align*}

If $1+2\alpha_s < 2\beta_s$ then
	\begin{align*}
h(\nu)&\asymp \nu^{-2-2\alpha_s} \sum_{i \leq i_\nu}  i^{1+2\alpha_s}    \eta^2_{0,i} = \nu^{-2-2\alpha_s}  \sum_{i \leq i_\nu}  i^{1+2\alpha_s-2\beta_s} i^{2\beta_s} \eta^2_{0,i} \\
&\geq \nu^{-2-2\alpha_s}  i_\nu^{1+2\alpha_s-2\beta_s} c_0 \geq  c_0 \nu^{-1-2\beta_s}.
\end{align*}
Since for $\nu \geq \underline{\nu}$, $h(\nu) \leq L \epsilon^2$ and $h(\nu) \asymp \nu^{-1-2\beta_s}$, we have
$$
L \epsilon^2 \geq h(\underline{\nu}) \asymp \underline{\nu}^{-1-2\beta_s}
$$
which implies $\underline{\nu}  \gtrsim [L \epsilon^2]^{-1/(1+2\beta_s)}$.
\hfill $\square$ \newline


\hspace{2pc}{\it Proof of Proposition~\ref{prop:convEBtau}.}

For $\mu_0 \in S^{\beta}$ satisfying Assumption~\ref{Assume:SelfSim}, by Proposition~\ref{prop: bounds for h(nu) and overline(nu)}, with probability tending to 1,
$$
[L \epsilon^2]^{-1/(1+2\beta_s)} \lesssim \underline{\nu} \leq \hat\nu \leq \overline{\nu} \lesssim ( \epsilon^{-2}  l )^{-1/(1+ 2 \beta_s )},
$$
hence, applying the argument in the proof of the similar result in Section 4.4 of \citet{SzaboEmpiricalBayes}, we have $\sup_{\nu \geq \underline{\nu}} |M(\nu) - E M(\nu)| \to 0$ in probability as $\epsilon \to 0$ which implies that  $|E \ell'(\nu)|_{\nu = \hat\nu}| \to 0$ in probability. For $\nu \in [\underline{\nu}, \overline{\nu}]$ and $\epsilon \to 0$,
$$
E M(\nu)=- E \ell'(\nu) = (1/2+\alpha_s)[\epsilon^{-2} h(\nu) - c_{\alpha_s}(1+o(1))]
$$
where $c_{\alpha_s} = \int_0^{\infty} (1_x^{2\alpha_s+1}+1)^{-2} dx$.
 Since $h(\nu)$ is a continuous monotonically decreasing function of $\nu$, there exists a unique solution to $E \ell'(\nu) =0$ denoted by $\nu^\star$ and hence $\hat\nu \to \nu^\star$ in probability. Since $\nu^\star = \arg \max_{\nu} E \ell(\nu)$ solves $h(\nu) - c_{\alpha_s}(1+o(1))$,  we have $ \underline{\nu} \leq \nu^\star \leq \overline{\nu}$ for sufficiently small $\ell$ and large $L$, hence  $\nu^\star \asymp [\epsilon^2]^{-1/(1+2\beta_s)}$.
 The proposition is proved.
\hfill $\square$ \newline



\hspace{2pc}{\it Proof of Proposition \ref{prop: posterior risk: hetero}.}
	
For $\mu_0\neq 0$, from Proposition~\ref{prop: bounds for h(nu) and overline(nu)}, $\underline{\nu} \to \infty$ as $\epsilon \to 0$.
	It is easy to see that with $\nu^{1+2\alpha_s} = \epsilon^{-2}\tau$,
	\begin{align*}
	\mu_i|Y_i  &\sim N(\frac{\nu^{1+2\alpha_s}}{\nu^{1+2\alpha_s}+i^{1+2\alpha_s}}\frac{Y_i}{\kappa_i},\frac{\nu^{1+2\alpha_s} }{\nu^{1+2\alpha_s}+i^{1+2\alpha_s}} \epsilon^2 \kappa_i^{-2}\sigma_i^2).
	\end{align*}
	
	Consequently, similarly to the proof of Theorem 1, 
 using the parametrisation $\nu^{1+2\alpha_s} = \epsilon^{-2}\tau$, we have
	\begin{align*}
	\mathbb{E}_{\mu_0}  \sum_{i=1}^{\infty}\mathbb{E}\left[ (\mu_i-\mu_{0i})^2\mid Y \right]
	&=\sum_i  \frac{ \epsilon^{-2} \sigma_i^2 \kappa_i^2\lambda_i^2}{ [\epsilon^{-2}\lambda_i\kappa_i^2 +  \sigma^2_i]^2}+   \frac{  \sigma^4_i \mu_{0i}^2 }{[\epsilon^{-2} \lambda_i\kappa_i^2 +  \sigma^2_i ]^2}  + \frac{\sigma^2_i\lambda_i}{\epsilon^{-2} \lambda_i\kappa_i^2 + \sigma^2_i}\\
 &\asymp \sum_i
	\frac{\nu^{2+4\alpha_s} }{[\nu^{1+2\alpha_s}+i^{1+2\alpha_s}]^2} \epsilon^2 i^{2(\alpha_s-\alpha)}
	+ \frac{i^{2+4\alpha_s} \mu_{0i}^2 }{[\nu^{1+2\alpha_s}+i^{1+2\alpha_s}]^2} \\
	& \qquad+ \frac{\nu^{1+2\alpha_s} }{\nu^{1+2\alpha_s}+i^{1+2\alpha_s}} \epsilon^2 i^{2(\alpha_s-\alpha)}.
\end{align*}

The series above that are independent of $\mu_0$ can be bounded as follows:
	\begin{align*}
	\sum_{i=1}^\infty \frac{\nu^{1+2\alpha_s} }{\nu^{1+2\alpha_s}+i^{1+2\alpha_s}} \epsilon^2 i^{2(\alpha_s-\alpha)}
	& \asymp \epsilon^2  \sum_{i \leq i_\nu}  i^{-1+2(0.5+p+\gamma)}
	+
	\epsilon^2  \nu^{1+2\alpha_s}   \sum_{i > i_\nu}  i^{-1 -2\alpha} \\
	& \asymp \epsilon^2  i_\nu^{2(0.5+p+\gamma)_+} [ \ln(i_\nu)]^{I(0.5+p+\gamma=0)}  +
	\epsilon^2  \nu^{1+2\alpha_s} i_\nu^{-2\alpha}\\
	& \asymp \epsilon^2   \nu^{2(0.5+p+\gamma)_+} [ \ln(i_\nu)]^{I(0.5+p+\gamma=0)},
	\end{align*}
where $i_\nu :=\max \{i: i \leq \nu \} $, $\nu \geq 1$ and $1+2\alpha_s>0$.
	
	Similarly, for $\nu \geq 1$,
	\begin{align*}
	\sum_{i=1}^\infty \frac{\nu^{2+4\alpha_s} }{[\nu^{1+2\alpha_s}+i^{1+2\alpha_s}]^2} \epsilon^2 i^{2(\alpha_s-\alpha)}
	& \asymp \epsilon^2  \sum_{i \leq i_\nu}  i^{2(\alpha_s-\alpha)}
	+
	\epsilon^2 \sum_{i > i_\nu}  i^{-2(1+\alpha_s+\alpha)}  \nu^{2+4\alpha_s}  \\
	& \asymp \epsilon^2  [\nu^{2(\frac{1}{2}+\alpha_s-\alpha)_+} \vee \ln(\nu)I_{2(\frac{1}{2}+\alpha_s-\alpha)=0}   ]\\
	& \qquad +
	\epsilon^2  \nu^{2(\frac{1}{2}+\alpha_s-\alpha)}.
	\end{align*}

	Thus, for $\nu \in [\underline{\nu}, \overline{\nu}]$, these two terms are bounded by
	\begin{eqnarray}\label{eq:EBtauSum2}
	&&\sum_i \left[\frac{\nu^{2+4\alpha_s} }{[\nu^{1+2\alpha_s}+i^{1+2\alpha_s}]^2} \epsilon^2 \kappa_i^{-2}\sigma_i^2
	+ \frac{\nu^{1+2\alpha_s} }{\nu^{1+2\alpha_s}+i^{1+2\alpha_s}} \epsilon^2 \kappa_i^{-2}\sigma_i^2\right]\notag\\
	&\lesssim& \epsilon^2  [\overline{\nu}^{2(\frac{1}{2}+p+\gamma)_+} [\ln(\overline{\nu})]^{I(0.5+p+\gamma=0)}\notag\\
	& \lesssim&  [\epsilon^2]^{1 - 2(0.5+p+\gamma)_+/[2\min( \beta,0.5+\alpha) + 2(0.5+p+\gamma)]}  [\ln(1/\epsilon)]^{I(0.5+p+\gamma=0)}
	\end{eqnarray}
using the upper bound on $\overline{\nu}$ given by (\ref{eq:UpprBoundOverlineNu}).


	The final term
	\[
	\sum_i \frac{i^{2+4\alpha_s} \mu_{0i}^2 }{[\nu^{1+2\alpha_s}+i^{1+2\alpha_s}]^2} 
	\]
	will be bounded by splitting the summation indices into groups $i\leq \underline{\nu}$ and $i>\underline{\nu} $.

For  $\nu \geq \underline{\nu}\geq 1$ and $\mu_0 \in S^{\beta}(A)$,
\begin{eqnarray*}
	\sum_{i>\underline{\nu}} \frac{i^{2+4\alpha_s} \mu_{0,i}^2 }{[\nu^{1+2\alpha_s}+i^{1+2\alpha_s}]^2} &\lesssim&
	\sum_{i>\underline{\nu}} \frac{i^{2+4\alpha_s} \mu_{0,i}^2 }{[\underline{\nu}^{1+2\alpha_s}+i^{1+2\alpha_s}]^2} \\
&\lesssim& 	\sum_{i>\underline{\nu}}   \mu_{0,i}^2 \lesssim  \underline{\nu}^{-2\beta} ||\mu_0||_{S^\beta}^2.
\end{eqnarray*}
Also,
\begin{eqnarray*}
	\sum_{1 \leq  i\leq \underline{\nu}} \frac{i^{2+4\alpha_s} \mu_{0,i}^2 }{[\nu^{1+2\alpha_s}+i^{1+2\alpha_s}]^2} &\lesssim&
	\sum_{1 \leq  i\leq \underline{\nu}} \frac{i^{2+4\alpha_s} \mu_{0,i}^2 }{[\underline{\nu}^{1+2\alpha_s}+i^{1+2\alpha_s}]^2} \\
&\lesssim& \underline{\nu}^{-2(1+2\alpha_s)}\sum_{1 \leq  i\leq \underline{\nu}} i^{2+4\alpha_s -2\beta} i^{2\beta}\mu_{0,i}^2  \\
&\lesssim& \underline{\nu}^{-2(1+2\alpha_s) +(2+4\alpha_s -2\beta)_+}  ||\mu_0||_{S^\beta}^2.
\end{eqnarray*}
Hence,
\begin{eqnarray*}
	\sum_{i=1}^{\infty} \frac{i^{2+4\alpha_s} \mu_{0,i}^2 }{[\nu^{1+2\alpha_s}+i^{1+2\alpha_s}]^2}
&\lesssim& 	\left[ \underline{\nu}^{-2\beta} + \underline{\nu}^{-2(1+2\alpha_s) +(2+4\alpha_s -2\beta)_+} \right] ||\mu_0||_{S^\beta}^2\\
&\lesssim& 	\left[ [\epsilon^{2}]^{2\beta/(2\alpha_s+2)} + [\epsilon^{2}]^{(1+2\alpha_s)/(1+\alpha_s)}I(1+2\alpha_s -\beta<0) \right] ||\mu_0||_{S^\beta}^2\\
&\lesssim& 	  [\epsilon^{2}]^{2\min(\beta, 1+2\alpha_s)/(2\alpha_s+2)}  ||\mu_0||_{S^\beta}^2
\end{eqnarray*}
using the lower bound on $\underline{\nu}$ in Proposition~\ref{prop: bounds for h(nu) and overline(nu)}.

Hence, for any $\nu \in [\underline{\nu}, \overline{\nu}]$,
\begin{eqnarray*}
	\mathbb{E}_{\mu_0}  \sum_{i=1}^{\infty}\mathbb{E}\left[ (\mu_i-\mu_{0i})^2\mid Y \right]
	& \lesssim&  [\epsilon^2]^{1 - 2(0.5+p+\gamma)_+/[2\min( \beta,0.5+\alpha) + 2(0.5+p+\gamma)]}  [\ln(1/\epsilon)]^{I(0.5+p+\gamma=0)}\\
&&+ [\epsilon^{2}]^{2\min(\beta, 1+2\alpha_s)/(2\alpha_s+2)}  ||\mu_0||_{S^\beta}^2.
\end{eqnarray*}
If $\beta\leq \alpha+0.5$, then
\begin{eqnarray*}
	\mathbb{E}_{\mu_0}  \sum_{i=1}^{\infty}\mathbb{E}\left[ (\mu_i-\mu_{0i})^2\mid Y \right]
	& \lesssim&  [\epsilon^2]^{1 - 2(0.5+p+\gamma)_+/[2\beta + 1+2(p+\gamma)]}  [\ln(1/\epsilon)]^{I(0.5+p+\gamma=0)}\\
&&+ [\epsilon^{2}]^{2\min(\beta, 1+2\alpha_s)/(2\alpha_s+2)}  ||\mu_0||_{S^\beta}^2.
\end{eqnarray*}
If $1+2\alpha_s <\beta $ , the last term has suboptimal rate. Hence the optimal rate is achieved under  additional condition $1+2\alpha_s \geq \beta $.
	
If $\mu_0 =0$, then $\hat \nu =O_P(1)$ and the final sum is 0, and hence, using (\ref{eq:EBtauSum2}),
$$
\mathbb{E}_{\mu_0}  \sum_{i=1}^{\infty}\mathbb{E}\left[ (\mu_i-\mu_{0i})^2\mid Y, \hat \nu\right] \lesssim \epsilon^2.
$$
\hfill $\square$ \newline

\end{document}